\documentclass{article}

\usepackage{graphicx}
\RequirePackage{amsmath,amsfonts,amssymb,mathtools}
\newcommand{\defeq}{\vcentcolon=}
\usepackage{hyperref}
\usepackage{physics}
\usepackage{color}
\RequirePackage{bm}
\RequirePackage{siunitx}

\graphicspath{{./fig/}{./}}

\newcommand{\Var}{\mathbb{V}\hspace{-.2mm}{\rm ar}}
\newcommand{\Esp}{\mathbb{E}}
\newcommand{\Cov}{\mathbb{C}\hspace{-.2mm}{\rm ov}}

\begin{document}

\title{Uncertainty quantification for mineral precipitation and dissolution in fractured porous media}

\author{%
  Michele Botti\and
        Alessio Fumagalli        \and
        Anna Scotti
}
\date{\small Dipartimento di Matematica, Politecnico di Milano, piazza Leonardo da Vinci 32, 20133 Milano, Italy}

%
%
%

\maketitle

\begin{abstract}
    In this work, we present an uncertainty quantification analysis to determine
    the influence and importance of some physical parameters in a reactive
    transport model in fractured porous media. An accurate description of
    flow and transport in the fractures is key to obtain reliable simulations, however, fractures geometry and physical characteristics pose several
    challenges from both the modeling and implementation side. We adopt a
    mixed-dimensional approximation, where fractures and their intersections
    are represented as objects of lower dimension. To simplify the
    presentation, we consider only two chemical species: one solute,
    transported by water, and one precipitate attached to the solid
    skeleton. A global sensitivity analysis to uncertain input data is performed exploiting the Polynomial Chaos expansion along with spectral projection methods on sparse grids.
\end{abstract}

\section{Introduction}\label{intro}

The Paris agreement, adopted by 196 parties in 2015, aims at limiting global
warming to below 2$^\circ C$, preferably to 1.5$^\circ C$,  compared to
pre-industrial levels. The reduction of greenhouse gas emissions in the
atmosphere is crucial to achieve such long term goal, and requires the transition
towards renewable energies, and the subsequent need for effective energy
storage; moreover, a safe long term sequestration of $CO_2$ is considered a
promising strategy to reduce emissions into the atmosphere.  Many of the
aforementioned strategies entail a massive use of the subsurface for fluids
injection, storage and production. If, on one hand, it is necessary to guarantee
the mechanical integrity of the subsurface to avoid unwanted fracturing and induced
seismicity, it is also important to evaluate the effect of chemical reactions on
the hydraulic properties of the porous media. Indeed, the injection of water at
different temperature and with different solutes concentration with respect to
the formation can cause dissolution and precipitation of minerals, with an
impact on porosity and permeability. In some cases, as in the case of \emph{in
situ carbon mineralization}, these reactions can be exploited to our advantage to
obtain a better storage of $CO_2$ by 1) chemically trapping carbon by binding
it into existing mineral and, at the same time 2) thanks to the change of
specific volume associated with the transformation of the minerals, creating a
more effective cap-rock, \cite{Ling2021}. Even if this process occurs over very long time scales
in natural conditions, its controlled exploitation is an interesting
technological challenge.

A realistic mathematical model of these phenomena encompasses a model for flow
in porous media (we focus on single phase flow, assuming small concentrations of
gases, if present), coupled with the transport of mobile species and the
modeling of both kinetic and equilibrium reactions. Moreover, since reaction
rates are usually influenced by temperature, and in view of the possible
application of the model to low temperature geothermal plants, the coupled model
is completed by the heat equation, \cite{Fumagalli2020a}.

Since fractures are ubiquitous in porous media and have a major impact on flow
and transport (both of solutes and heat) this work is focused on the modeling of
fractured porous media, where fractures are modeled as lower dimensional objects,
which can be much more permeable than the surrounding medium, or nearly
impermeable if, for instance, precipitation occurs reducing their aperture, \cite{Fumagalli2020f}.

Many physical, geometrical and geochemical parameters involved in the model are
affected by uncertainty, therefore a purely deterministic evaluation of the
model is useless without a suitable analysis of the solution variance. In this
work we propose an uncertainty quantification workflow based on Polynomial Chaos  \cite{Wiener:38,Ghanem:91}
and sparse grids  \cite{Conrad:13} for the sampling of the parameters space: this choice allows us
to obtain accuracy with a manageable number of evaluations of the model, which,
being coupled and time dependent, has a non-negligible computational cost. The
goal is to compute the Sobol indices associated with some input parameters to
quantify the impact of the uncertainty of such quantities on a variable of
interest, typically the medium porosity. The space distribution of the Sobol
indices, or the partial variances, can also give interesting insight into the
problem, which, in spite of being reduced to the minimum possible complexity,
already exhibits a non-trivial, fully coupled behavior.

The paper is structured as follows. In Section 2 we present the model equation in a homogeneous porous medium; then, in Section 3, we extend the equations to the hybrid dimensional case to account for the coupling with fractures.
Section 4 briefly discusses the proposed numerical approximation schemes and some implementation details.
Section 5 describes the uncertainty quantification workflow and Section 6 is devoted to the presentation of a a complete set of test cases. Finally, conclusions are drawn in Section 7.

\section{Problem description and mathematical model}
\label{sec:1}

We consider the coupled problem of single-phase flow and reactive transport in
porous media, accounting for the porosity changes linked to mineral reactions.
Single phase flow can be assumed, even in the presence of $CO_2$, at low gas
concentrations, for instance at the edge of the plume, away from injection
wells. Typically, a large number of species is needed for a realistic modeling
of reactive transport, including solid mineral species and aqueous complexes,
and the chemical reactions include equilibrium and kinetic ones. In the
following we will discuss the general case and then specialize the model to the
simple case of the dissolution/precipitation of a single mineral species.

\subsection{Single phase flow}

Under the assumption of single-phase flow, the fluid pressure $p$ and Darcy velocity
$\bm{q}$ can be computed as the solution of a Darcy problem
\begin{gather}\label{eq:darcy}
    \begin{aligned}
        &\bm{q} + \dfrac{k(\phi)}{\mu}(\nabla p - \rho_w g \bm{e}_z ) = 0\\
        &\nabla \cdot \bm{q} = -\partial_t \phi
    \end{aligned}
\end{gather}
complemented by suitable boundary conditions prescribing a given pressure or
normal flux on the boundary. Here $k$ is the permeability, $\mu$ the viscosity, $\rho_w$ is the fluid density and $g$ is the gravity
acceleration. Note that in \eqref{eq:darcy} we have assumed constant density,
but we are accounting for porosity changes through a source term. Indeed,
porosity depends on the changes of the mineral volume fraction as detailed in
Section \ref{sec:min}.

The intrinsic permeability, which we
assume to be isotropic in the bulk porous medium, can be modeled as a function
of porosity. We consider the following law
\begin{gather*}
    k(\phi)=k_0\left(\dfrac{\phi}{\phi_0}\right)^2
\end{gather*}
where $k_0$ is the reference value at the initial porosity $\phi_0$. Note that
this dependence will introduce a non-linear coupling among the model equations.

\subsection{Advection-diffusion-reaction equations for mobile species}

Let $u_i$ denote the molar concentration of the $i$-th \emph{mobile} species, for instance a dissolved mineral or gas. The governing equation
reads
\begin{gather}\label{eq:ADR}
    \partial_t (\phi u_i) + \nabla \cdot (\bm{q} u_i - \phi d\nabla u_i) =
    \phi \displaystyle\sum_{r=1}^{N_r} \nu_{ir} R_r
\end{gather}
where $d$ is the diffusion tensor (which may account for mechanical
dispersion, and is considered to be the same for all species), $\bm{q}$ is the Darcy velocity, $N_r$ is the number of reactions
and $\nu_{ir}$ is the stoichiometric coefficient of species $i$ in reaction $r$.

In particular, we consider a set of $N_r$ reactions in the form
\begin{gather*}
    \sum_i \nu_{ir} X_i \rightleftarrows 0\qquad r=1,\ldots,N_r
\end{gather*}
where $X_i$ is a generic species, solid or mobile, $\nu_{ir}$ is zero if $X_i$
is not involved in the reaction $r$, negative if $X_i$ is a reactant in the forward
reaction, and positive if it is a product. The reaction rate $R_r$ should be
interpreted as the \emph{net} rate, i.e.
\begin{gather*}
    R_r= R_r^{+} - R_r^{-}
\end{gather*}
where $R_r^+$ is the forward reaction rate and $R_r^-$ is the backward one.
Their expressions are often empirical and problem dependent; the one considered
in this work is illustrated in Section \ref{sec:ourcase}.  Finally, equation
\eqref{eq:ADR} should be complemented by initial conditions $u_i(\bm{x},0) =
u_{i,0}(\bm{x})$, and boundary conditions on the concentration or normal flux.

\subsection{Mineral species}\label{sec:min}

In the case the species is \emph{immobile} the evolution equation for its concentration, denoted as $w_i$, reduces to an ordinary differential equation since we do not have the transport and diffusion terms, i.e.
\begin{gather}\label{eq:R}
    d_t (\phi w_i)= \phi \displaystyle\sum_{r=1}^{N_r} \nu_{ir} R_r
\end{gather}
where an initial condition is supplemented as $w_i(\bm{x},0) = w_{i,0}(\bm{x})$.
However, it may be more convenient to consider a different measure of
concentration for solid species, in particular we want to compute the
\emph{solid volume fractions} $\phi_i$ as the ratio of the volume of mineral $i$
for a unit volume of rock. The volume fractions can be obtained from the
concentrations as

\begin{gather}\label{eq:conversion}
    \phi_i=w_i \eta_i \phi
\end{gather}
where $\eta_i$ is the molar volume of the mineral. If we let $\phi_I$ be the
volume fraction of inert minerals (i.e. species that are not affected by
chemical reactions) we have that
\begin{gather*}
    \begin{aligned}
        &\phi=1-\phi_I - \displaystyle \sum_{i = 1}^{N_s} \phi_i\\
        &d_t \phi = - \displaystyle \sum_{i = 1}^{N_s}\eta_i d_t (\phi  w_i)
    \end{aligned}
\end{gather*}
where $N_s$ is the number of solid, immobile species.

\subsection{Heat equation}

Since reaction rates usually depend on temperature, we include in our model the
heat equation for thermal conduction (based on Fourier's law) and convection in
the porous media. We assume thermal equilibrium between rock and water, so that
only one equation can be considered instead of two.

The temperature field is indicated as $\theta$ and its evolution is described by
\begin{align}\label{eq:heat}
    &\begin{aligned}
        &\partial_t (c(\phi) \theta) +
        \nabla \cdot ( \rho_w c_w \bm{q} \theta - \Lambda(\phi) \nabla \theta  )    + j = 0.
    \end{aligned}
\end{align}
Here $c$ is the effective thermal capacity defined as the  porosity-weighted
average of the water $c_w$ and solid $c_s$ specific thermal capacities,
\begin{gather*}
    c(\phi) = \phi \rho_w c_w + (1-\phi) \rho_s c_s,
\end{gather*}
where $\rho_w$ and $\rho_s$ are the densities of the water and solid
phase respectively. The effective thermal conductivity  $\Lambda$ is computed in a similar way, as
\begin{gather*}
    \Lambda(\phi) = \Lambda_w^\phi \Lambda_s^{1-\phi},
\end{gather*}
where $\Lambda_s$ and $\Lambda_s$ are the water and solid thermal conductivity.
Finally, $j$  models a source or sink of heat in the system. Equation
\eqref{eq:heat} is completed by initial conditions,
$\theta(\bm{x},0)=\theta_0(\bm{x})$, and boundary conditions on the temperature
or heat flux.

\subsection{Simplified dissolution-precipitation model}\label{sec:ourcase}

In this work we consider a single, simple kinetic reaction in the form
\begin{gather*}
    U+V -W\leftrightarrow 0
\end{gather*}
where $U$ and $V$ are the positive and negative ion respectively, and $W$ is the
precipitate they can form, \cite{Radu2010}. The net reaction rate results from the difference
between precipitation (forward reaction) and dissolution (backward reaction). Under
the hypotheses of electrical equilibrium, we can assume that the concentration of
$U$ and $V$ are equal and denoted by $u$, which from now on is the molar
concentration of the mobile species, whereas $w$ is the molar concentration of
the precipitate.
The reaction rates are modeled as
\begin{gather*}
    \begin{aligned}
        &r_p = \lambda_-(\theta) u^2 \\
        &r_d = \lambda_+(\theta)
    \end{aligned}
\end{gather*}
moreover, at equilibrium we have that $r_d=r_p$, thus, $\lambda^+=\lambda$ and $\lambda_-= \lambda u_{e}^{-2}$ where $u_e$ is the equilibrium solute
concentration. The net reaction rate is then
\begin{gather*}
    r_w(u,w,\theta) =
    \begin{dcases*}
        \displaystyle \lambda \left(\frac{u^2}{u_e^2}-1\right) & if $w>0$\\
        \lambda \frac{u^2}{u_e^2} & if $w\leq 0$.
    \end{dcases*}
\end{gather*}
Note that precipitation proceeds at a constant rate until $w=0$, and is then set to zero.

\subsection{The complete model}
The complete model, in the case of our interest, computes $(\bm{q}, p, \theta,
u, w, \phi)$ by solving the following system of non-linear and coupled equations
given by
\begin{gather}\label{eq:model}
    \begin{aligned}
        & \bm{q} + \dfrac{k(\phi)}{\mu}(\nabla p - \rho g \bm{e}_z ) = 0\\
        & \nabla \cdot \bm{q} = -\partial_t \phi\\
        & \partial_t [c (\phi) \theta] + \nabla \cdot ( \rho_w c_w \bm{q} \theta -
        \Lambda(\phi) \nabla \theta  )    + j = 0\\
        &\partial_t (\phi u) + \nabla \cdot (\bm{q} u - \phi d\nabla u) = -\phi r_w(u,w,\theta)\\
        &d_t (\phi w) = \phi r_w(u,w,\theta)\\
        &d_t \phi = - \eta d_t( \phi w)
    \end{aligned}
\end{gather}
complemented by constitutive laws and suitable initial and boundary conditions.

\section{Hybrid dimensional model for fractured porous media}

\begin{figure}[tb]
    \centering
    \resizebox{0.33\textwidth}{!}{\fontsize{1cm}{2cm}\selectfont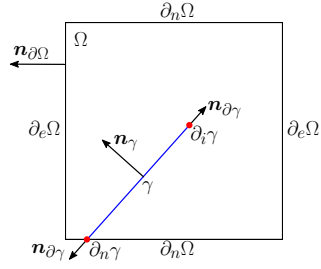}
    \caption{Fractured domain and notation.}%
    \label{fig:fracture_domain}
\end{figure}

In the following we introduce the extension of the coupled model \eqref{eq:model} to the case of fractured porous media, following \cite{Fumagalli2020f}.
For the sake of simplicity we will present the model in the case of a single fracture, geometrically reduced to its centerline, see Figure \ref{fig:fracture_domain}. This model reduction strategy is often adopted in the simulation of fractured porous media to reduce the computational cost by avoiding excessive mesh refinement; moreover, in our case it is particularly convenient since the aperture can change in time due to reactions. For more references on this approach see
\cite{Martin2005,Jaffre2011,Sandve2012,Nordbotten2018,Chave2018,Ahmed2018,Berre2019b,Fumagalli2020a} and references therein.

Let $\tilde{\Omega}$ be the fractured domain. Following \cite{Formaggia2012}, we define  $\gamma$ as a non self-intersecting $C^2$ curve (if
$n=2$) or surface (if $n=3$). In an equi-dimensional setting a fracture can be defined as the following set of points
\begin{gather*}
 \Gamma=   \left\{ \bm{x} \in \tilde{\Omega}: \,\bm{x} = \bm{s} + r \bm{n}_\gamma, \, \bm{s} \in
    \gamma,\, r \in \left(-\frac{\epsilon_\gamma(\bm{s})}{2},
    \frac{\epsilon_\gamma(\bm{s})}{2}\right) \right\}.
\end{gather*}
Thus, we have a subdomain of $\tilde{\Omega}$, separated by the surrounding porous medium $\Omega$ by the interfaces
$\gamma_+$ and $\gamma_-$ with associated normal vectors $\bm{n}_+$ and $\bm{n}_-$. We replace the fracture $\Gamma$ with its centerline $\gamma$ and we assign a {unique} normal $\bm{n}_\gamma=\bm{n}_{\gamma,+}$ to the fracture, see Figure \ref{fig:fracture_domain}.
We assume that the fracture is open, with porosity $\phi_\gamma=1$.

\subsection{Reduced variables}

In the following we denote with the subscript $\gamma$ the variables in the
fracture, and with $\Omega$ the variables in the porous medium. Note that, after
geometrical reduction, we denote with $\Omega$ the domain
$\tilde{\Omega}\setminus \gamma$. Reduced vectors variables in the fracture are defined
as the integral of the tangential components of the corresponding equi-dimensional variables,
\begin{gather*}
    \bm{q}_\gamma(\bm{x}) \defeq \int_{\epsilon_\gamma(\bm{x})} T(\bm{x}) \bm{q}(\bm{x}, s) d s
\end{gather*}
with $T\defeq I - N$ and $N \defeq \bm{n}_\gamma \otimes \bm{n}_\gamma$ the
tangential and normal projection
matrices, respectively. The reduced scalar variables are instead the integral
average, for each section of the fracture, of the equi-dimensional counterparts
\begin{gather*}
    p_\gamma(\bm{x}) \defeq \dfrac{1}{\epsilon_\gamma(\bm{x})}\int_{\epsilon_\gamma(\bm{x})}
    p(\bm{x}, s) ds \quad
    \theta_\gamma(\bm{x}) \defeq \dfrac{1}{\epsilon_\gamma(\bm{x})}
    \int_{\epsilon_\gamma(\bm{x})} \theta(\bm{x}, s) ds\\
    u_\gamma(\bm{x}) \defeq \dfrac{1}{\epsilon_\gamma(\bm{x})}\int_{\epsilon_\gamma(\bm{x})}
    u(\bm{x}, s) ds \quad
    w_\gamma(\bm{x}) \defeq \dfrac{1}{\epsilon_\gamma(\bm{x})}\int_{\epsilon_\gamma(\bm{x})}
    w(\bm{x}, s) ds.
\end{gather*}

Moreover, we denote by $\nabla_\tau$, $\nabla_\tau \cdot$ the tangential gradient and divergence defined on the tangential space of the fracture.

Following \cite{Martin2005} we assume that the permeability $k$ and diffusivity $d$, can be decomposed in normal and tangential components as
\begin{gather}\label{eq:reduced_perm}
    k = \kappa_\gamma N + k_\gamma T
    \quad\text{and}\quad
    d = \delta_\gamma N + d_\gamma T.
\end{gather}

\subsection{Reduced Darcy flow model}


The reduced model for the Darcy flow, which describes the evolution of the reduced
Darcy velocity $\bm{q}_\gamma$ and pressure $p_\gamma$ in the fracture is obtained, following \cite{Martin2005}, by the following steps:
\begin{itemize}
 \item integration of the mass balance equation in each section of the fracture, to obtain a conservation equation for $\bm{q}_\gamma$;
 \item integration of the tangential component of Darcy law in each section of the fracture to obtain a relationship between $\bm{q}_\gamma$ and $p_\gamma$;
 \item integration, with a suitable approximation,  of the  normal component of Darcy law in each section of the fracture to obtain two coupling conditions between the fracture and the surrounding medium.
\end{itemize}
The coupled problem in $\Omega$ and $\gamma$ then reads:
\begin{align}\label{eq:reduced_problem_darcy}
    &\begin{aligned}
        &\mu \bm{q}_\Omega +  k_\Omega(\phi_\Omega) \nabla p_\Omega =\bm{0}\\
        &\partial_t \phi_\Omega + \nabla_\gamma \cdot \bm{q}_\Omega +
        f_\Omega =0
    \end{aligned}
    &\quad \text{in } \Omega\times (0, T),\\
    &\begin{aligned}
        &\mu \bm{q}_\gamma + \epsilon_\gamma k_\gamma(\epsilon_\gamma) \nabla p_\gamma =\bm{0}\\
        &\partial_t \epsilon_\gamma + \nabla_\gamma \cdot \bm{q}_\gamma + \bm{q}_\gamma\cdot\bm{n}_\gamma|_{\gamma_+}  - \bm{q}_\gamma\cdot\bm{n}_\gamma|_{\gamma_-} +
        \epsilon_\gamma f_\gamma =0
    \end{aligned}
   &\quad \text{in } \gamma \times (0, T),\\
   &    \mu \epsilon_\gamma  \bm{q}_\Omega \cdot \bm{n}_{\gamma}|_{\gamma_+} + \kappa_\gamma(\epsilon_\gamma) (p_\gamma - p_\Omega|_{\gamma_+}) =
    0 &\quad \text{on } \gamma_+ \times (0, T)\\
    &    \mu \epsilon_\gamma  \bm{q}_\Omega \cdot \bm{n}_{\gamma}|_{\gamma_-} + \kappa_\gamma(\epsilon_\gamma) (p_\gamma - p_\Omega|_{\gamma_-}) =
    0 &\quad \text{on } \gamma_- \times (0, T),
\end{align}
where $f_\gamma(\bm{x}) \defeq \epsilon_\gamma^{-1}(\bm{x})\int_{\epsilon_\gamma(\bm{x})} f(\bm{x}, s) ds$ is the reduced source or sink term.
Following lubrication theory, the fracture tangential
permeability $k_\gamma$ can be expressed as a function of the aperture, as described
in more detail in Subsection \ref{eq:reduced_permeability_aperture}.

The conditions on $\gamma_\pm$ model the fact that the flux exchange between the fracture
and the surrounding porous media is related to the pressure jump via
$\kappa_\gamma$. Note that, since $\kappa_\gamma$, as $k_\gamma$, can be modeled as a function of $\epsilon_\gamma$, if the aperture goes to zero the flux exchange vanishes.

\subsection{Reduced heat model}

The reduced model that describes the evolution of temperature $\theta$ is
obtained, similarly to the Darcy problem, by integrating the  conservation
equation in each section of the fracture; the coupling conditions however, which
stem from a suitable approximation of the total normal heat flux,  should take
into account the different nature of the advective and diffusive fluxes in the coupling. The coupled model in $\Omega$ and $\gamma$ reads:
\begin{align}\label{eq:reduced_theat}
    \begin{aligned}
        &\partial_t (c(\phi) \theta_\Omega)+ \nabla \cdot (\rho_w c_w \bm{q} - \Lambda(\phi) \nabla \theta_\Omega) + j=0
        &&\text{in } \Omega\times (0, T),\\
        &\partial_t (\epsilon_\gamma c_w \theta_\gamma)+ \nabla_\tau \cdot (\rho_w c_w \bm{q}_\gamma - \Lambda_w\epsilon_\gamma \nabla_\tau \theta_\gamma) + \psi^+ - \psi^- +j_\gamma = 0
        &&\text{in } \gamma \times (0, T).
    \end{aligned}
\end{align}
where the conservation equation in the fracture accounts for heat flux exchanged with the fracture on both sides, through the terms $\psi^\pm$, defined as
\begin{align*}
 & \psi^+ = \rho_w c_w\bm{q}_\Omega\cdot\bm{n}_\gamma|_{\gamma_+} \tilde{\theta}_+ + \dfrac{2 \Lambda_w}{\epsilon_\gamma}(\theta_\Omega|_{\gamma_+} - \theta_\gamma)\\
 & \psi^- = \rho_w c_w\bm{q}_\Omega\cdot\bm{n}_\gamma|_{\gamma_-} \tilde{\theta}_- + \dfrac{2\Lambda_w}{\epsilon_\gamma}( \theta_\gamma - \theta_\Omega|_{\gamma_-})
\end{align*}
where $\tilde{\theta}_\pm$ is selected as the upwind value, i.e.
\begin{align*}
 &\tilde{\theta}_+=\begin{cases}
  \theta_\gamma &\quad\mathrm{if}\,\bm{q}_\Omega\cdot\bm{n}_\gamma|_{\gamma_+}>0\\
  \theta_\Omega|_{\gamma_+} &\quad\mathrm{if}\,\bm{q}_\Omega\cdot\bm{n}_\gamma|_{\gamma_+}<0
 \end{cases} \qquad
 &\tilde{\theta}_-=\begin{cases}
  \theta_\gamma&\quad\mathrm{if}\,\bm{q}_\Omega\cdot\bm{n}_\gamma|_{\gamma_-}<0\\
  \theta_\Omega|_{\gamma_-} &\quad\mathrm{if}\,\bm{q}_\Omega\cdot\bm{n}_\gamma|_{\gamma_-}>0.
 \end{cases}
\end{align*}

\subsection{Reduced solute and precipitate model}

The reduced model that describes the evolution of the solute $u_\gamma$ is similar to the reduced heat equation. By integrating the solute conservation equation in the fracture we obtain the reduced equation, such that
\begin{align}\label{eq:reduced_transport}
    \begin{aligned}
        &\partial_t (\phi_\Omega u_\Omega) + \nabla \cdot (\bm{q}_\Omega u_\Omega - d_\Omega \nabla u_\Omega) = \phi_\Omega r_w
        &&\text{in } \Omega\times (0, T),\\
        &\partial_t (\epsilon_\gamma u_\gamma) + \nabla_\tau \cdot (\bm{q}_\gamma u_\gamma - \epsilon_\gamma d_\gamma \nabla_\tau u_\gamma) + \chi^+ - \chi^-  = \epsilon_\gamma r_w
        &&\text{in } \gamma \times (0, T).
    \end{aligned}
\end{align}
Note that the balance equation in the fracture accounts for exchanges with the porous medium, in particular we have that
 \begin{align*}
  & \chi^+ = \bm{q}_\Omega\cdot\bm{n}_\gamma|_{\gamma_+} \tilde{u}_+ + \dfrac{2\delta_\gamma}{\epsilon_\gamma}(u_\Omega|_{\gamma_+} - u_\gamma)\\
  & \chi^- = \bm{q}_\Omega\cdot\bm{n}_\gamma|_{\gamma_-} \tilde{u}_- + \dfrac{2\delta_\gamma}{\epsilon_\gamma}( u_\gamma - u_\Omega|_{\gamma_-})
 \end{align*}
where once again $\tilde{u}_\pm$ is the upwind value, i.e.
 \begin{align*}
  &\tilde{u}_+=\begin{cases}
   u_\gamma &\quad\mathrm{if}\,\bm{q}_\Omega\cdot\bm{n}_\gamma|_{\gamma_+}>0\\
   u_\Omega|_{\gamma_+} &\quad\mathrm{if}\,\bm{q}_\Omega\cdot\bm{n}_\gamma|_{\gamma_+}<0
  \end{cases} \qquad
  &\tilde{u}_-=\begin{cases}
   u_\gamma&\quad\mathrm{if}\,\bm{q}_\Omega\cdot\bm{n}_\gamma|_{\gamma_-}<0\\
   u_\Omega|_{\gamma_-} &\quad\mathrm{if}\,\bm{q}_\Omega\cdot\bm{n}_\gamma|_{\gamma_-}>0.
  \end{cases}
 \end{align*}

For the precipitate in the fracture $w_\gamma$, being the original model an
ordinary differential equation valid for each point of the domain, the reduced
model in the fracture becomes simply
\begin{align}\label{eq:reduced_reaction}
    &\partial_t (\epsilon_\gamma w_\gamma) - \epsilon_\gamma r_w(u_\gamma,
    w_\gamma, \theta_\gamma) = 0 &&
    \text{in } \gamma\times (0, T).
\end{align}
Note that \eqref{eq:reduced_reaction} is not directly coupled with the corresponding equation in the porous matrix since both describe local phenomena.

\subsection{Permeability and aperture model}\label{eq:reduced_permeability_aperture}

As already mentioned we assume that both components of the permeability $k$ in the fracture follow a law which relates them to the aperture, more precisely
\begin{gather}\label{eq:permeability_aperture}
    k_\gamma(\epsilon_\gamma) = k_{\gamma,0}
    \dfrac{\epsilon_\gamma^2}{\epsilon_{\gamma, 0}^2}
    \quad \text{and} \quad
    \kappa_\gamma(\epsilon_\gamma) =
    \kappa_{\gamma, 0}\dfrac{\epsilon_\gamma^2}{\epsilon_{\gamma, 0}^2},
\end{gather}
where $k_{\gamma,0}$ and $\kappa_{\gamma, 0}$
are reference coefficients along and across the
fracture, respectively, and $\epsilon_{\gamma, 0} > 0$  is the initial aperture.
As the porosity changes with mineral precipitation, we consider a similar law to describe the
evolution of the fracture aperture $\epsilon_\gamma$. We have
\begin{gather}\label{eq:reduced_aperture}
    \begin{aligned}
        &\partial_t \epsilon_\gamma  + \eta_\gamma \partial_t(
         \epsilon_\gamma w_\gamma )= 0&&
        \text{in } \gamma \times (0, T)\\
        &\epsilon_\gamma(t=0) = \epsilon_{\gamma, 0} && \text{in } \gamma \times \{0\}
    \end{aligned},
\end{gather}
where $\eta_\gamma$ represents the molar volume of the mineral as it precipitates on the fracture
walls.

\section{Numerical approximation}

The numerical schemes for the solution of the deterministic problem \eqref{eq:model} are
implemented in the PorePy library \cite{Keilegavlen2019} which provides support for
multidimensional coupling, allowing for an easy implementation of the problem in
fractured media.

\subsection{Time integration and splitting strategy}

Problem \eqref{eq:model} is fully coupled in a non-linear way. For the sake of
computational efficiency, in this work its solution  is based on a non-iterative
splitting strategy, with the underlying assumption that the changes to the flow
parameters due to chemical reactions are relatively slow. In particular, at each
time step we follow the scheme proposed in \cite{Fumagalli2020f} and:
\begin{enumerate}
    \item we first solve the Darcy problem to obtain the advective fields
    $\bm{q}_\Omega$, $\bm{q}_\gamma$. The flow problem is discretized in time by
    approximating the time derivative $\partial_t \phi$ by finite differences as
    \begin{gather*}
        \partial_t\phi \simeq \dfrac{\phi^*-\phi^n}{\Delta t}
    \end{gather*}
    where
    $\phi^*=2\phi^n-\phi^{n-1}$ is the extrapolated value;
    \item with $\bm{q}_\Omega$, $\bm{q}_\gamma$ we solve the heat equation,
    discretized in time with the Implicit Euler method;
    \item then, given the temperature field we solve the
    advection-diffusion-reaction problem which is in turn split into
    \begin{enumerate}
        \item the advective step, discretized in time with the Implicit Euler
        scheme, which gives and intermediate solute concentration $u^*$;
        \item the reaction step to compute the final solute and precipitate
        concentrations (note that the precipitate is not affected by transport).
        This step is integrated explicitly in time, with the addition of an
        event location procedure to avoid negative precipitate concentrations.
    \end{enumerate}
    \item Finally, we update porosity and permeability for the next step.
\end{enumerate}

\subsection{Space discretization}

Space discretization is based on a standard, conforming approach where fractures
are honoured by the computational grid and each element of the fractures grid is
a face of the porous media grid. However, this assumption could be relaxed
allowing for different grid resolutions with the use of mortar variables.
Finally, since equations are in mixed-dimensions, all the numerical schemes are
applied in different dimensions, i.e. in 2D and 1D.

Since the Darcy flux is involved in the advective terms of the transport and heat
equations it is of fundamental importance that local mass conservation is
fulfilled. Therefore, we approximate the Darcy problem in its mixed formulation
and employ a suitable pair of discrete spaces for the pressure and the Darcy
flux. In particular we employ the lowest order Raviart-Thomas element pair $\mathbb{RT}_0$, $\mathbb{P}_0$.

For the numerical solution of the heat equation and the advection-diffusion-reaction equation we apply, in hybrid dimensions, the Finite
Volume method. Consistently with the continuous model we consider an upwind
approximation of the advective term, whereas the diffusive term is approximated
with the two point flux approximation (TPFA), see \cite{Eymard2000}, \cite{Faille2002}, \cite{Droniou2013}. Since we are
considering constrained triangulations to honour the fractures the grid may in
principle not be orthogonal. However, we assume that the distortion is small
enough to obtain a reliable approximation even with a simple TPFA scheme.

\section{Sensitivity analysis workflow}

In this section we present the algorithm employed to approximate stochastic
quantities by means of Polynomial Chaos (PC) expansions~\cite{Wiener:38,Ghanem:91}. This
technique will allow us to compute the sensitivity Sobol indices and to obtain a
surrogate model of the problem for a quick evaluation of the quantities of
interest. PC expansion have been used to treat a large variety of problems,
including elliptic models (see, e.g.,~\cite{Babuska:02,Matthies.Keese:05}),
fluid mechanics problems~\cite{LeMaitre:02}, and flow and transport in porous
media (see, e.g.,~\cite{Ghanem.Dham:98,Botti.ea:20}).

The sampling of the uncertain parameters space is performed with pseudo-spectral
projection on \emph{sparse grids}~\cite{Conrad:13}, thus obtaining an accurate
estimate with a limited number of evaluation of the deterministic model. This is
particularly important since the problem is time dependent and, moreover, the
presence of chemical reactions can introduce a fast time scale, constraining the
time step amplitude with an increase in the computational cost for each evaluation.

\subsection{Polynomial Chaos expansion}

Let $N$ be the number of parameters $\bm{\xi}= (\xi_i)_{1\le i \le N}$
and $\Xi$ the space of possible realizations. For the sake of simplicity the
parameters are rescaled so that $\Xi=[0,1]^N$. Moreover, given a probability
measure $\rho:\Xi\to\mathbb{R}^+$, the inner product of two second-order random
variables $X(\bm{\xi})$ and $Y(\bm{\xi})$ is defined as
\begin{gather*}
    \langle X, Y \rangle = \int_\Xi X(\bm{\xi}) Y(\bm{\xi}) \rho(\bm{\xi}) \mathrm{d}\bm{\xi}.
\end{gather*}
The Polynomial Chaos (PC) expansion of a variable $X(\bm{\xi})$ reads
\begin{gather}\label{eq:PCexp}
    X(\bm{\xi}) = \sum_{\bm{k}\in\mathbb{N}^N} X_{\bm{k}} \phi_{\bm{k}}(\bm{\xi}),
\end{gather}
where $\{X_{\bm{k}}=\langle X, \phi_{\bm{k}}\rangle  : \bm{k}\in \mathbb{N}^N\}$
are the spectral modes of $X$, the basis functions  $\{ \phi_{\bm{k}}(\bm{\xi})
: \bm{k}\in \mathbb{N}^N\}$ are multi-variate polynomials chosen to be
orthogonal with respect to the product $\langle\cdot,\cdot\rangle$, and the
multi-index $\bm{k}=(k_1,\ldots, k_N)$ denotes the polynomial degree with
respect to the parameters $\xi_i$.

The PC approximation is then obtained by truncating the expansion in \eqref{eq:PCexp} to a finite set $\mathcal{K}\subset\mathbb{N}^N$, which determines the quality of the approximation:
\begin{gather}\label{eq:PCdef}
    X_{\mathcal{K}}(\bm{\xi}) := \sum_{\bm{k}\in\mathcal{K}} X_{\bm{k}}\phi_{\bm{k}}(\bm{\xi}).
\end{gather}
The statistical moments of the variables of interest are easily obtained from the PC approximations; e.g., mean, variance, and covariance are given by
\begin{gather}\label{eq:Esp_Var}
    \langle X_{\mathcal{K}} \rangle = X_{\bm0}, \quad
    \Var(X_{\mathcal{K}}) =\sum_{\bm{k}\in\mathcal{K}\setminus\bm0} X_{\bm{k}}^2, \quad
    \Cov(X_{\mathcal{K}}, Y_{\mathcal{K}}) = \sum_{\bm{k}\in\mathcal{K}\setminus\bm0} X_{\bm{k}} Y_{\bm{k}}.
\end{gather}

\subsection{Spectral projection method}\label{sec:PSP}

The coefficients of the PC expansion in \eqref{eq:PCdef} can be computed in
different ways (see, e.g., \cite{LeMaitre:10,Constantine:12}). This work is
based on non-intrusive pseudo-spectral projection, which in our opinion,
provides the best trade-off between complexity and precision. The numerical
quadrature schemes are constructed as sparse tensorization of a one dimensional
formula~\cite{Gerstner:03}. Then, given $M$ quadrature points and the respective
weights  $\{w^{(q)}\}_{1\le q \le M}$, the modes
$(X_{\bm{k}})_{\bm{k}\in\mathcal{K}}$ are computed as
\begin{gather*}
    X_{\bm{k}} = \langle X,\phi_{\bm{k}}\rangle
    \simeq \sum_{q=1}^{M} {w^{(q)} X(\bm{\xi}^{(q)}) \phi_{\bm{k}}(\bm{\xi}^{(q)}) }.
\end{gather*}

\begin{figure}[tbp]
    \centering
    \includegraphics[width=0.3\textwidth]{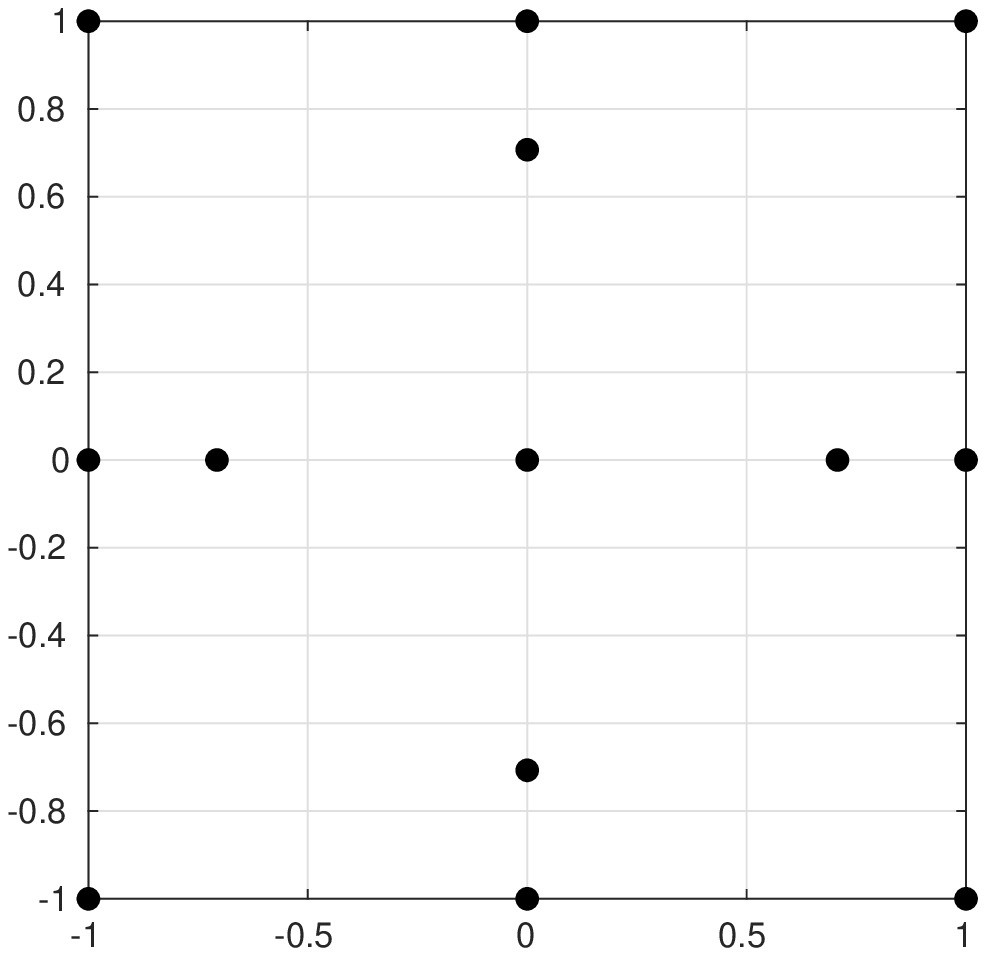}%
    \hspace*{0.05\textwidth}
    \includegraphics[width=0.3\textwidth]{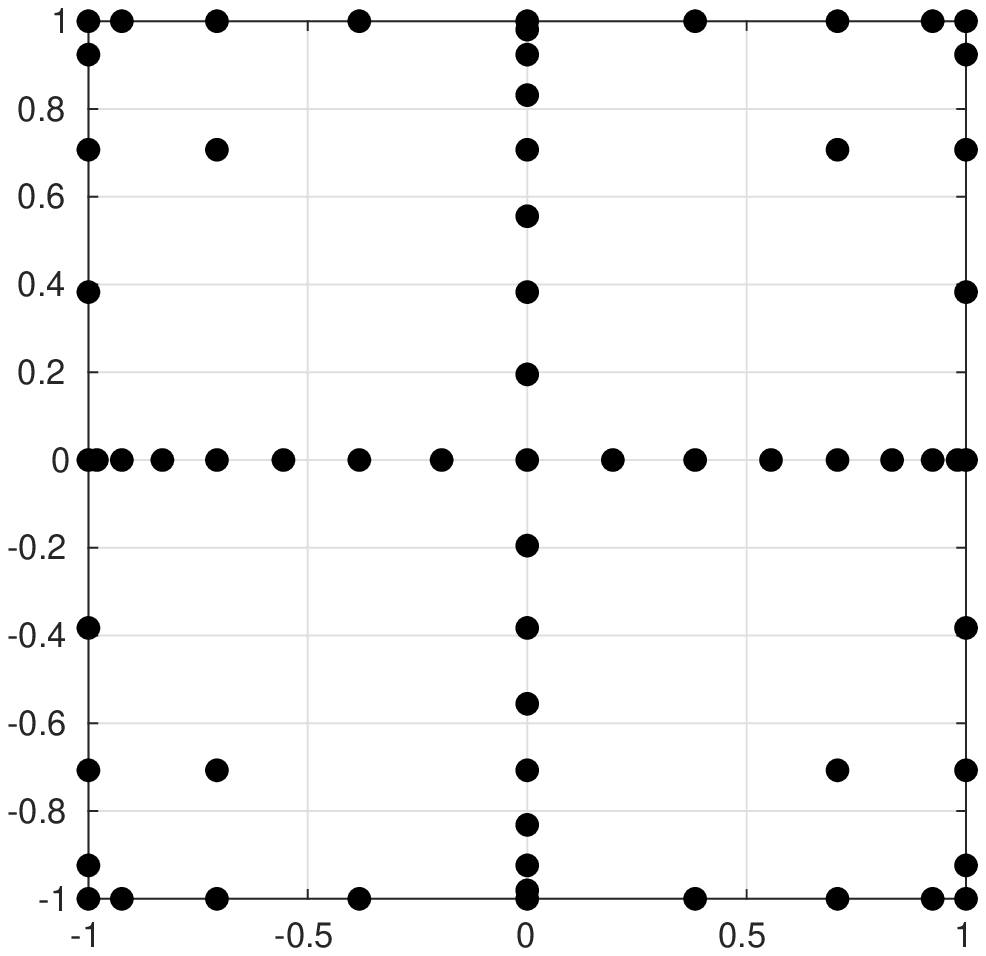}%
    \hspace*{0.05\textwidth}
    \includegraphics[width=0.3\textwidth]{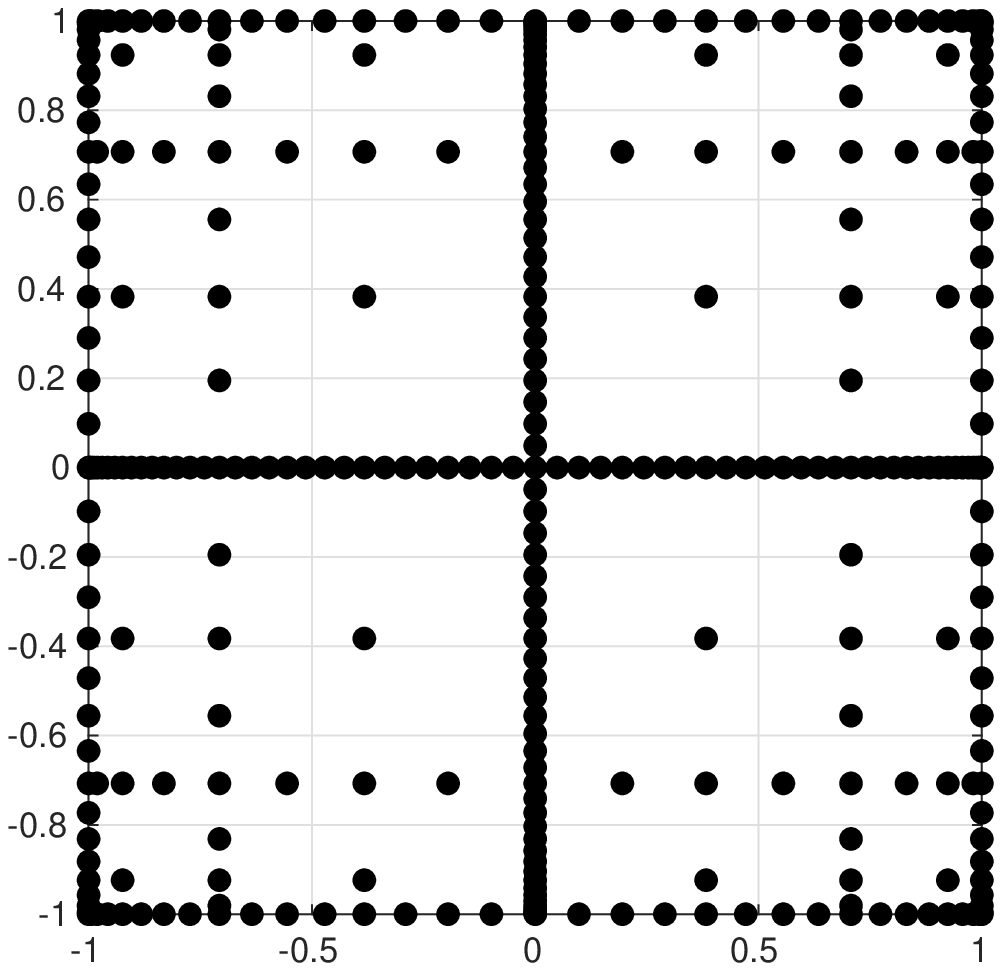}
    \caption{{Sparse grid} for two parameters of level 2 on the left, level 4 on the centre
    and level 6 on the right.}%
    \label{fig:SG}
\end{figure}

The complexity of the method is governed by the number $M$ of evaluations of the
deterministic problem, while the accuracy depends on the PC basis
$\{\phi_{\bm{k}}\}_{k\in\mathcal{K}}$. In order to maximize the accuracy with
respect to the computational effort, we adopt a sparse method (cf. Figure
\ref{fig:SG}) hinging on the application of Smolyak's formula~\cite{Smolyak:63}
directly on the projection operator, rather than on the integration operator.
Specifically, the set $\mathcal K$ is defined as the largest possible one such
that the discrete projection is exact for any function spanned by
$\{\phi_{\bm{k}}\}_{k\in\mathcal{K}}$, i.e.:
\begin{gather*}
   \forall \bm{k}, \bm{l} \in \mathcal{K}, \:
  \sum_{q=1}^{M} {w^{(q)} \phi_{\bm{k}}(\bm{\xi}^{(q)})
  \phi_{\bm{l}}(\bm{\xi}^{(q)}) } =
  \begin{dcases*}
    1 &if $\bm{k} = \bm{l}$\\
    0 &otherwise
  \end{dcases*}.
\end{gather*}
In this work, we have opted for an isotropic sparse tensorization of nested
Clenshaw--Curtis quadrature rules, where the level $l\in \mathbb N$ of the
sparse grid is the only parameter controlling the quality of the approximation.
As $l$ increases, both the number of nodes $M$ and the multi-index set
$\mathcal{K}$ increase.

\subsection{Sensitivity analysis}

The sensitivity analysis consists in the evaluation of the different
contributions of the input parameters on the variance of the solution. This is achieved by the computation of the Sobol indices, defined as the ratio between
the partial variance corresponding to the input parameter under investigation
$\xi_i$, and the total variance of the quantity of interest $X$
\begin{gather}\label{eq:sobol1}
    S_{1,i} := \frac{\Var(\mathbb{E}(X(\bm{\xi})|\xi_i))}{\Var(X(\bm{\xi}))}, \qquad \forall 1\le i\le N,
\end{gather}
where $\mathbb{E}(X(\bm{\xi})|\xi_i)$  denotes the conditional expected value of
$X$ given $\xi_i$. The indices $(S_{1,i})_{1\le i\le N}$ in \eqref{eq:sobol1}
are known as principal or first-order Sobol indices and measure the individual
contribution of the coefficient $\xi_i$ to the variance.

Higher-order Sobol indices measure the effect of the concurrent variation of
more variables. For instance, second-order Sobol indices read
\begin{gather*}
    S_{i,j} := \frac{\Var(\mathbb{E}(X(\bm{\xi})|\xi_i, \xi_j))}{\Var(X(\bm{\xi}))} - S_{1,i} - S_{1,j}, \qquad \forall 1\le i<j\le N.
\end{gather*}
Finally, the total Sobol index is obtained as the sum of the indices involving parameter $\xi_i$
\begin{gather*}
S_{T_i} : = 1 - \frac{\Var(\Esp(X(\bm{\xi})|\bm{\xi}_{\setminus i}))}{\Var(X(\bm{\xi}))},
\end{gather*}
where the vector $\bm{\xi}_{\setminus i}=\left( \xi_{j\ne i}\right)$ contains
all uncertain variables except $\xi_i$.

The computation of the Sobol indices results directly from the PC approximation
\eqref{eq:PCdef} of the variable of interest. Indeed, the partial variances can
be explicitly expressed as functions of the spectral modes similarly to the
statistical quantities in \eqref{eq:Esp_Var}. We refer the reader
to~\cite{Crestaux.Le-Maitre.ea:09} for all the details.

\section{Numerical examples}\label{sec:numerical_examples}

In this section, we present two test cases to validate the proposed approach. In
both cases we show the reference numerical solution and discuss the uncertainty
quantification related to three parameters affected by uncertainty.

Due to the complexity of the problem, all data in the two cases are the same
except for the number of fractures in the network.
In the first case, presented in
Subsection \ref{subsec:example_single_fracture}, a single fracture touching one boundary is
considered, while in Subsection
\ref{subsec:example_multiple_fracture}, ten intersecting fractures are considered.

The data for the Darcy problem, the heat equation, and the precipitattion-dissolution process are given in Table \ref{tab:data_case_adv}.
Furthermore, we assume that the following three
parameters are affected by uncertainty and uniformly distributed with mean and
variance given by
\begin{gather*}
    \eta_\gamma \sim \mathcal{U}(2, 0.17)
    \qquad
    E \sim \mathcal{U}(4, 0.35)
    \qquad
    \theta^{\rm inflow} \sim \mathcal{U}(1.5, 0.11)
\end{gather*}
where $\theta^{\rm inflow}$ denotes the inflow temperature on the bottom boundary.
Being the construction of the sparse grids dependent only on the chosen level
and the number of uncertain parameters, the number of
simulations needed to construct the PC expansion are: 31 runs for level 2, 111 runs for level 3, 351 for level 4, 1023 for and level 5.
Once constructed, the evaluation of the PC expansion takes a small fraction of the time used by the full order model to evaluate the solutions for different times and for different values of
the uncertain parameters. Hence, we will use the PC expansion as a surrogate model and evaluate its performances and accuracy properties.
\begin{table}[tbp]
    \centering
    \begin{tabular}{|c|c|c|c|c|}
        \hline
        $\mu=1$ & $g=0$ & $\rho_w=1$ & $k_{\Omega,0} = 1$ & $\phi_{\Omega,0} = 0.2$ \\\hline
        $d_\Omega = 0.1$ & $c_w=1$  & $c_s=1$  & $\Lambda_w=1$ &  $\Lambda_s=0.1$  \\\hline
        $j =0$ & $\theta_0=0$ &
        $\lambda=10 \exp(-\frac{E}{\theta})$ & $u_0 = 0$ & $w_0 = 0.3$ \\\hline
        $u_e=1$ &
        $p_{\partial\Omega}^{\rm inflow}=1$ & $p_{\partial\Omega}^{\rm outflow}=0$ & $p_{\partial\gamma}^{\rm inflow}=1$ & $p_{\partial\gamma}^{\rm outflow}=0$
        \\ \hline
        $\epsilon_{\gamma, 0}=10^{-2}$ &  $k_{\gamma, 0} = 10^2$ & $\kappa_{\gamma, 0} =10^2$  & $f_{\Omega}=0$ & $f_{\gamma}=0$\\ \hline
        $\theta_{\partial \Omega}^{\rm inflow} = 1.5$ &   $\theta_{\partial \gamma}^{\rm inflow} = 1.5$ &
        $u_{\partial \Omega}^{\rm inflow} = 2$ &  $u_{\partial \gamma}^{\rm inflow} = 2$ & $j_{\gamma} =0$ \\  \hline
        $d_\gamma = 0.1$ & $\delta_\gamma = 0.1 $ &  $\eta=0.5$  & & \\ \hline
    \end{tabular}
    \caption{Common data for the advection-reaction problem, examples in Section
    \ref{sec:numerical_examples}. For $\bullet\in\{ k_\Omega; \phi_\Omega; \Theta; u; w; \varepsilon_\gamma; k_\gamma; \kappa_\gamma \}$, the notation $(\bullet)_0$ is used for the reference value of the quantity $\bullet$ prescribed as initial condition.}%
    \label{tab:data_case_adv}
\end{table}
In both cases, we will discuss the convergence properties of the
PC expansion by increasing the level of the sparse grid considered, the analysis of the
Sobol indices, conditioned variances and covariances for selected solutions and
finally the probability distribution functions.

The simulations are developed with
the library PorePy, a simulation tool for fractured and deformable porous media
written in Python, see \cite{Keilegavlen2019}.

\subsection{Single fracture network}\label{subsec:example_single_fracture}

\begin{figure}
    \centering
    \resizebox{0.33\textwidth}{!}{\fontsize{1cm}{2cm}\selectfont
\begingroup%
  \makeatletter%
  \providecommand\color[2][]{%
    \errmessage{(Inkscape) Color is used for the text in Inkscape, but the package 'color.sty' is not loaded}%
    \renewcommand\color[2][]{}%
  }%
  \providecommand\transparent[1]{%
    \errmessage{(Inkscape) Transparency is used (non-zero) for the text in Inkscape, but the package 'transparent.sty' is not loaded}%
    \renewcommand\transparent[1]{}%
  }%
  \providecommand\rotatebox[2]{#2}%
  \newcommand*\fsize{\dimexpr\f@size pt\relax}%
  \newcommand*\lineheight[1]{\fontsize{\fsize}{#1\fsize}\selectfont}%
  \ifx\svgwidth\undefined%
    \setlength{\unitlength}{311.96932983bp}%
    \ifx\svgscale\undefined%
      \relax%
    \else%
      \setlength{\unitlength}{\unitlength * \real{\svgscale}}%
    \fi%
  \else%
    \setlength{\unitlength}{\svgwidth}%
  \fi%
  \global\let\svgwidth\undefined%
  \global\let\svgscale\undefined%
  \makeatother%
  \begin{picture}(1,0.98417801)%
    \lineheight{1}%
    \setlength\tabcolsep{0pt}%
    \put(0,0){\includegraphics[width=\unitlength,page=1]{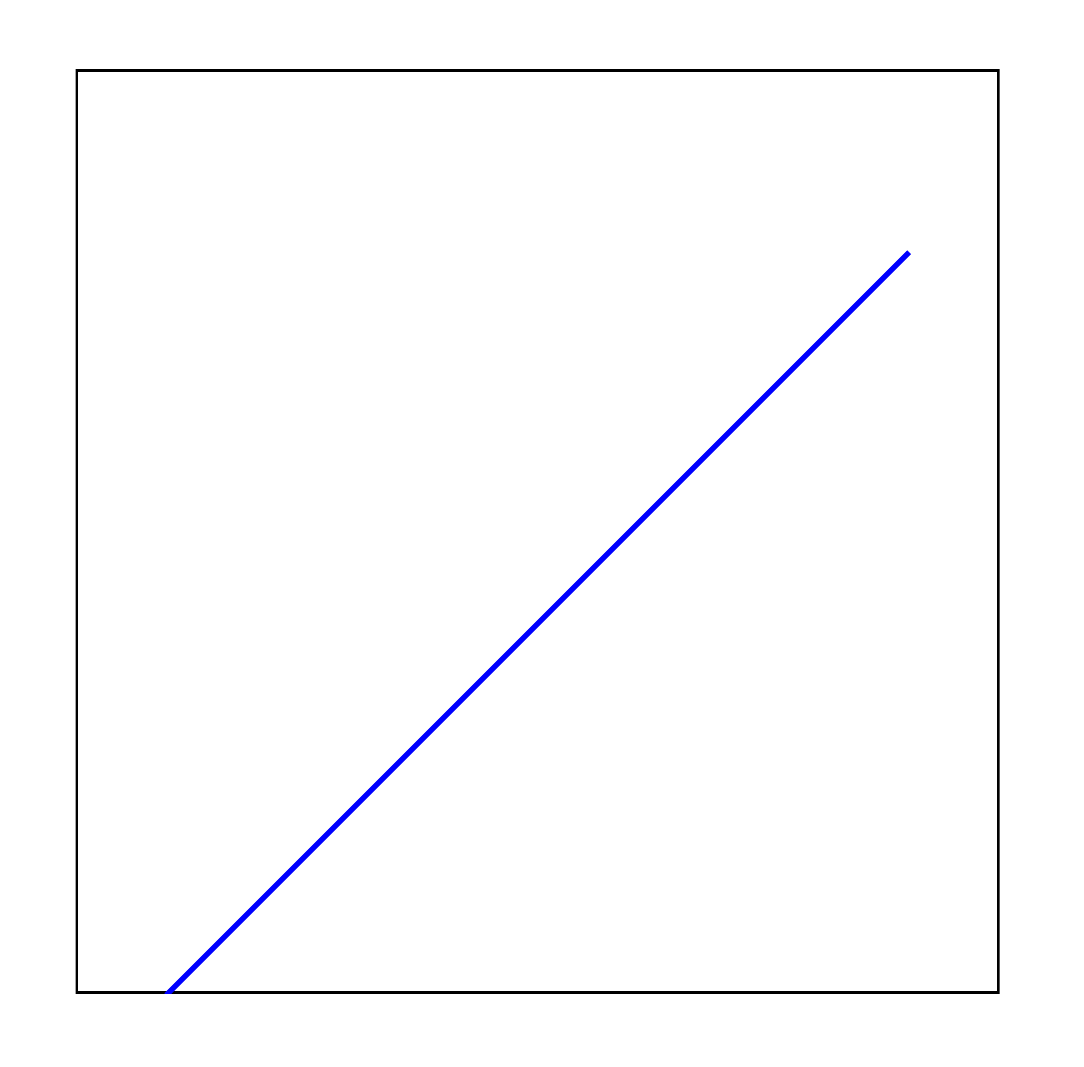}}%
    \put(0.40318924,0.93547031){\color[rgb]{0,0,0}\makebox(0,0)[lt]{\lineheight{1.25}\smash{\begin{tabular}[t]{l}out-flow\end{tabular}}}}%
    \put(0.40068498,0.0009078){\color[rgb]{0,0,0}\makebox(0,0)[lt]{\lineheight{1.25}\smash{\begin{tabular}[t]{l}in-flow\end{tabular}}}}%
    \put(0.9990922,0.40283447){\color[rgb]{0,0,0}\rotatebox{90}{\makebox(0,0)[lt]{\lineheight{1.25}\smash{\begin{tabular}[t]{l}no-flow\end{tabular}}}}}%
    \put(0.00090779,0.64233492){\color[rgb]{0,0,0}\rotatebox{-90}{\makebox(0,0)[lt]{\lineheight{1.25}\smash{\begin{tabular}[t]{l}no-flow\end{tabular}}}}}%
    \put(0.50403155,0.36449245){\color[rgb]{0,0,0}\makebox(0,0)[lt]{\lineheight{1.25}\smash{\begin{tabular}[t]{l}$\gamma$\end{tabular}}}}%
    \put(0.12781315,0.81059003){\color[rgb]{0,0,0}\makebox(0,0)[lt]{\lineheight{1.25}\smash{\begin{tabular}[t]{l}$\Omega$\end{tabular}}}}%
  \end{picture}%
\endgroup%
}
    \caption{Domain $\Omega$ and fracture $\gamma$ for the example of Subsection
    \ref{subsec:example_single_fracture}.}%
    \label{fig:example_single_fracture_domain}
\end{figure}
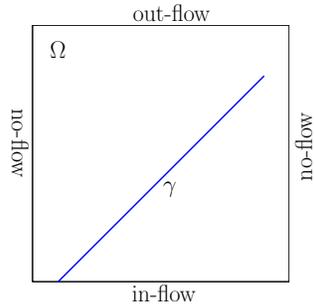

Let us consider a domain $\Omega=(0, 1)^2$ with a single immersed fracture
defined by the following vertices: $(0.1, 0)$ and $(0.9, 0.8)$. The fracture thus
touches the bottom boundary of $\Omega$ as depicted in Figure \ref{fig:example_single_fracture_domain}.
Data and uncertain parameters are
reported in the beginning of Section \ref{sec:numerical_examples}.
We point out that the fracture is permeable at the beginning of the simulation and due to the solute precipitation its aperture diminishes in time.
As a result the effective fracture permeability decreases and until the fracture behaves as a
barrier and not any more as a preferential path.

In the following parts, we detail some aspects related to the uncertainty
quantification analysis. In Subsection \ref{subsec:case1_conv} a convergence study is
carried out, in Subsection \ref{subsec:case1_var} we discuss the variances and covariances of the solutions and in Subsection \ref{subsec:case1_pdf} we
introduce the computed probability distribution functions of some of the components of the solutions along the fracture.

\begin{figure}[tbp]
    \centering
    \includegraphics[width=0.32\textwidth]{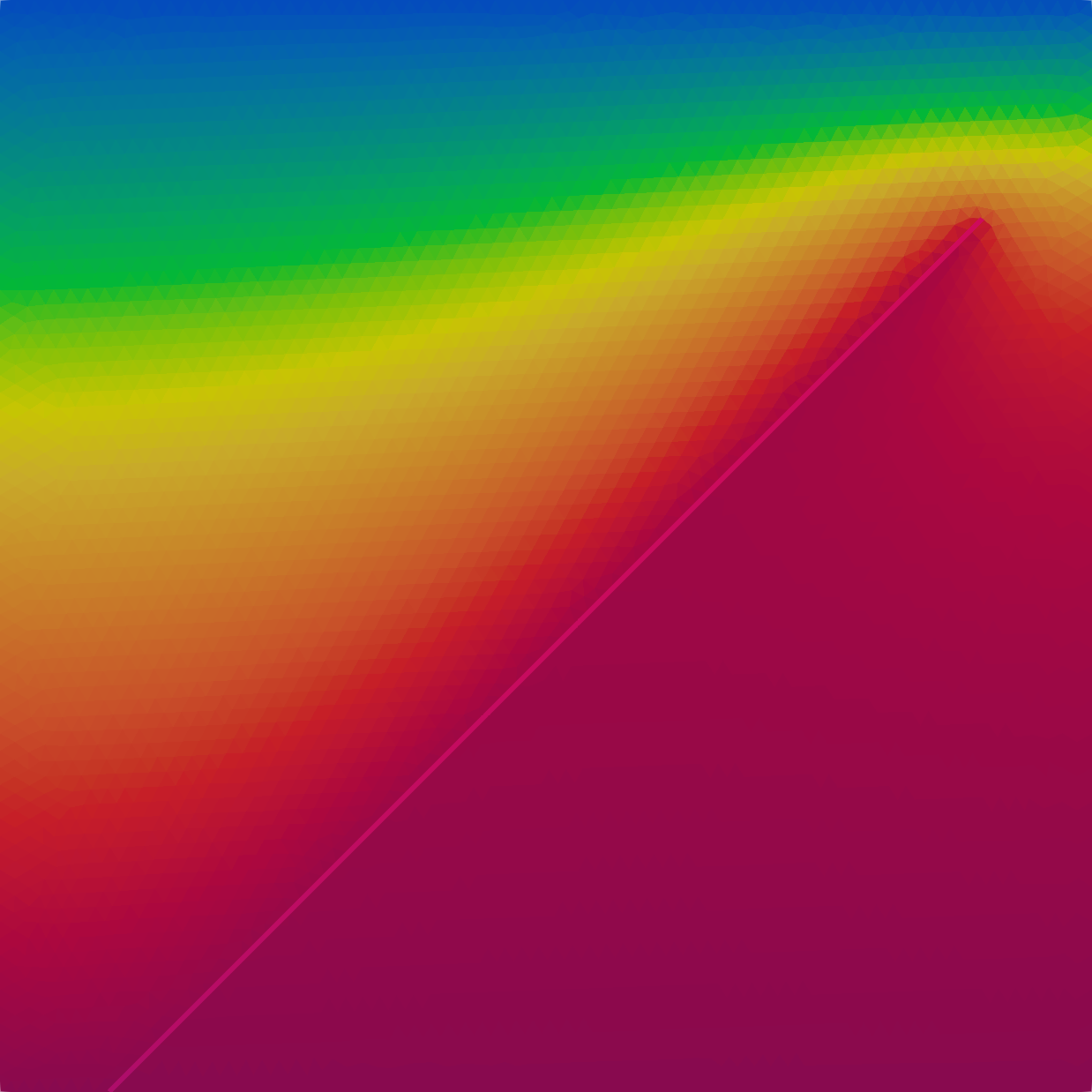}%
    \hspace*{0.02\textwidth}%
    \includegraphics[width=0.32\textwidth]{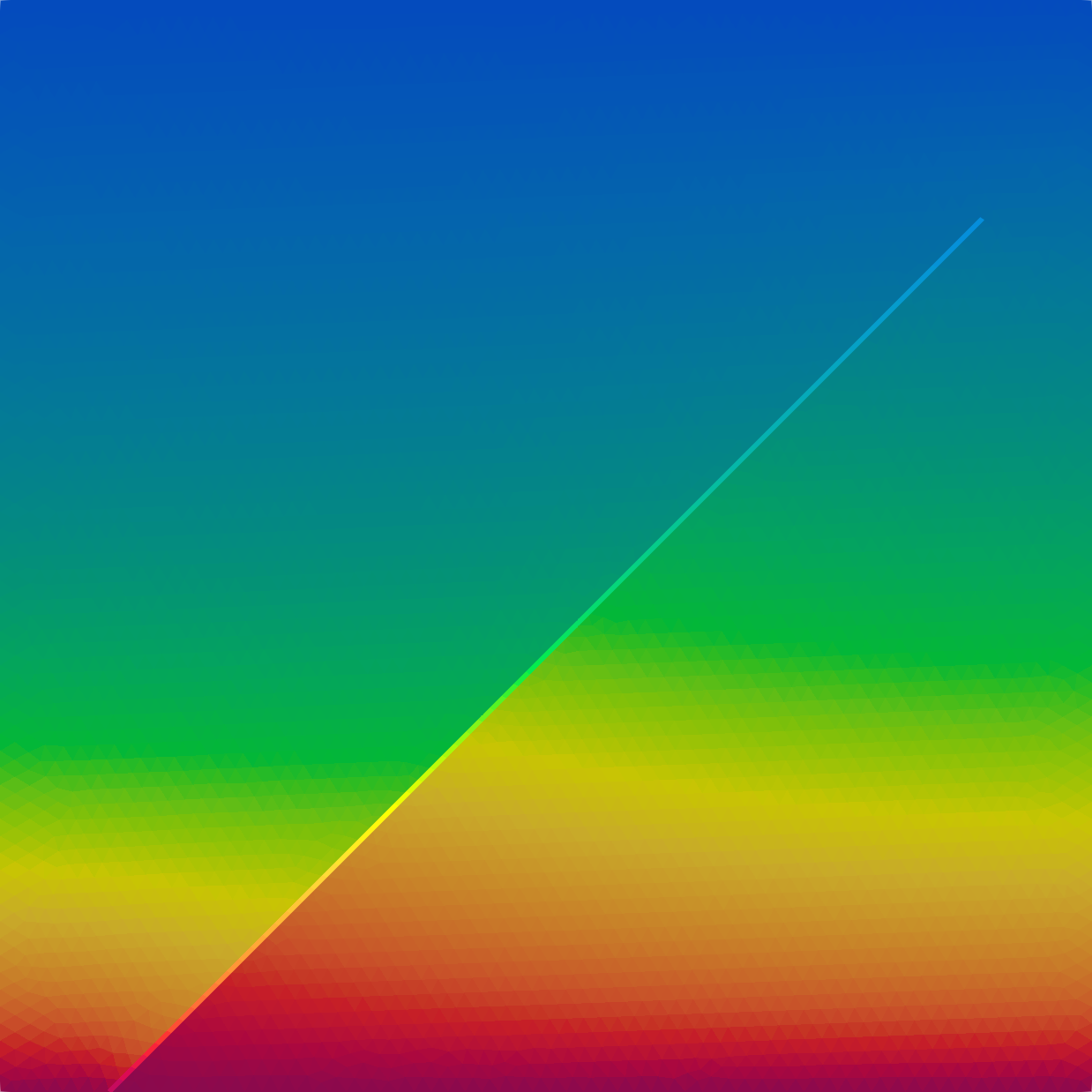}%
    \hspace*{0.02\textwidth}%
    \includegraphics[width=0.32\textwidth]{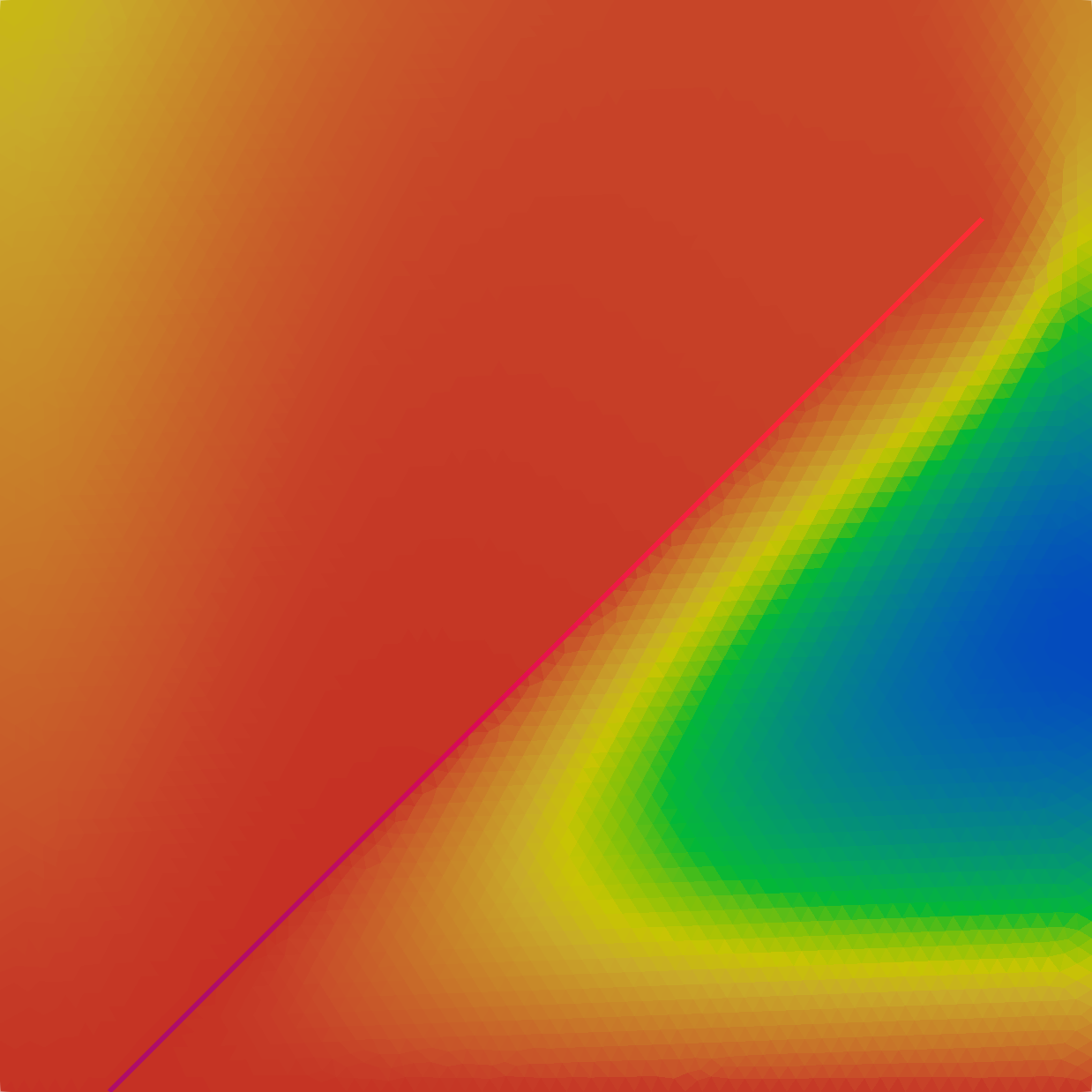}
    \hspace*{0.18\textwidth}%
    \includegraphics[width=0.32\textwidth]{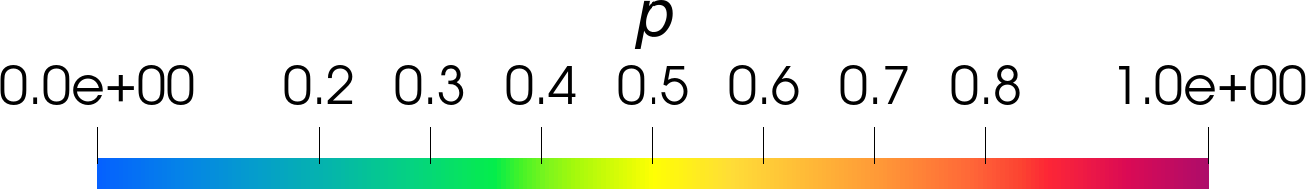}%
    \hspace*{0.18\textwidth}%
    \includegraphics[width=0.32\textwidth]{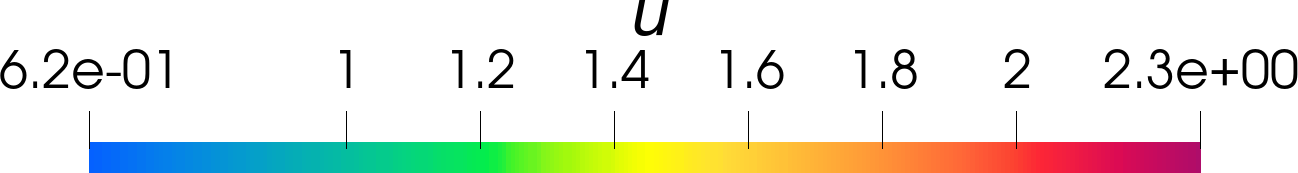}
    \caption{Reference solutions with mean value of the uncertain parameters.
    On the left pressure at time $t=0.1T$ and on the centre for $t=T$, on the
    right the solute for $t=0.1T$.  Test case of
    Subsection \ref{subsec:example_single_fracture}.}
    \label{fig:case1_solution}
\end{figure}

The reference solution, corresponding to the average input parameters, is reported in Figure
\ref{fig:case1_solution}, where it is possible to notice the variation of the pressure
distribution over time due to the sealing of the fractures and the transport of
the solute when the fractures are still highly permeable.

\subsubsection{Convergence}\label{subsec:case1_conv}

In this part, we discuss the convergence and accuracy properties of the
surrogate model built with the PC expansion. In fact, the latter can be used to make fast
simulations without the need of running the full order model. Since we are dealing with a time dependent problem, we analyse the PC expansion for two different simulation times: after few time steps $(t=t_1=0.1T)$ and at the end of the simulation $(t=T)$.

\begin{figure}[tbp]
    \centering
    \includegraphics[width=0.475\textwidth]{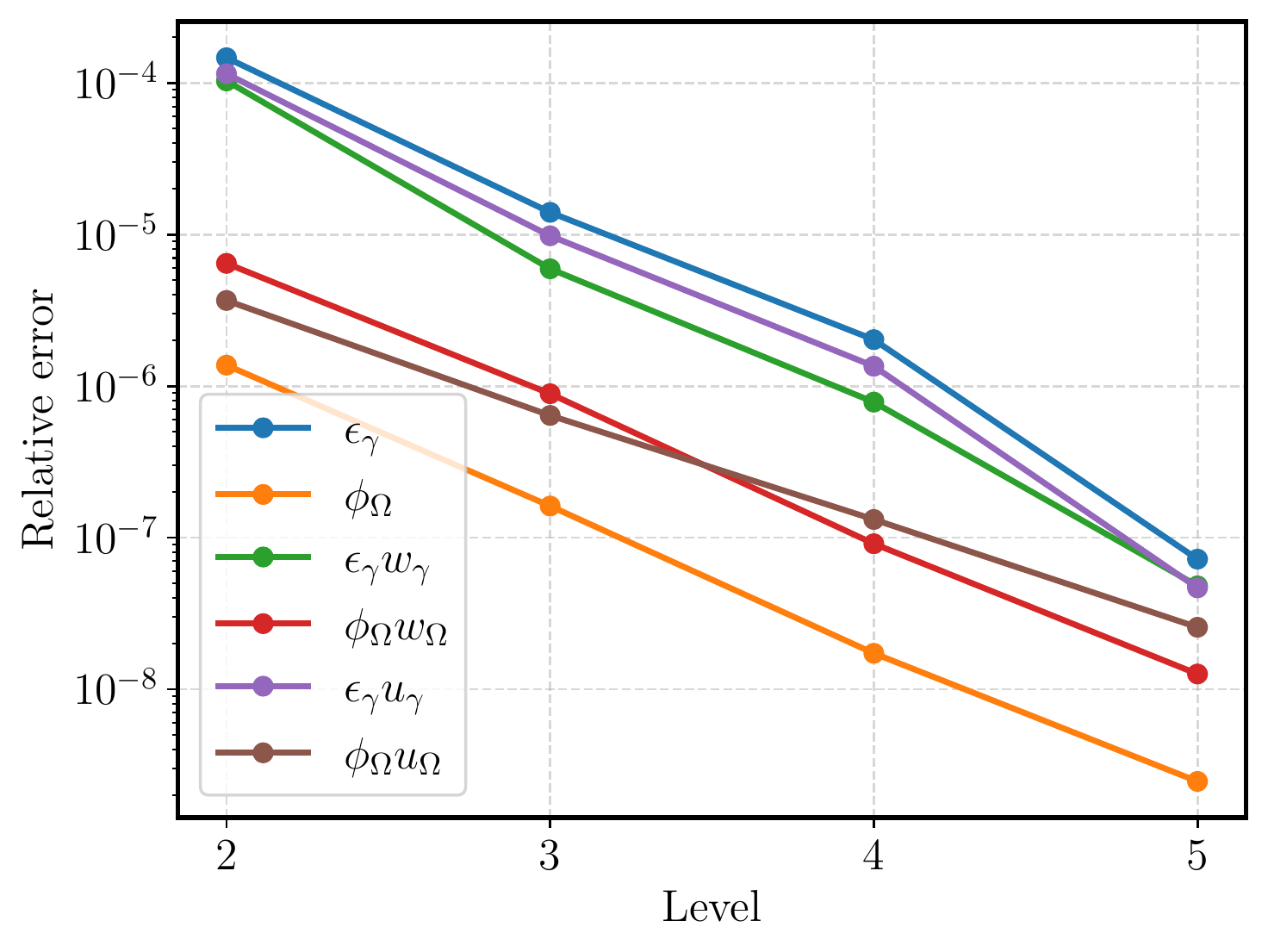}%
    \hspace*{0.05\textwidth}%
    \includegraphics[width=0.475\textwidth]{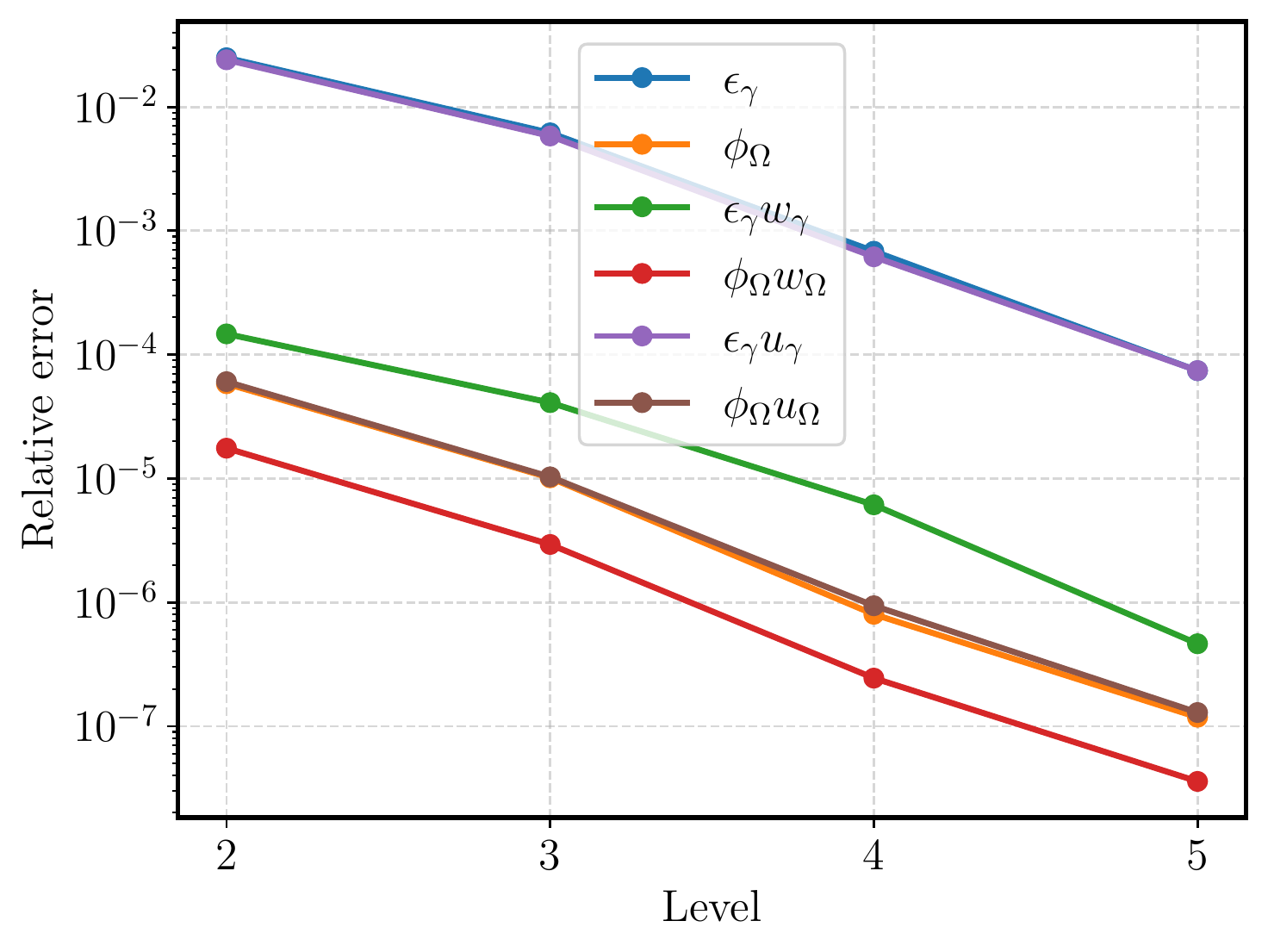}
    \caption{Error convergence for increasing level of the sparse grids, on the
    left at time $t=0.1T$ and on the right at final time $t=T$. Test case of
    Subsection \ref{subsec:case1_conv}.}
    \label{fig:case1_error}
\end{figure}

Figure \ref{fig:case1_error} presents, for both times, the error decay of the
computed solutions by increasing the level of the sparse grid. For a smooth
relation between uncertain parameters and the solutions, we expect exponential
decay of the error with respect to the level, which is the behaviour observed in
the figure. Additionally, for $t = t_1$ the error computed is much smaller compared to the
end of the simulation, showing a temporal dependence on the quality of the PC expansion.

\begin{figure}[tbp]
    \centering
    \includegraphics[width=0.32\textwidth]{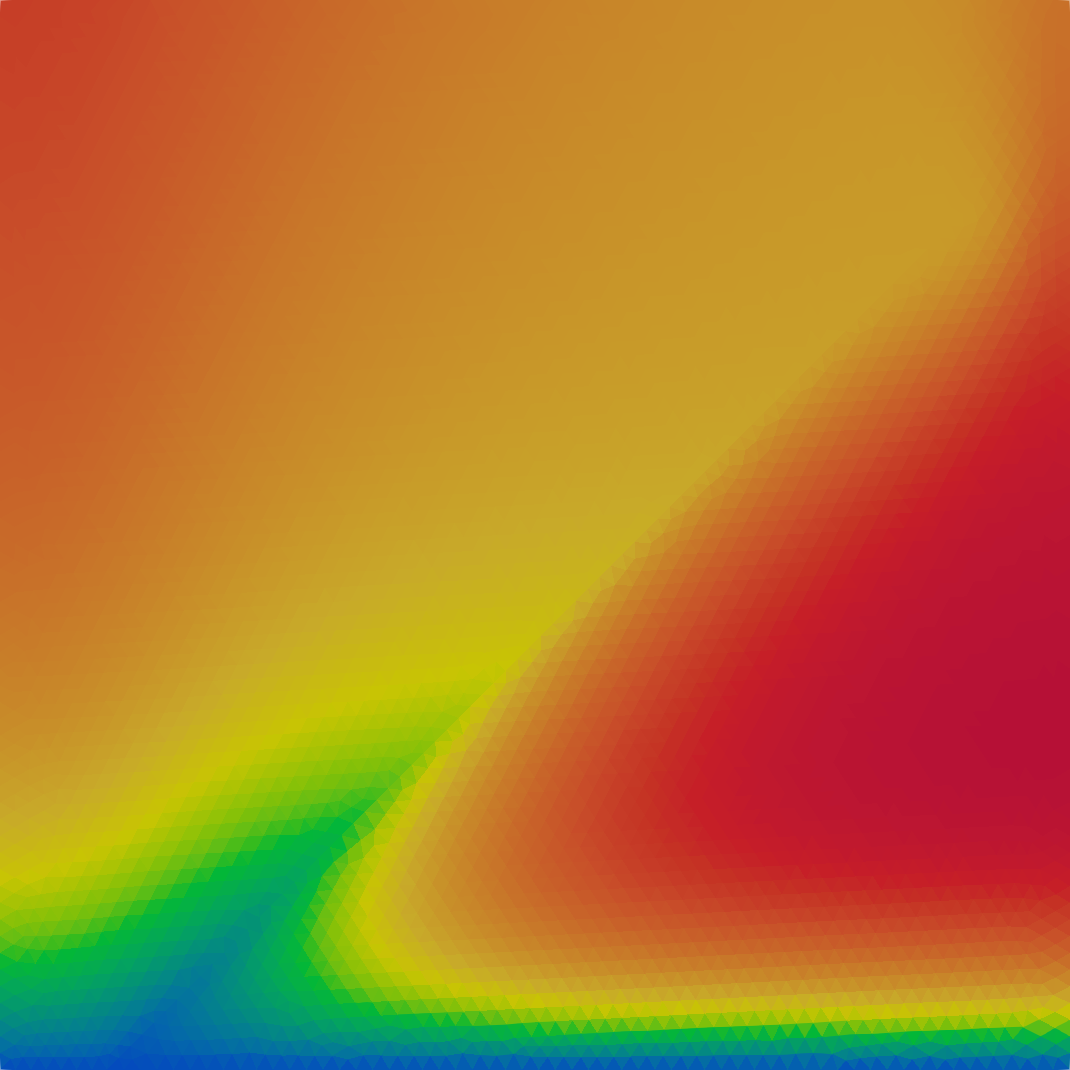}%
    \hspace*{0.02\textwidth}%
    \includegraphics[width=0.32\textwidth]{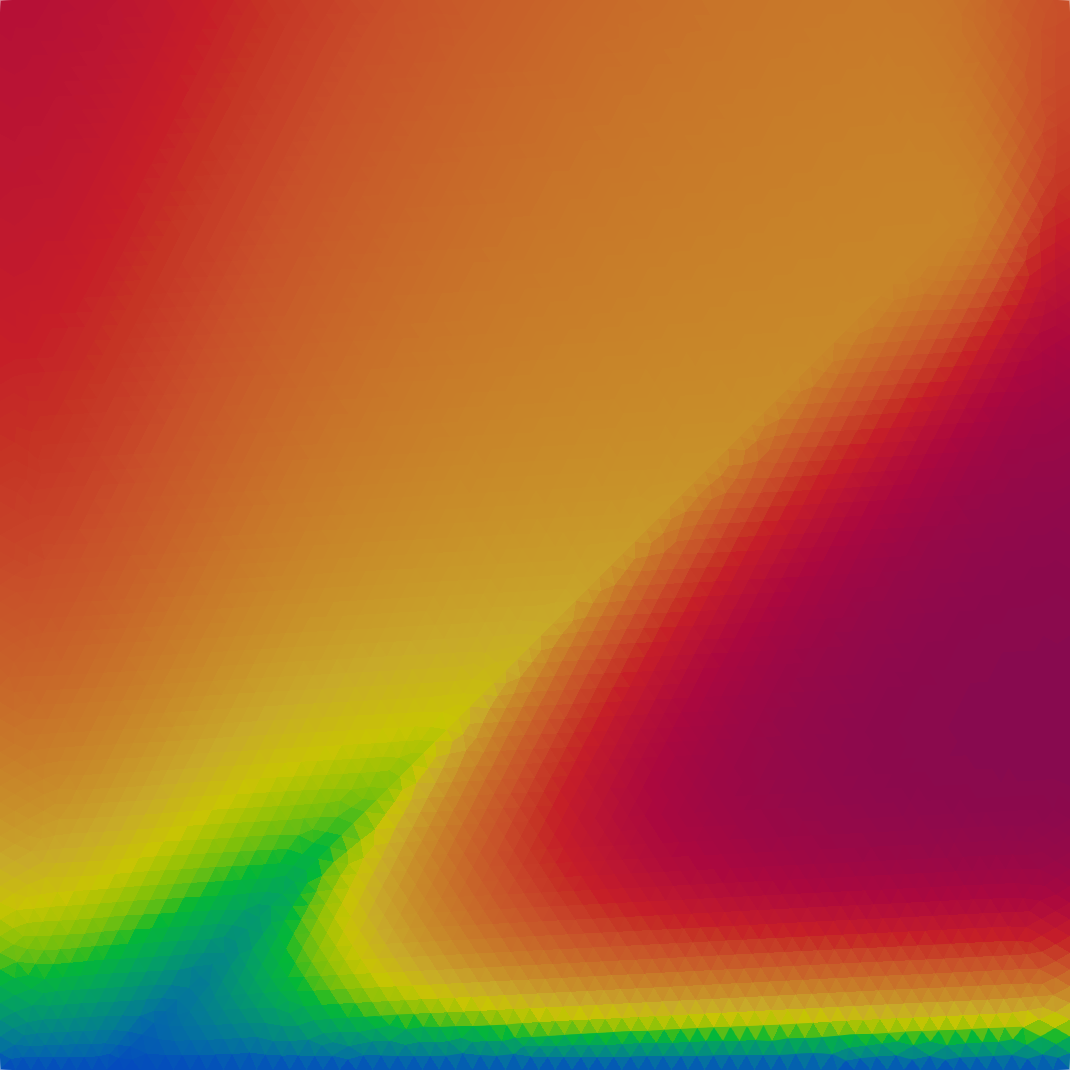}%
    \hspace*{0.02\textwidth}%
    \includegraphics[width=0.32\textwidth]{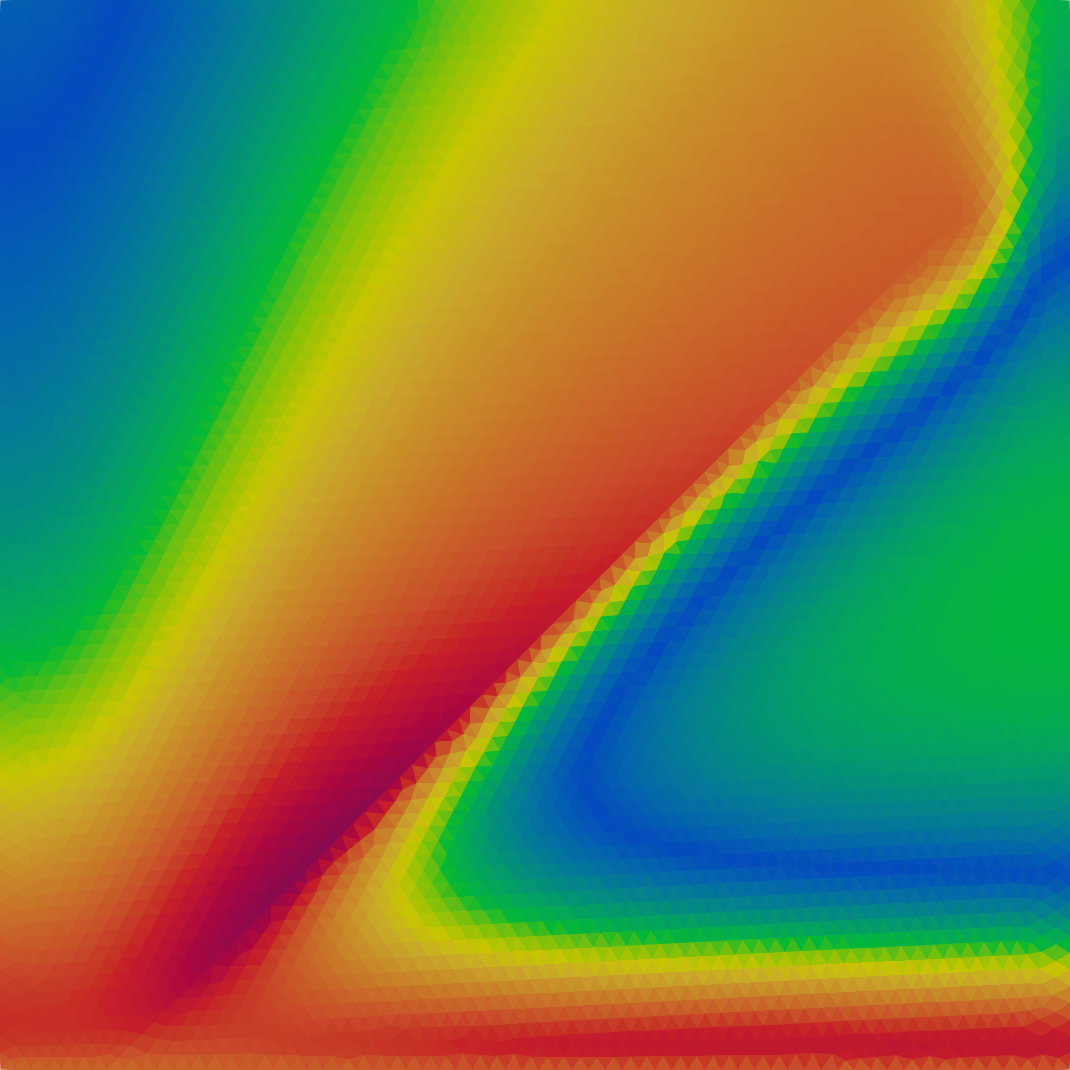}
    \includegraphics[width=0.32\textwidth]{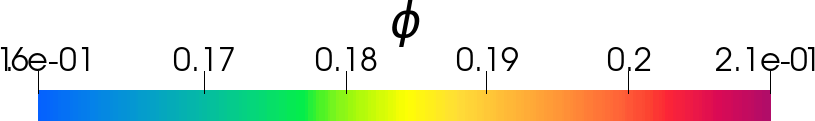}%
    \hspace*{0.02\textwidth}%
    \includegraphics[width=0.32\textwidth]{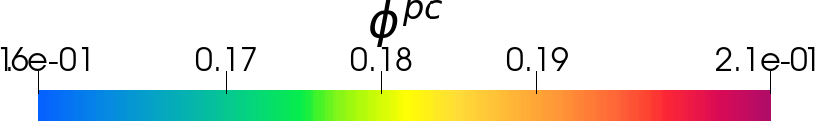}%
    \hspace*{0.02\textwidth}%
    \includegraphics[width=0.32\textwidth]{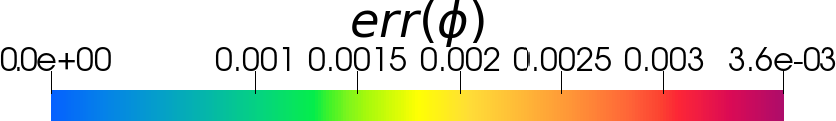}\\
    \vspace*{0.05\textwidth}
    \includegraphics[width=0.32\textwidth]{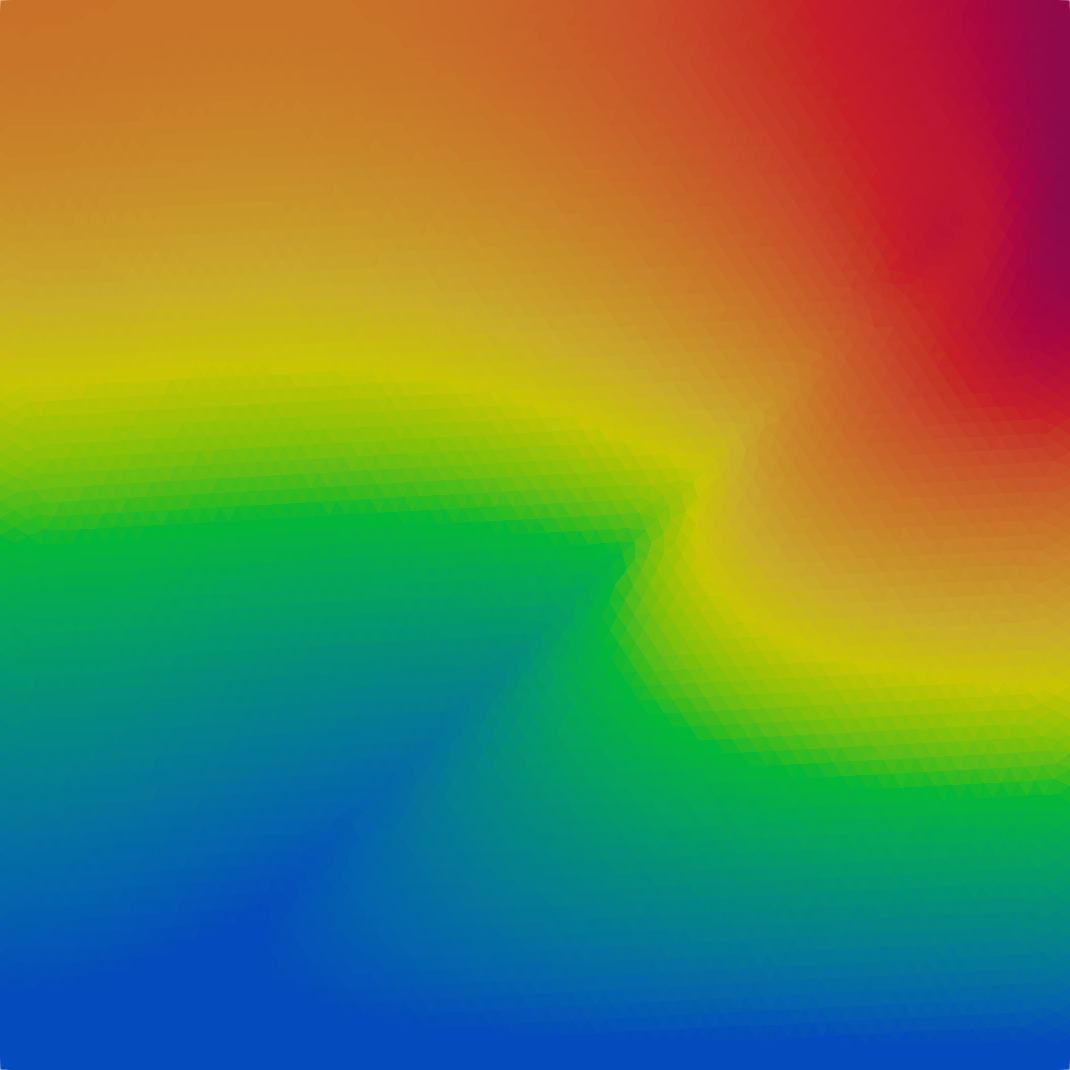}%
    \hspace*{0.02\textwidth}%
    \includegraphics[width=0.32\textwidth]{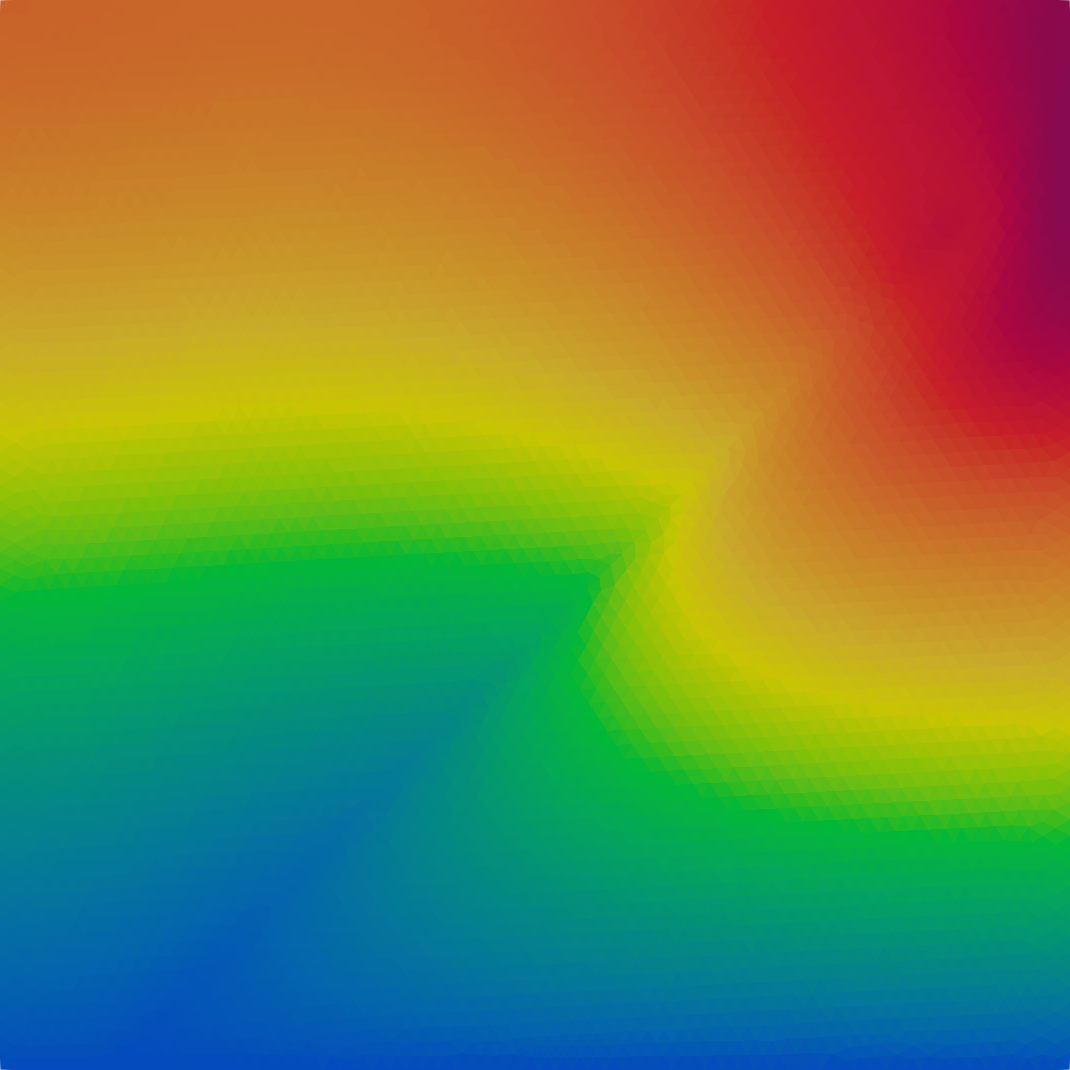}%
    \hspace*{0.02\textwidth}%
    \includegraphics[width=0.32\textwidth]{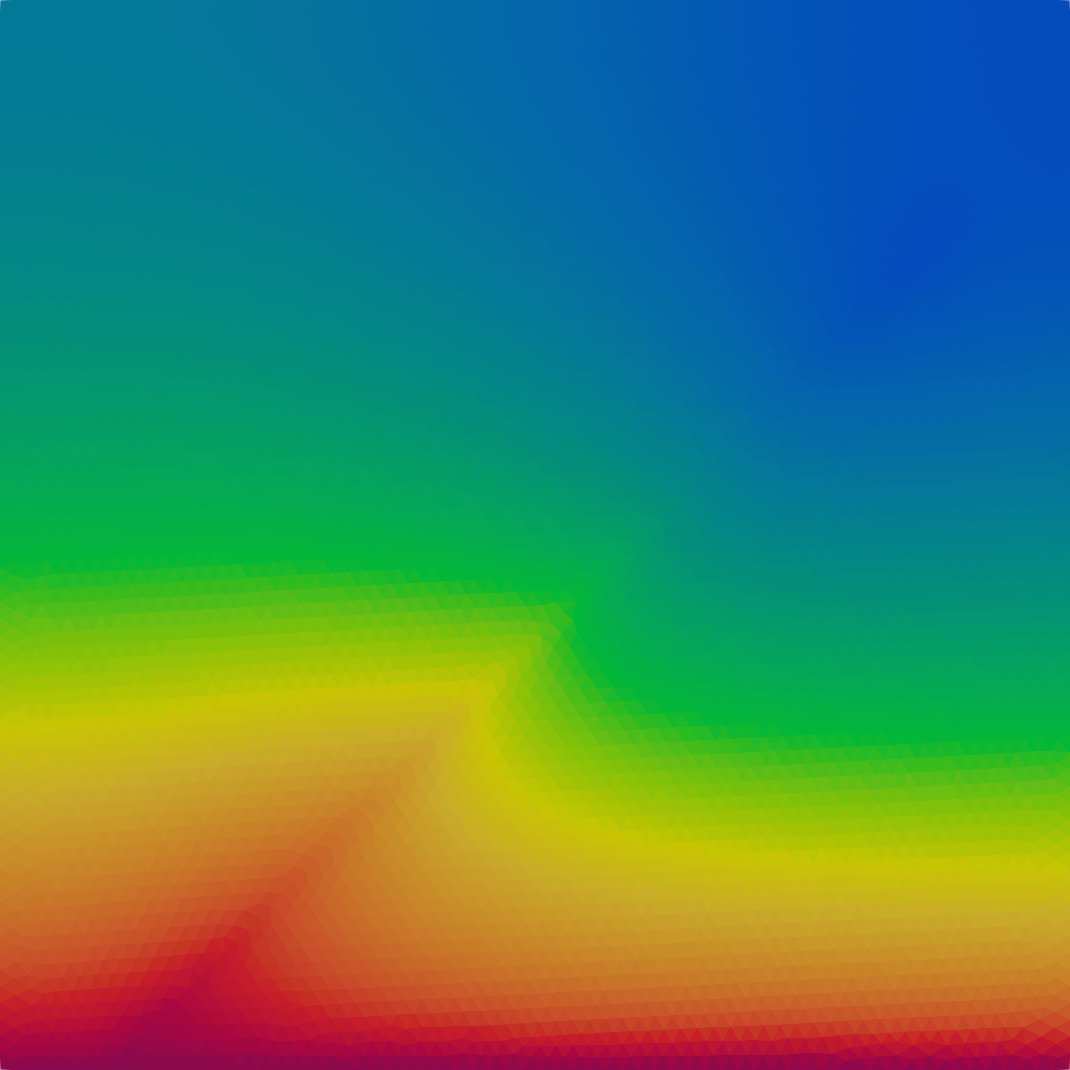}
    \includegraphics[width=0.32\textwidth]{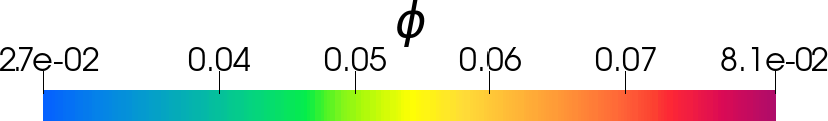}%
    \hspace*{0.02\textwidth}%
    \includegraphics[width=0.32\textwidth]{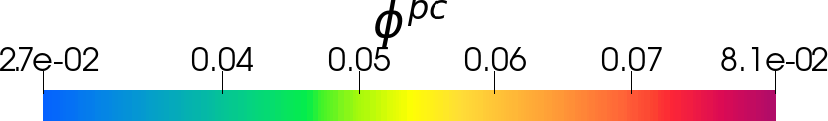}%
    \hspace*{0.02\textwidth}%
    \includegraphics[width=0.32\textwidth]{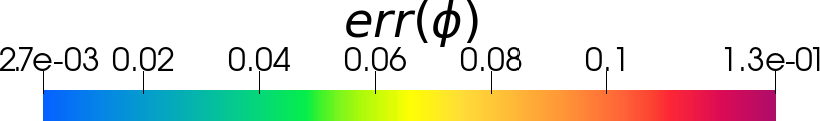}
    \caption{On the left porosity in the media computed with the original model
    and on the centre with the polynomial chaos expansion, on the right the
    error between them. On the
    top at time $t=0.1T$ and on the bottom at final time $t=T$. Test case of
    Subsection \ref{subsec:case1_conv}.}
    \label{fig:case1_porosity}
\end{figure}

In Figure \ref{fig:case1_porosity}, we compare the porosity $\phi$ computed with
the full order model with the one constructed by the PC expansion and the
corresponding relative error.
The two solutions are in good agreement for both times and the error is rather
low. As before, the latter is smaller at the beginning of the simulation and tends to
increase at the end, in particular near the inflow boundary at the bottom of
$\Omega$. Moreover, at the beginning of the simulation the fracture is highly
permeable and, consequently, we observe also a region in the proximity of $\gamma$ where the
error is higher due to stronger geochemical effects.

\begin{figure}[tbp]
    \centering
    \includegraphics[width=1\textwidth]{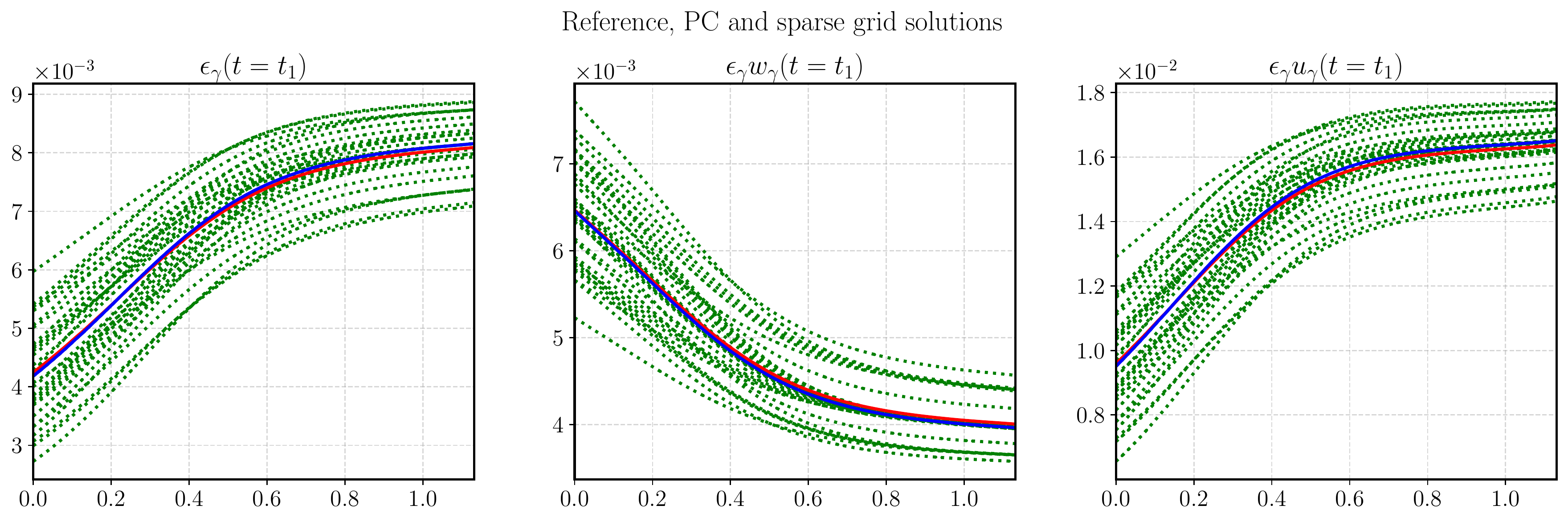}
    \includegraphics[width=1\textwidth]{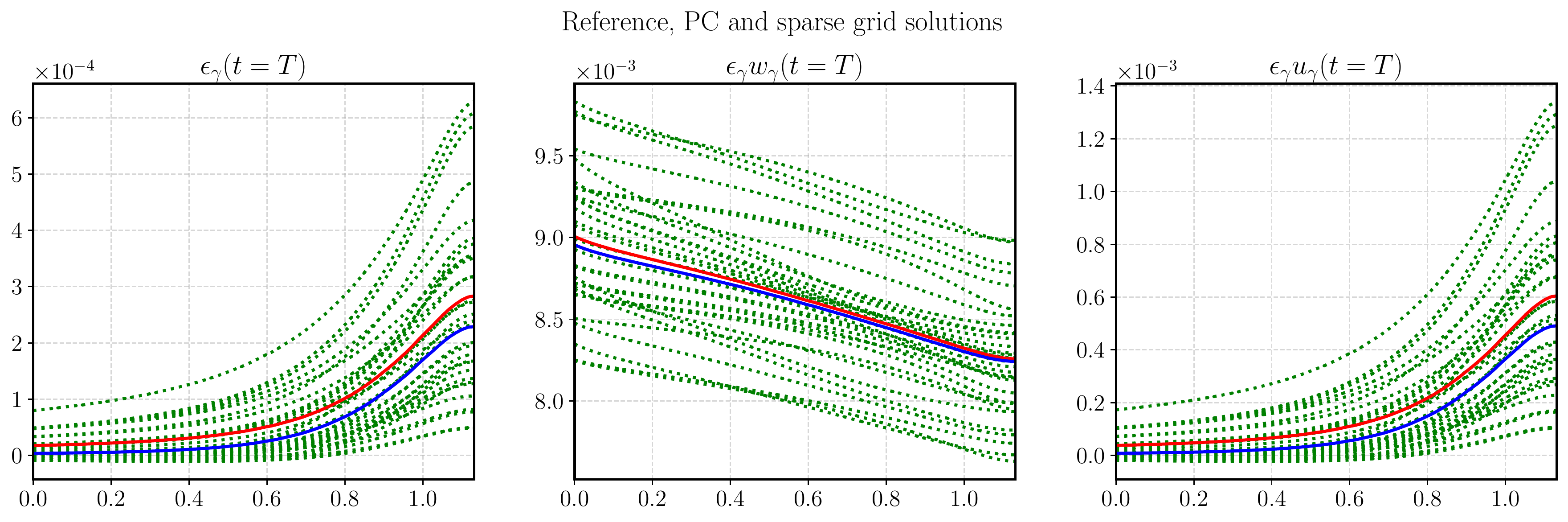}
    \caption{Solutions along $\gamma$: in blue computed with the original model,
    in red with the polynomial chaos expansion and in green computed by the original model
    on each sparse grid node  (sparse grid level 2). On the top at time $t=0.1T$
    and on the bottom at final time $t=T$. Test case of
    Subsection \ref{subsec:case1_conv}.}
    \label{fig:case1_1d}
\end{figure}

Finally, Figure \ref{fig:case1_1d} compares some of the variables in the
fractures computed with the full order model or reconstructed with the PC
expansion.
Also in this case, the quality of the latter is high and in good agreement with
the reference solution. Moreover, in Figure \ref{fig:case1_1d} the green dashed lines represent the solutions obtained for
each run to construct the PC expansion and can be useful to visualize the variance of the solution.

\subsubsection{Analysis of variance and correlations} \label{subsec:case1_var}

\begin{figure}[tbp]
    \centering
    \includegraphics[width=1\textwidth]{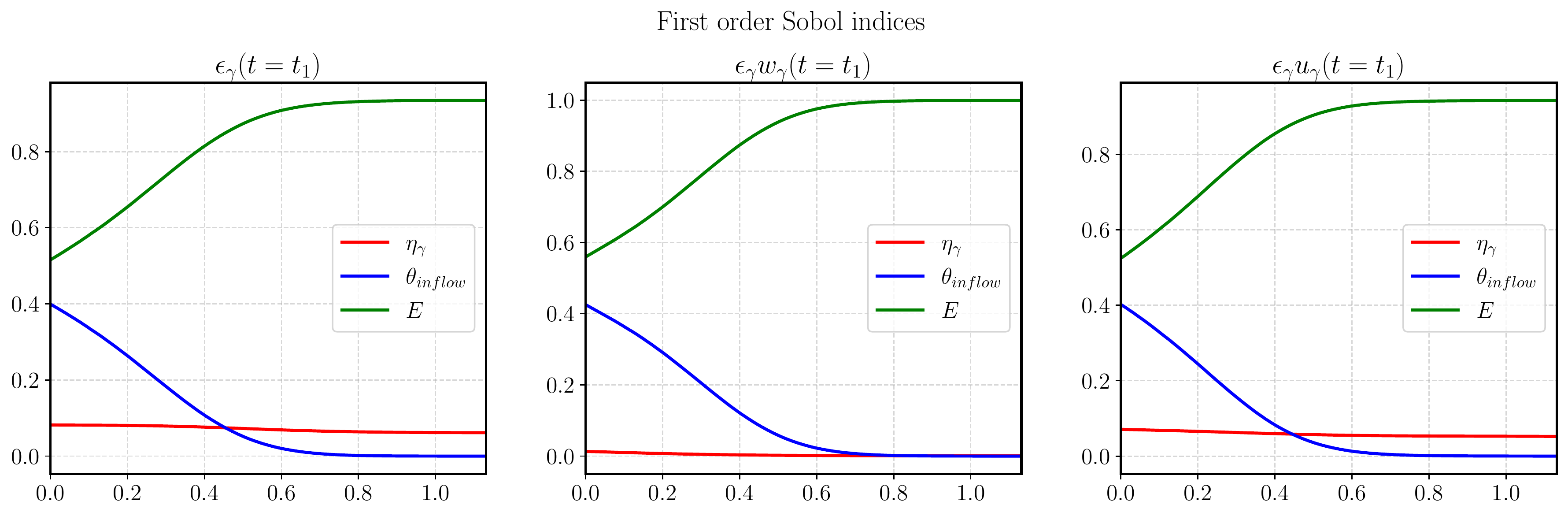}
    \includegraphics[width=1\textwidth]{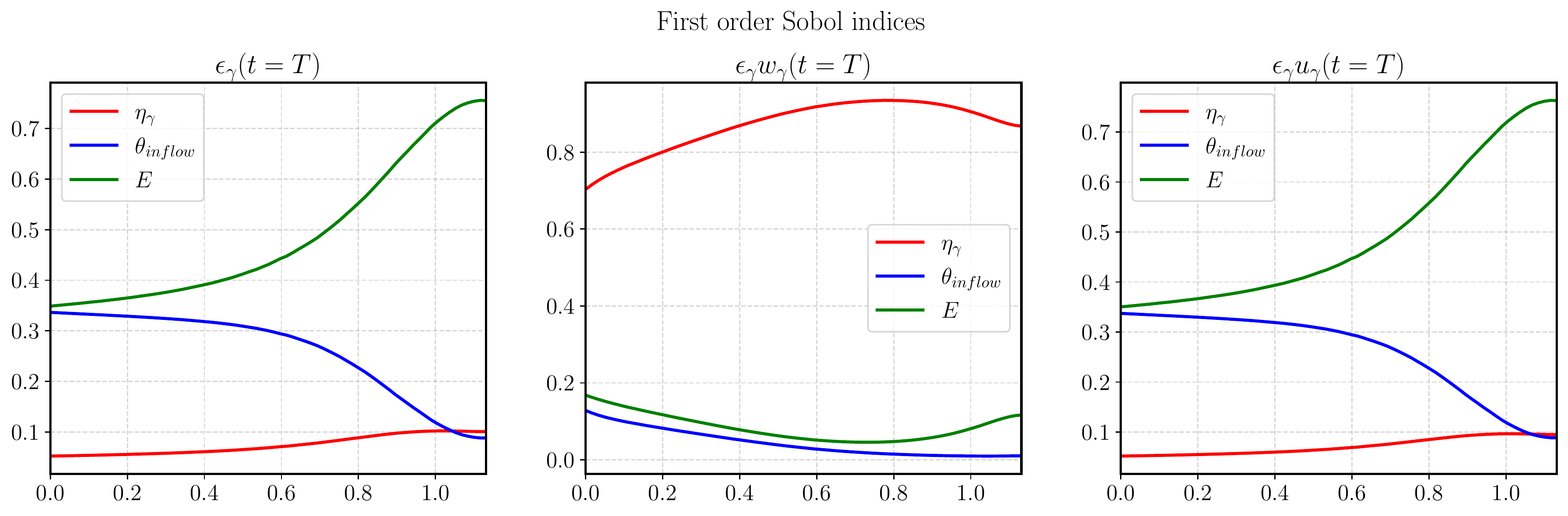}
    \caption{First order Sobol indexes for different unknowns along $\gamma$, on the
    top at time $t=t_1$ and on the bottom at final time $t=T$. The local
    coordinate of the fracture starts from the bottom boundary. Test case of
    Subsection \ref{subsec:case1_var}.}
    \label{fig:case1_sobol}
\end{figure}

An important factor is the impact on observed variables of the uncertain input data.
In Figure \ref{fig:case1_sobol} we report the Sobol indices for some variables
in the fracture for $t=t_1$ and $t=T$. We notice that the three variables
are influenced in a similar manner by the uncertain data, and the activation
energy is the most important factor, followed by the temperature at the inflow
boundary. We notice that while the high temperature front penetrates in the domain and in the fracture, the importance of the activation energy over the temperature
inflow tends to diminish and the latter becomes more important.
Note also that the Sobol indices of $\epsilon_\gamma w_\gamma$ at final time are different from
the other two considered variables since $\eta_\gamma$ dominates the induced uncertainty of the precipitate.

\begin{figure}[tbp]
    \centering
    \includegraphics[width=0.32\textwidth]{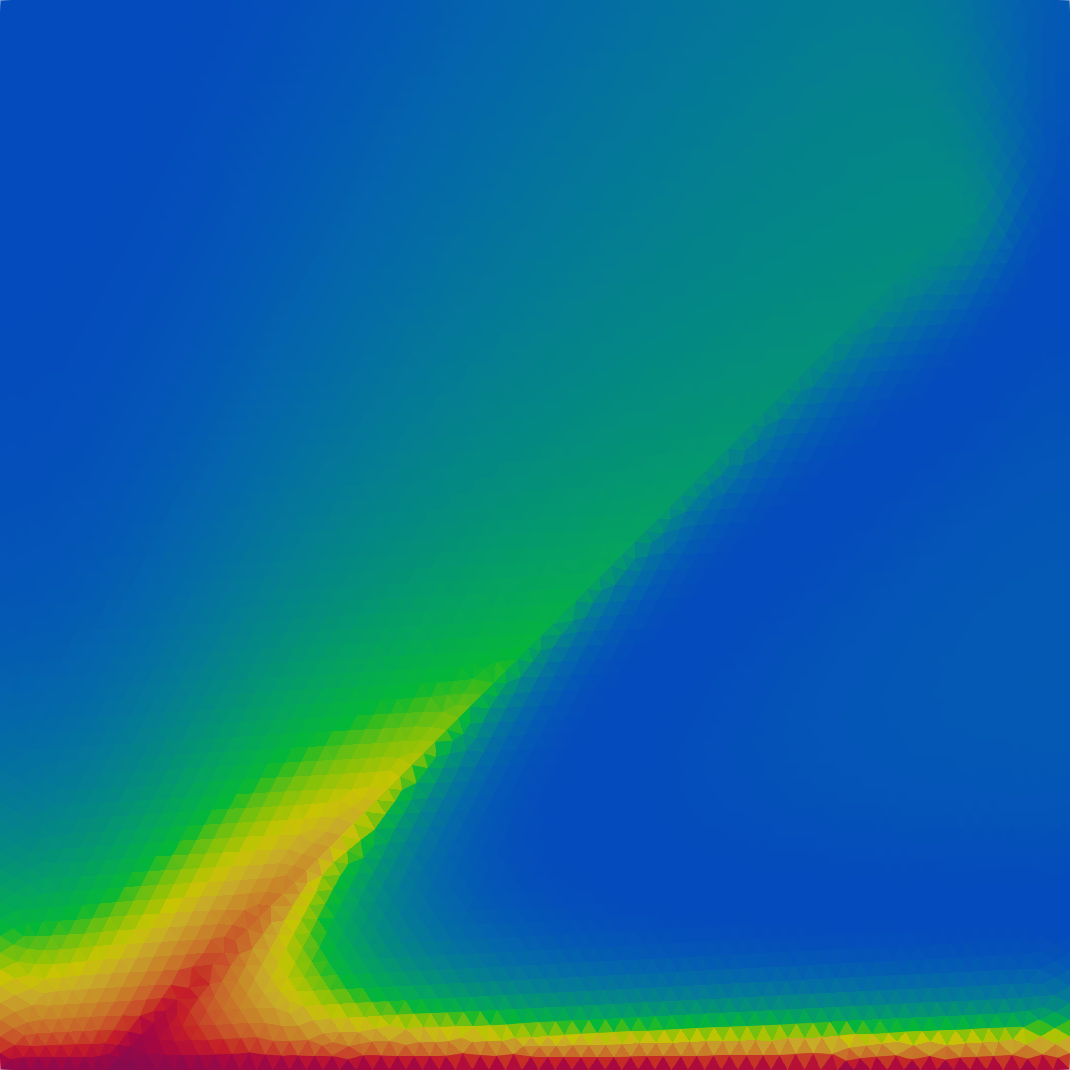}%
    \hspace*{0.02\textwidth}%
    \includegraphics[width=0.32\textwidth]{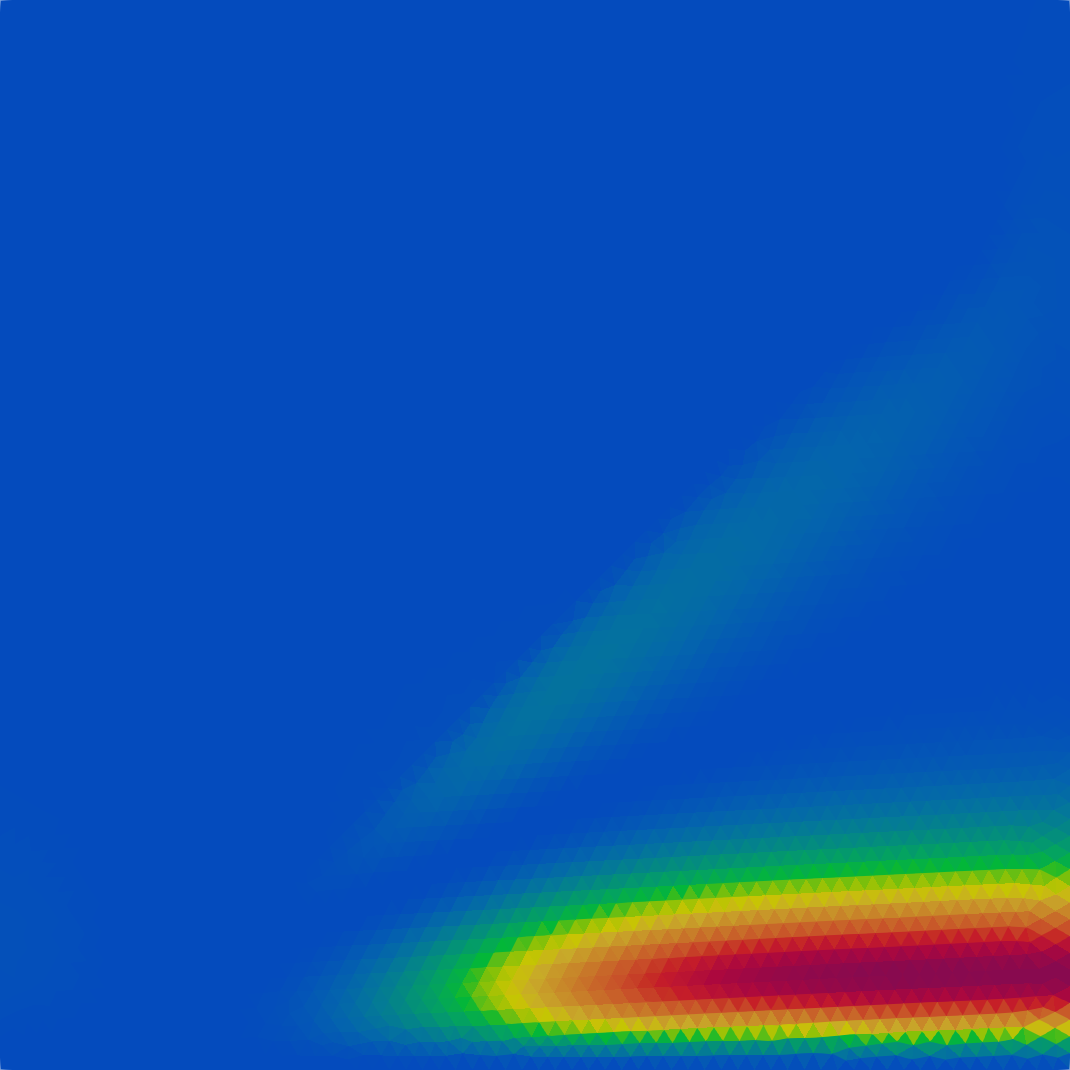}%
    \hspace*{0.02\textwidth}%
    \includegraphics[width=0.32\textwidth]{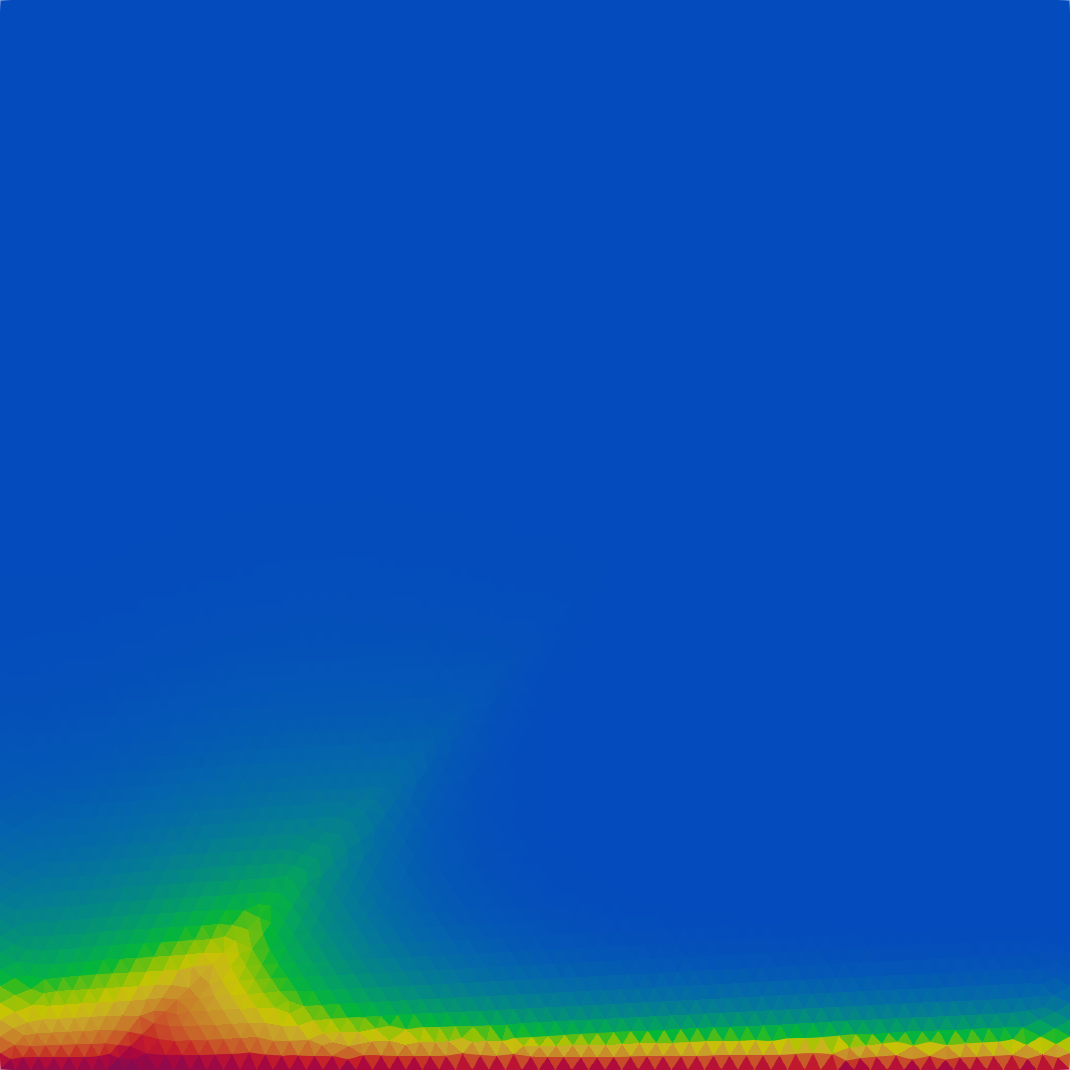}
    \includegraphics[width=0.32\textwidth]{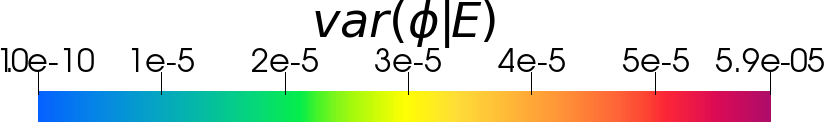}%
    \hspace*{0.02\textwidth}%
    \includegraphics[width=0.32\textwidth]{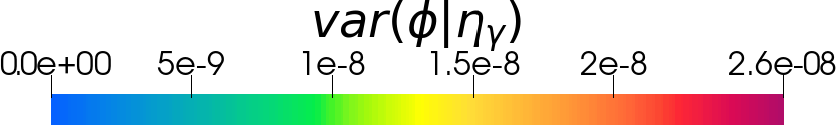}%
    \hspace*{0.02\textwidth}%
    \includegraphics[width=0.32\textwidth]{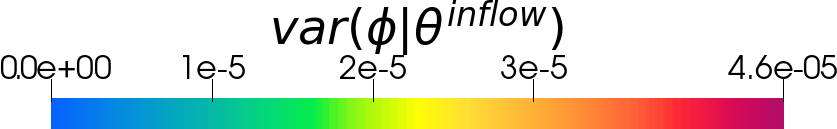}\\
    \vspace*{0.05\textwidth}
    \includegraphics[width=0.32\textwidth]{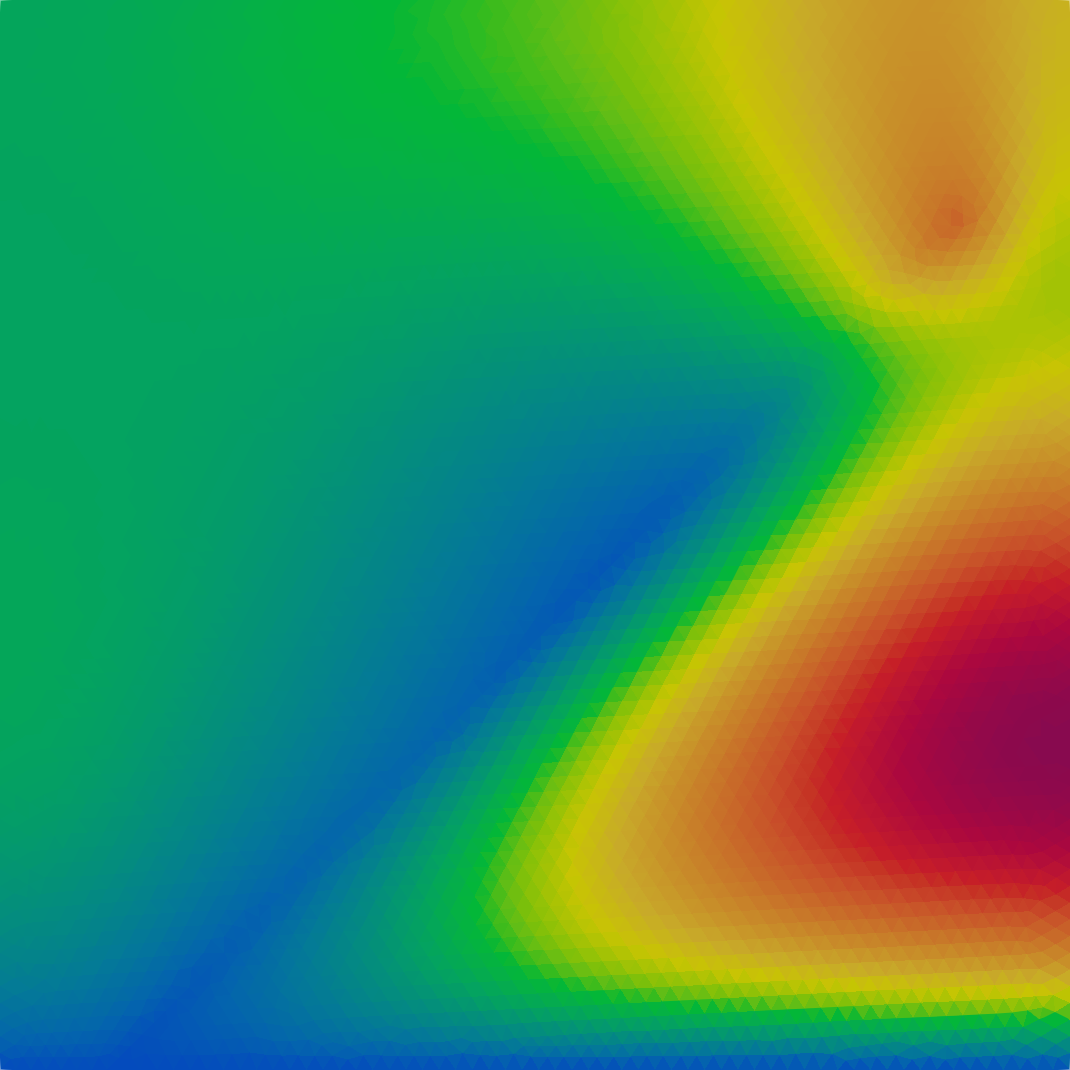}%
    \hspace*{0.02\textwidth}%
    \includegraphics[width=0.32\textwidth]{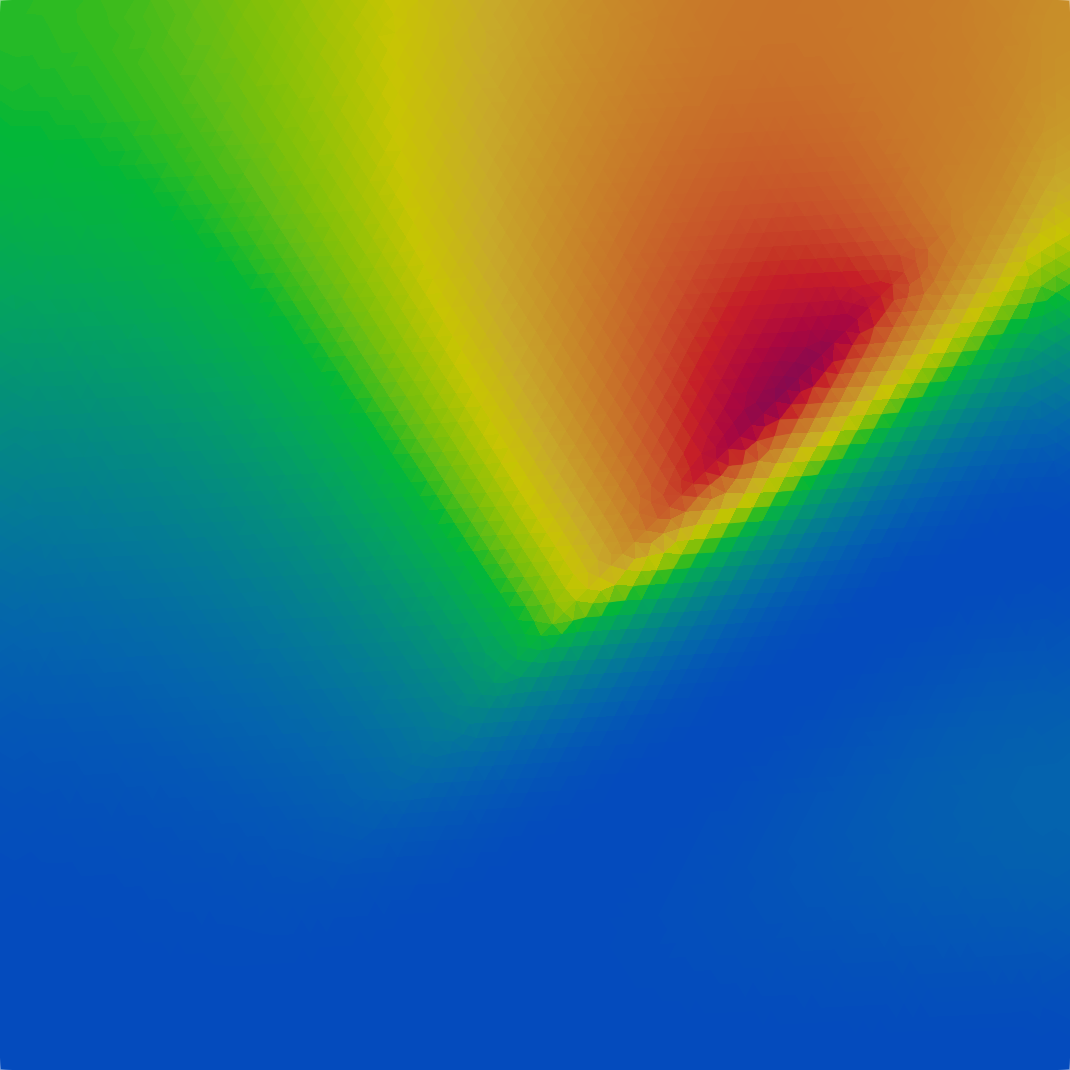}%
    \hspace*{0.02\textwidth}%
    \includegraphics[width=0.32\textwidth]{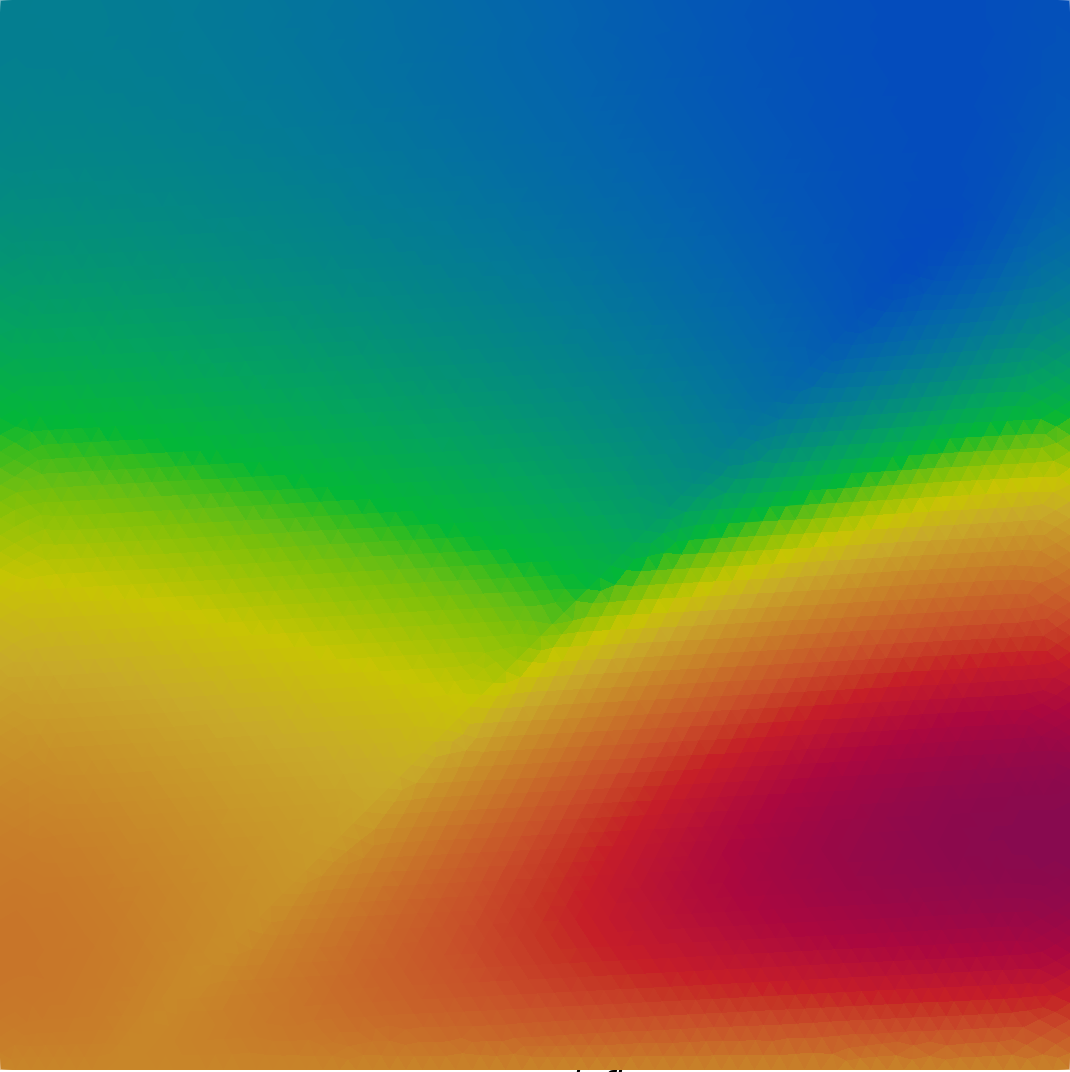}
    \includegraphics[width=0.32\textwidth]{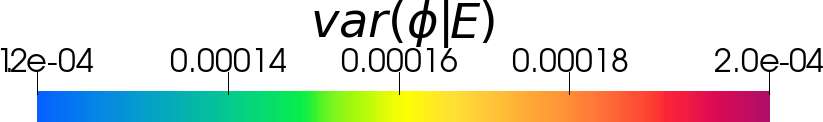}%
    \hspace*{0.02\textwidth}%
    \includegraphics[width=0.32\textwidth]{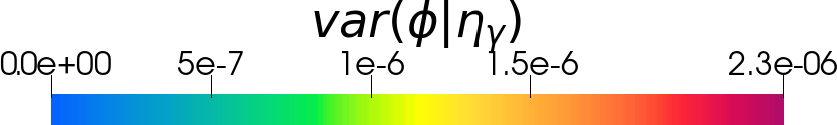}%
    \hspace*{0.02\textwidth}%
    \includegraphics[width=0.32\textwidth]{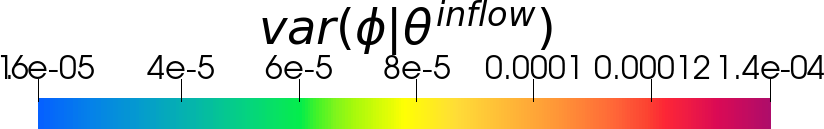}
    \caption{On the left the variance of porosity conditioned to the activation
    energy $E$, on the centre conditioned with $\eta_\gamma$ expansion, and on
    the right conditioned with the temperature inflow $\theta^{inflow}$. On the
    top at time $t=0.1T$ and on the bottom at final time $t=T$. Test case of
    Subsection \ref{subsec:case1_var}.}
    \label{fig:case1_var}
\end{figure}

In Figure \ref{fig:case1_var}, we compare the variances of the porosity in the
domain $\Omega$ induced by the uncertain data. Depending on the situation, it
might be more convenient to consider
the Sobol indices instead of the variances in particular when they span
different order of magnitudes. The importance of the activation energy over time
is rather interesting. In the beginning the fracture is highly permeable and
most of the flow is concentrated around the fracture. However, due to the deposition of new material the fracture becomes very low permeable; thus the water flow tends to avoid it and concentrates more in the right part of the domain, where
the solute is transported and becomes precipitate altering the porosity.
This effect, less evident, can be seen also in the spatial distribution of the Sobol index of $\eta_\gamma$. For the temperature,
we still notice the permeability change effect coupled with the inflow of higher
temperature that speeds up the precipitation process and, as a result, the porosity decay.

\begin{figure}[tbp]
    \centering
    \includegraphics[width=0.32\textwidth]{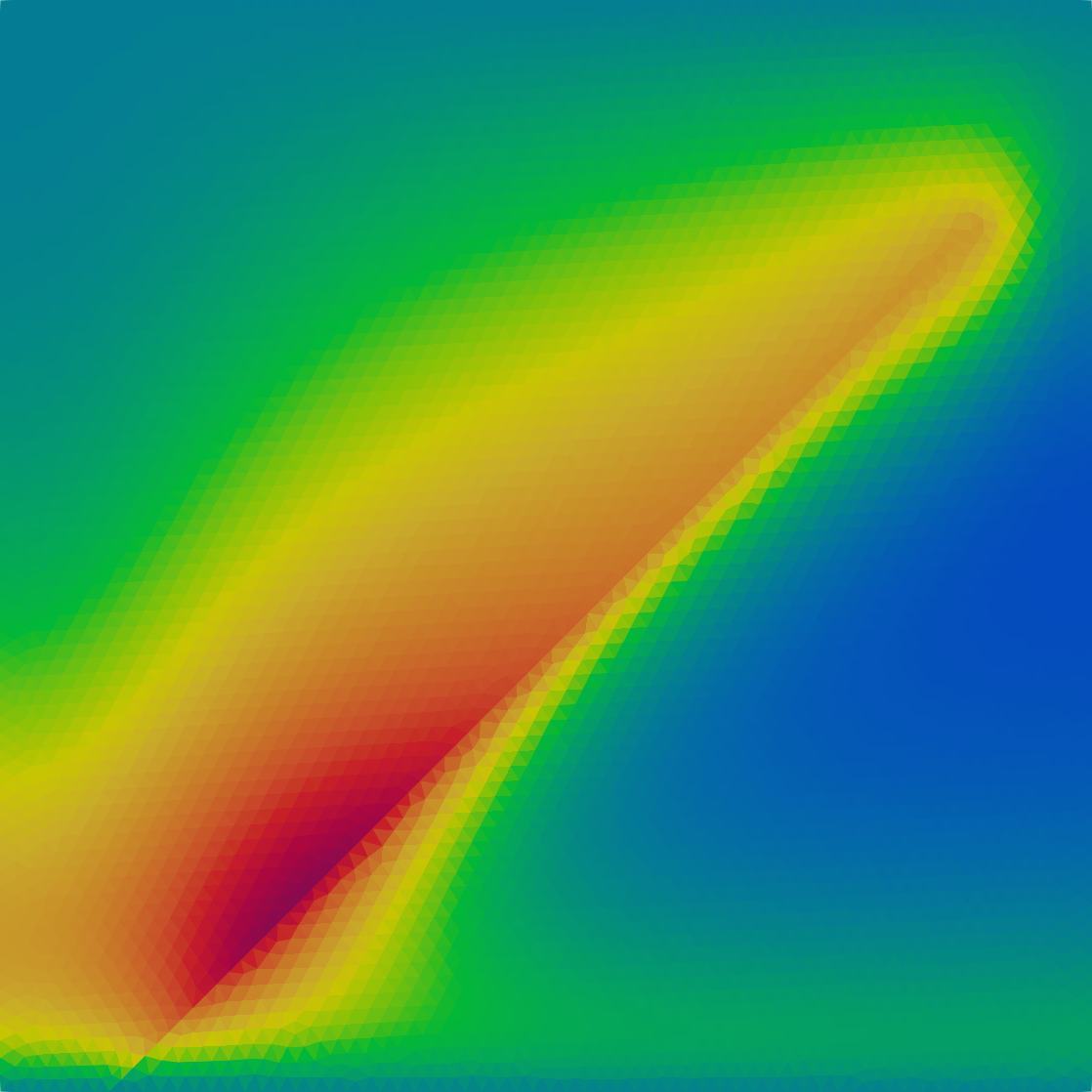}%
    \hspace*{0.02\textwidth}%
    \includegraphics[width=0.32\textwidth]{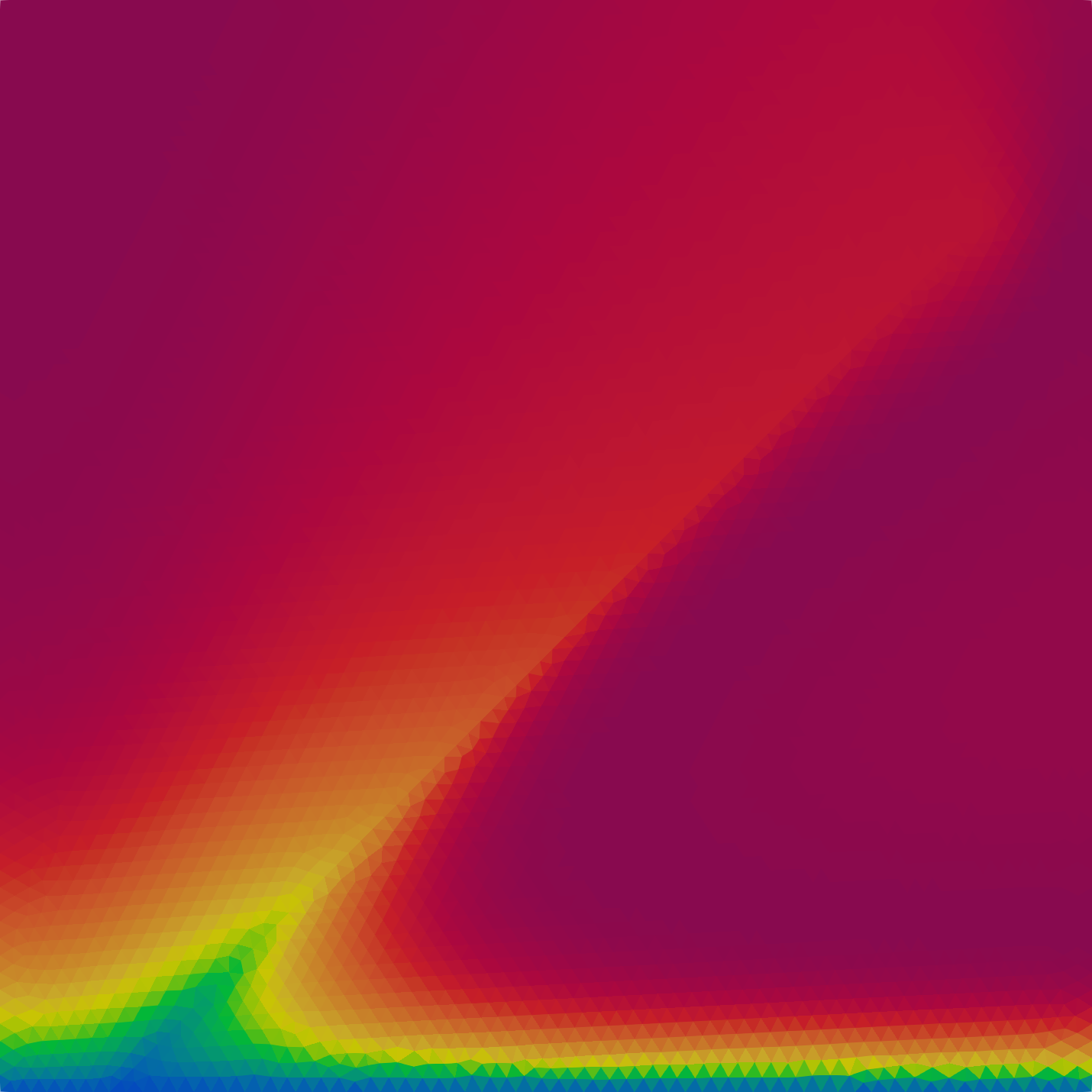}%
    \hspace*{0.02\textwidth}%
    \includegraphics[width=0.32\textwidth]{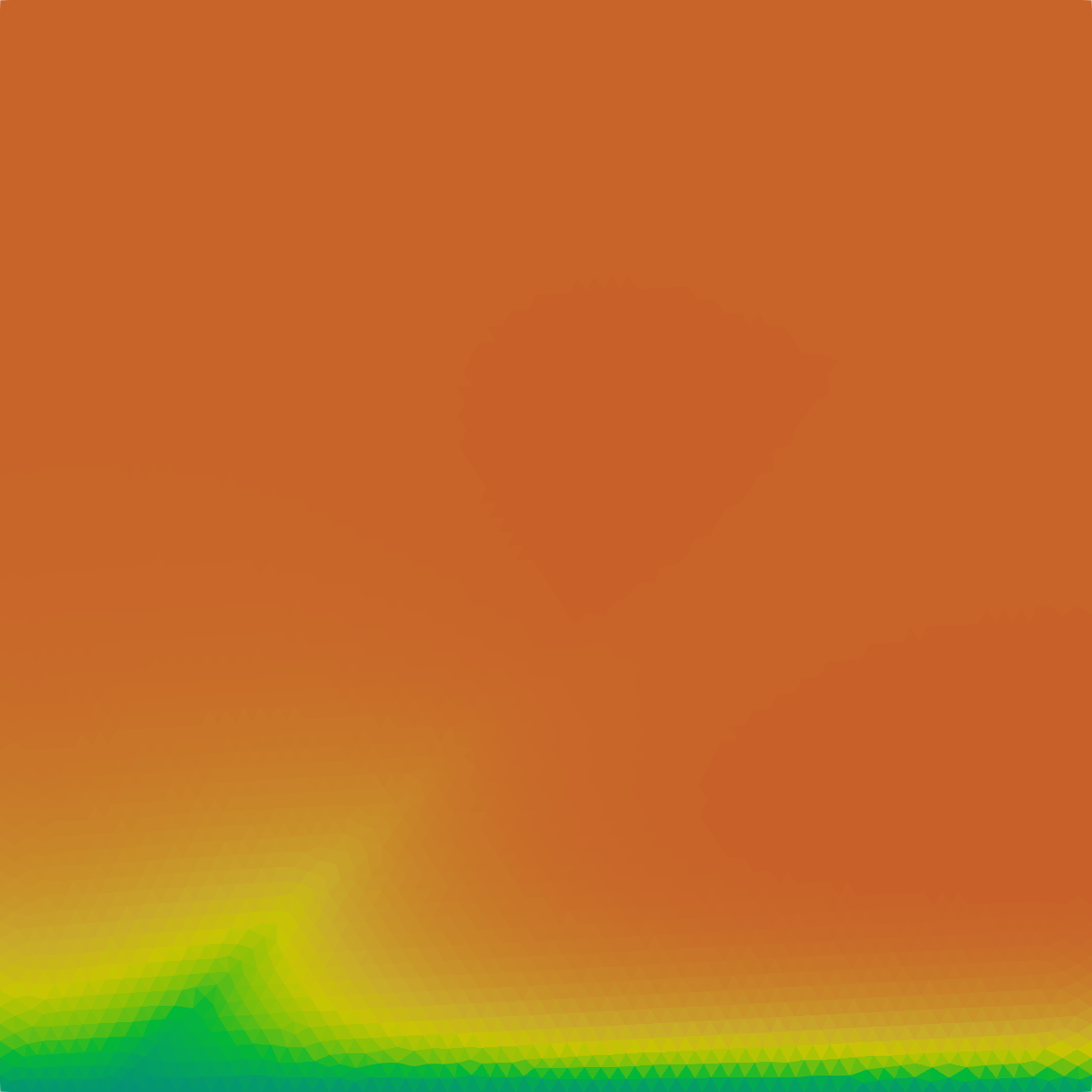}
    \includegraphics[width=0.32\textwidth]{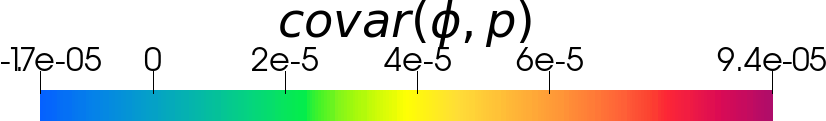}%
    \hspace*{0.02\textwidth}%
    \includegraphics[width=0.32\textwidth]{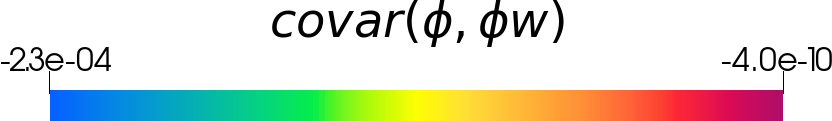}%
    \hspace*{0.02\textwidth}%
    \includegraphics[width=0.32\textwidth]{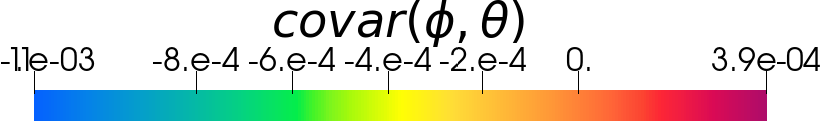}\\
    \vspace*{0.05\textwidth}
    \includegraphics[width=0.32\textwidth]{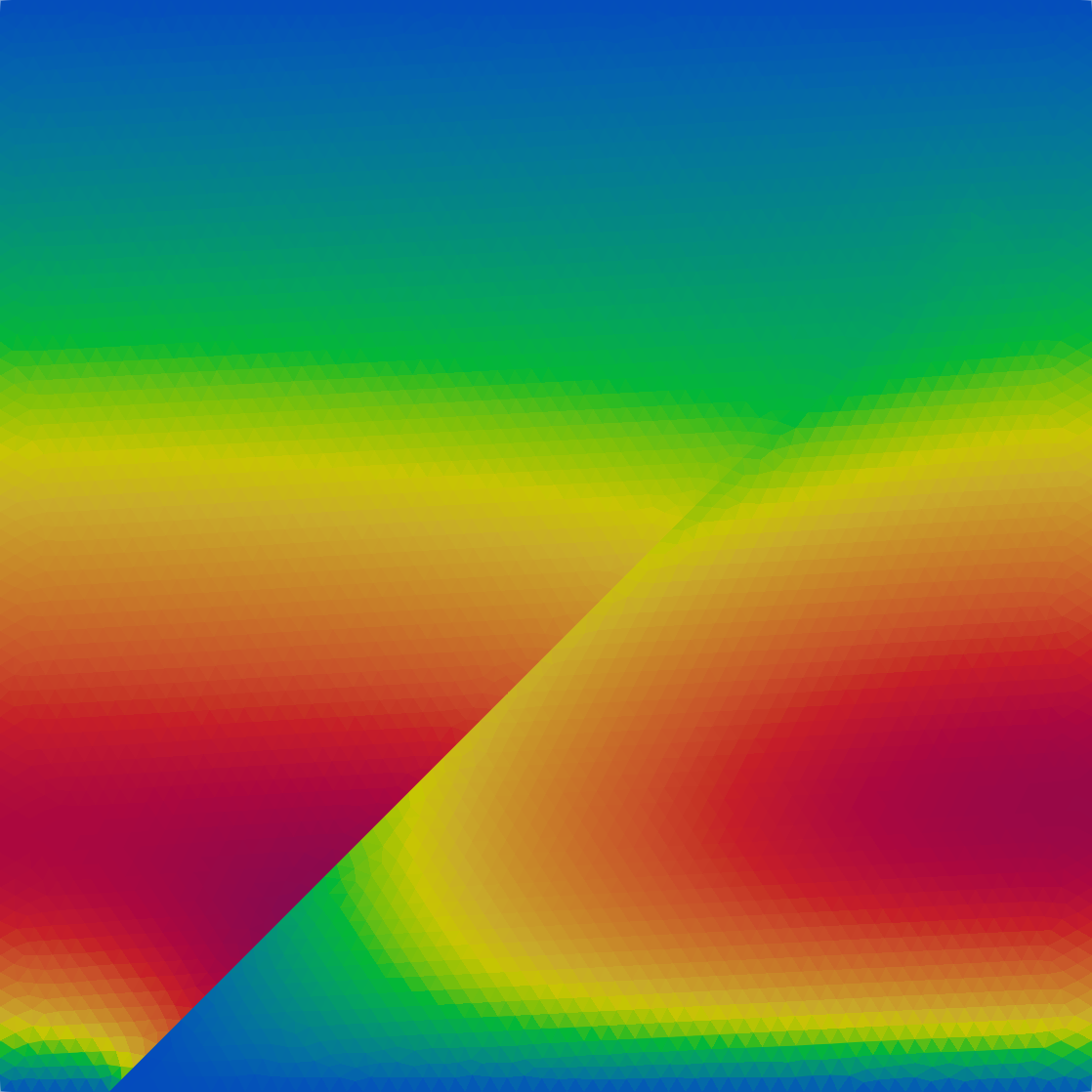}%
    \hspace*{0.02\textwidth}%
    \includegraphics[width=0.32\textwidth]{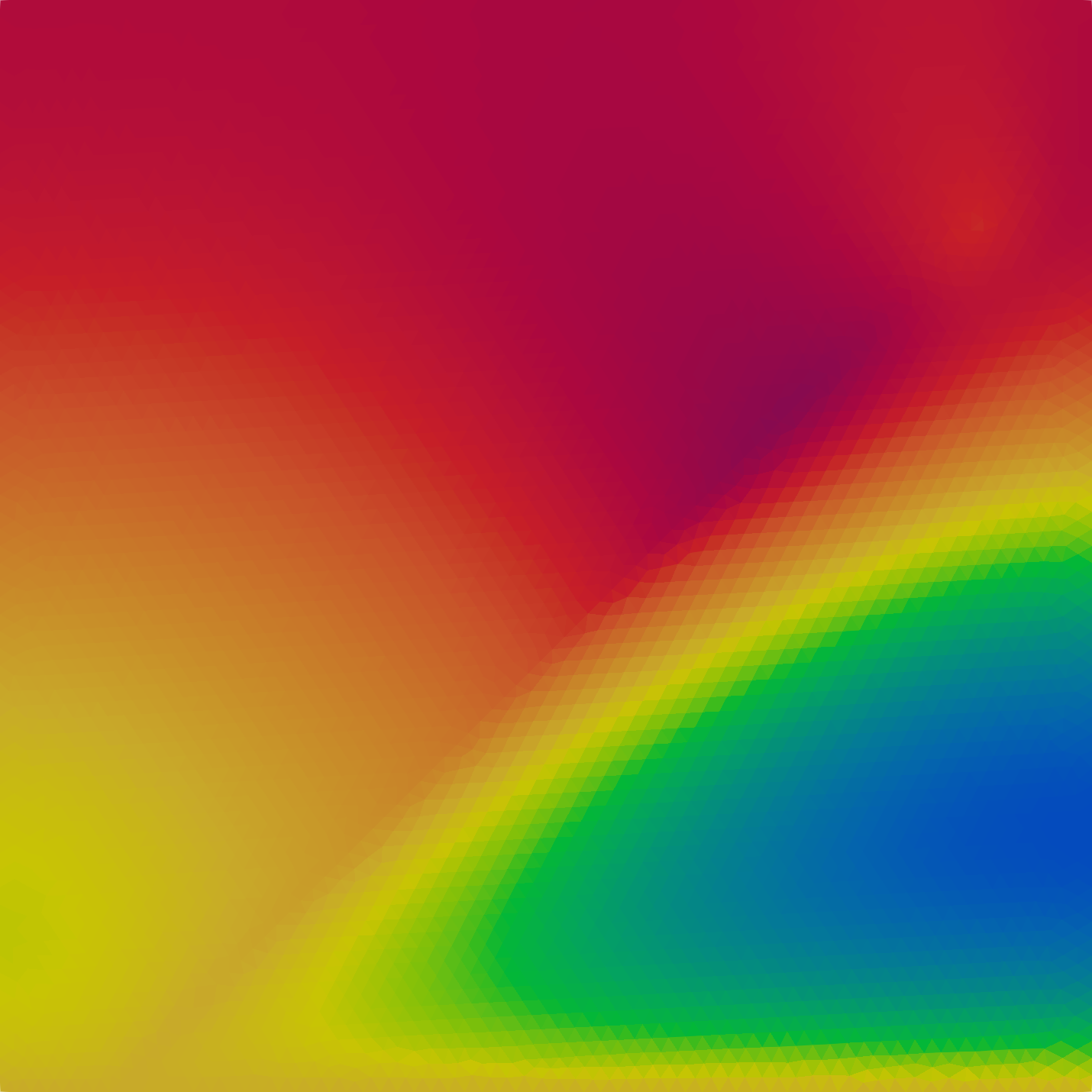}%
    \hspace*{0.02\textwidth}%
    \includegraphics[width=0.32\textwidth]{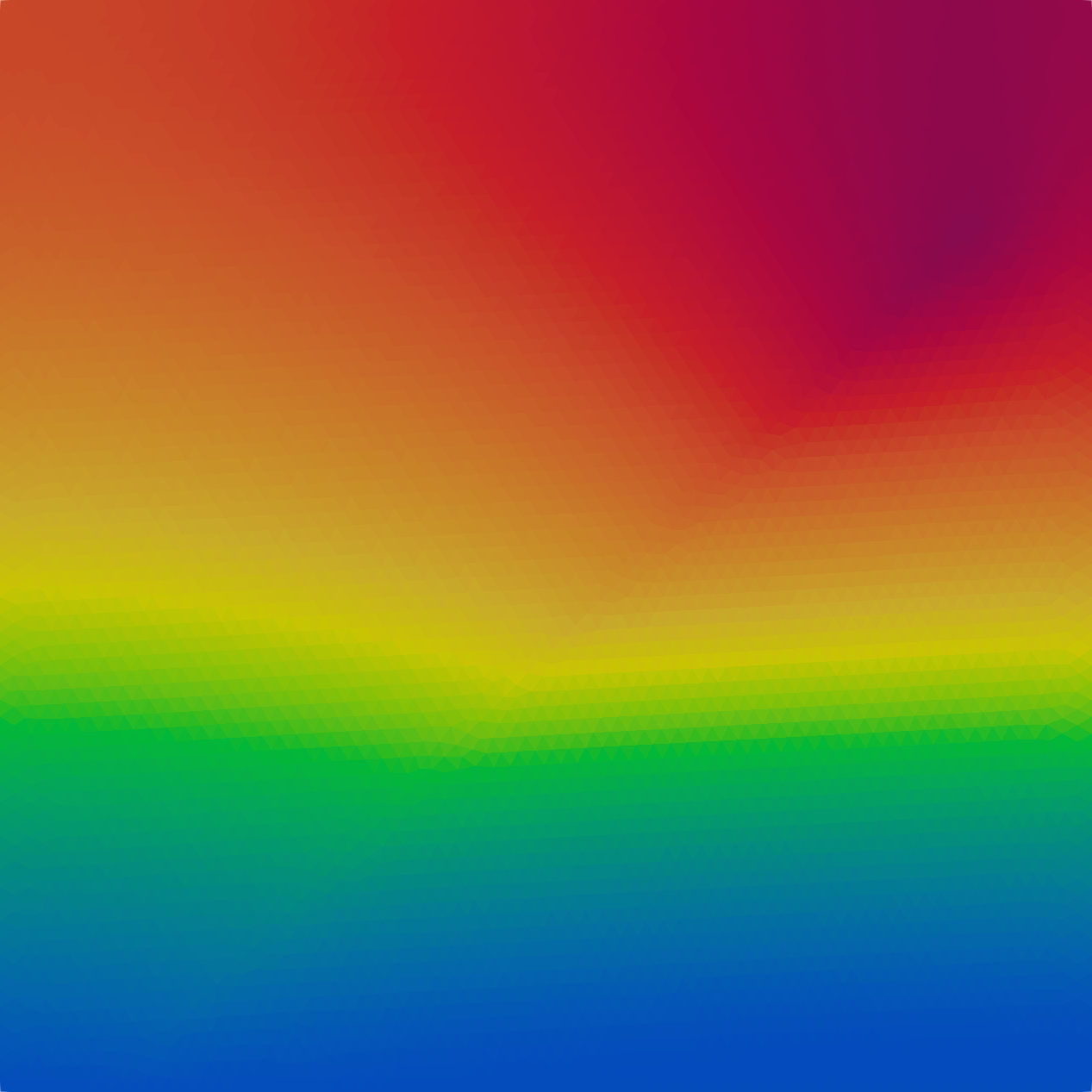}
    \includegraphics[width=0.32\textwidth]{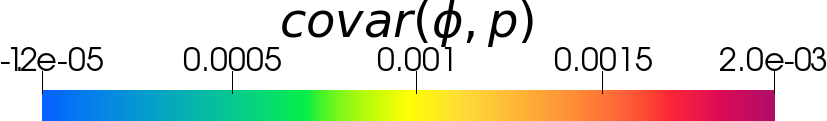}%
    \hspace*{0.02\textwidth}%
    \includegraphics[width=0.32\textwidth]{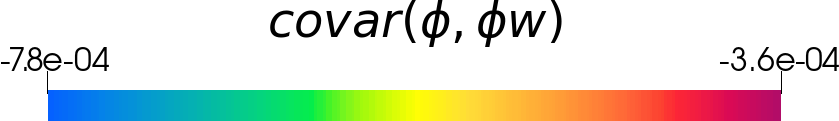}%
    \hspace*{0.02\textwidth}%
    \includegraphics[width=0.32\textwidth]{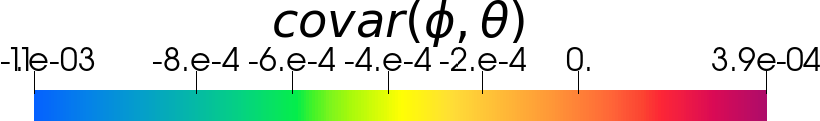}
    \caption{Covariances between different solutions in the porous media. On the
    top at time $t=0.1T$ and on the bottom at final time $t=T$. Test case of
    Subsection \ref{subsec:case1_var}.}
    \label{fig:case1_covar}
\end{figure}

Another important aspect is the interdependency of the output variables, which can be
expressed by their covariances. Figure \ref{fig:case1_covar} presents this
relation between the porosity in the media and several other variables at the
two considered times. In the pressure-porosity covariance we observe again the effect of the variation of the
fracture behavior in time, fromp permeable to impermeable. For the
porosity and $\phi_\Omega w_\Omega$ we observe a negative correlation, namely if more
precipitate is deposited in the porous media, less void space is left, and the
porosity diminishes. Finally, the temperature front can be seen in the plot of the 
covariance between $\phi_\Omega$ and $\theta_\Omega$, i.e. higher values of the
latter tend to facilitate the deposition of solute with concentration higher than the equilibrium. This increases the value of the precipitate and
consequently lowers the value of the porosity; conversely, in the top part of the domain the effect is the opposite due to the fact that thermal capacity and conductivity depend on the porosity in a complex way.

\subsubsection{Probability density functions}\label{subsec:case1_pdf}

\begin{figure}[tbp]
    \centering
    \includegraphics[width=1\textwidth]{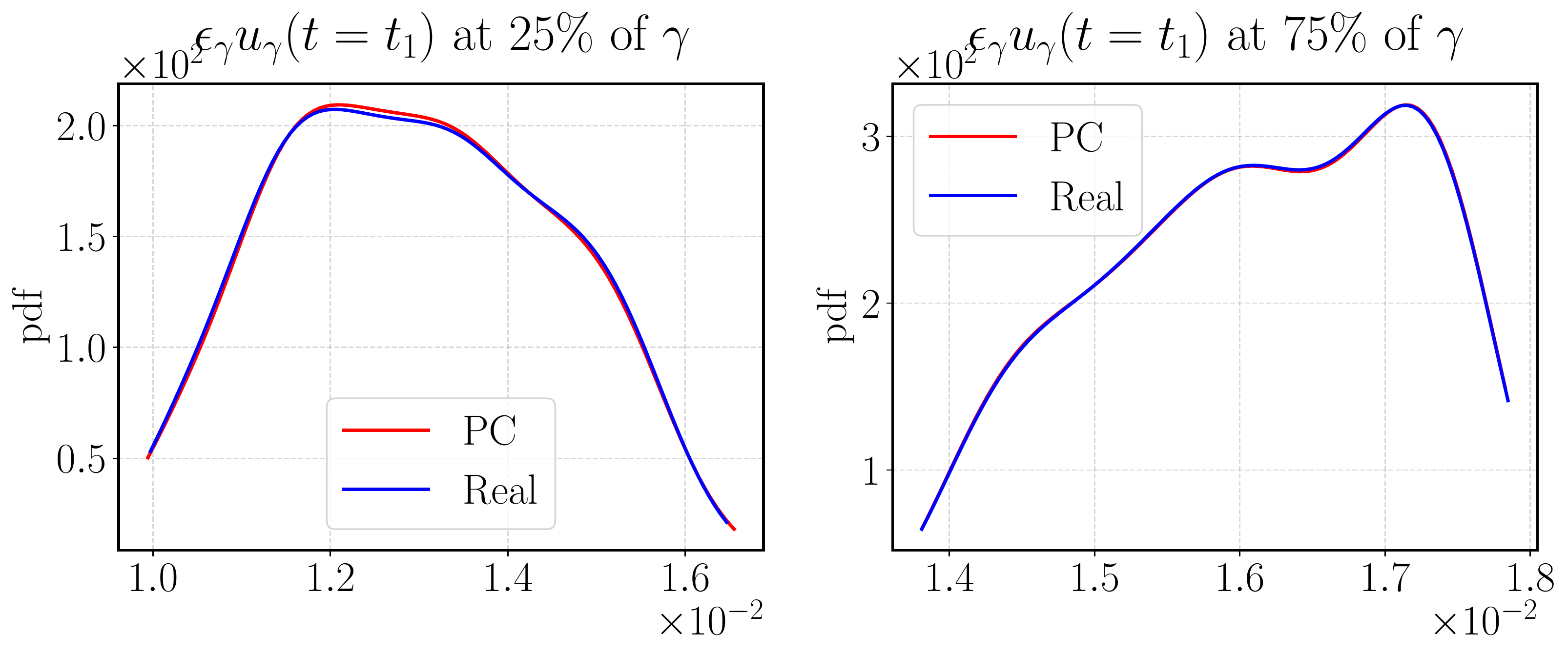}
    \includegraphics[width=1\textwidth]{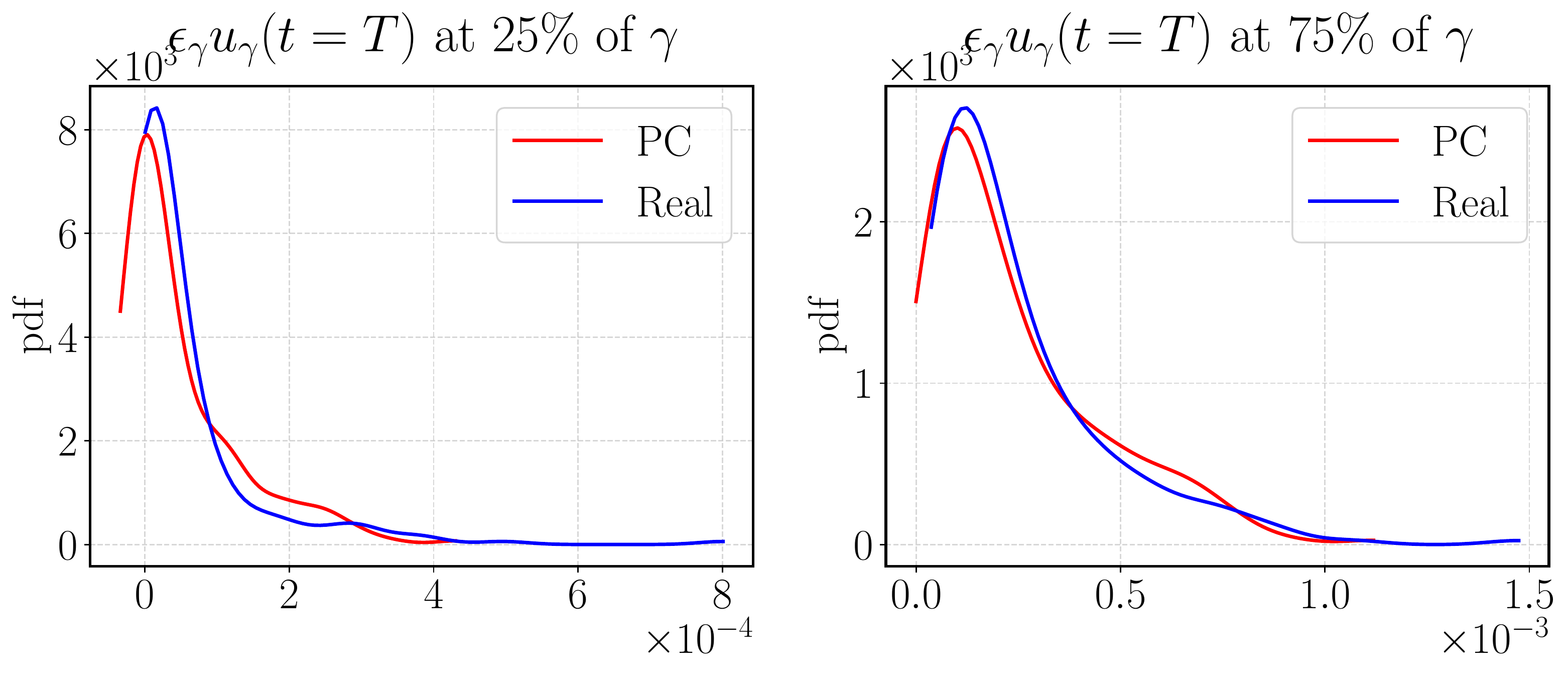}
    \caption{Probability distribution function of $\epsilon_\gamma u_\gamma$ for
    level 2 at two points in $\gamma$, on the
    top at time $t=0.1T$ and on the bottom at final time $t=T$. Test case of
    Subsection \ref{subsec:case1_pdf}.}
    \label{fig:case1_pdf_lev2}
\end{figure}

\begin{figure}[tbp]
    \centering
    \includegraphics[width=1\textwidth]{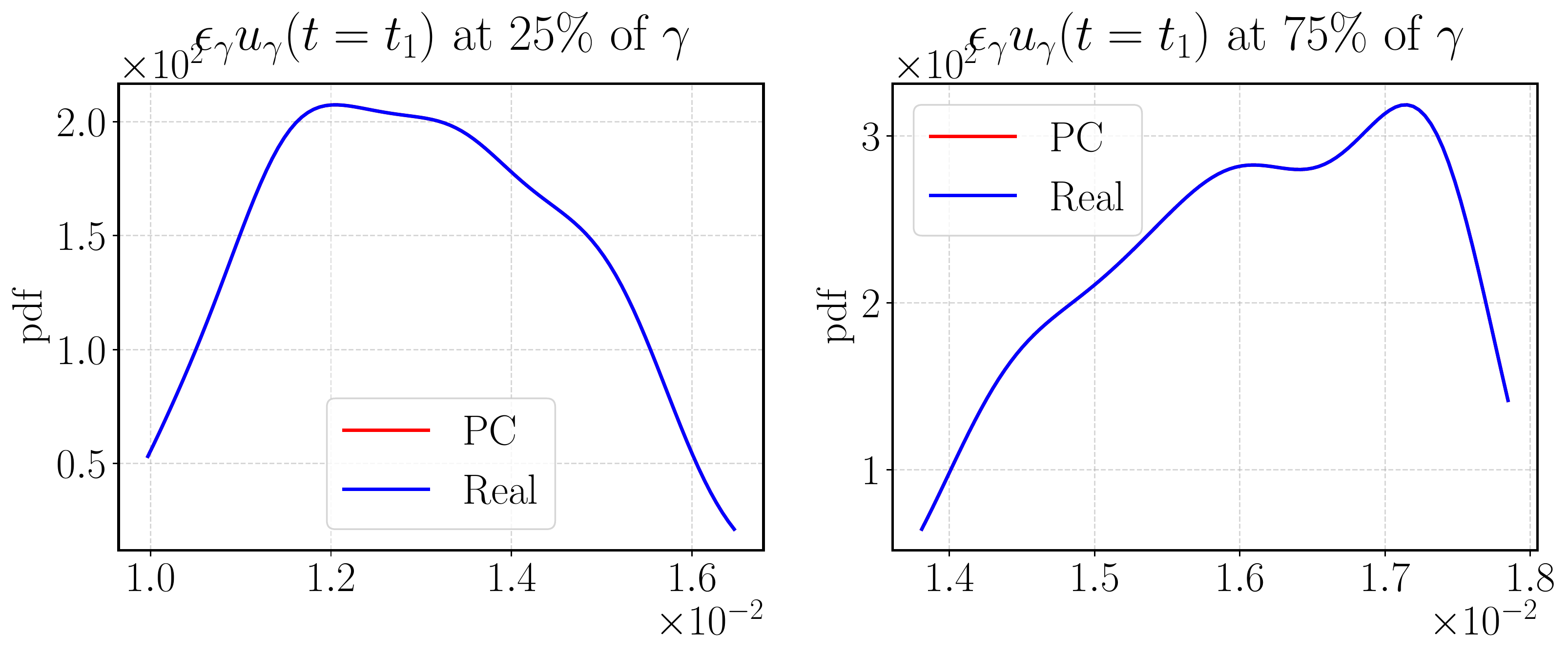}
    \includegraphics[width=1\textwidth]{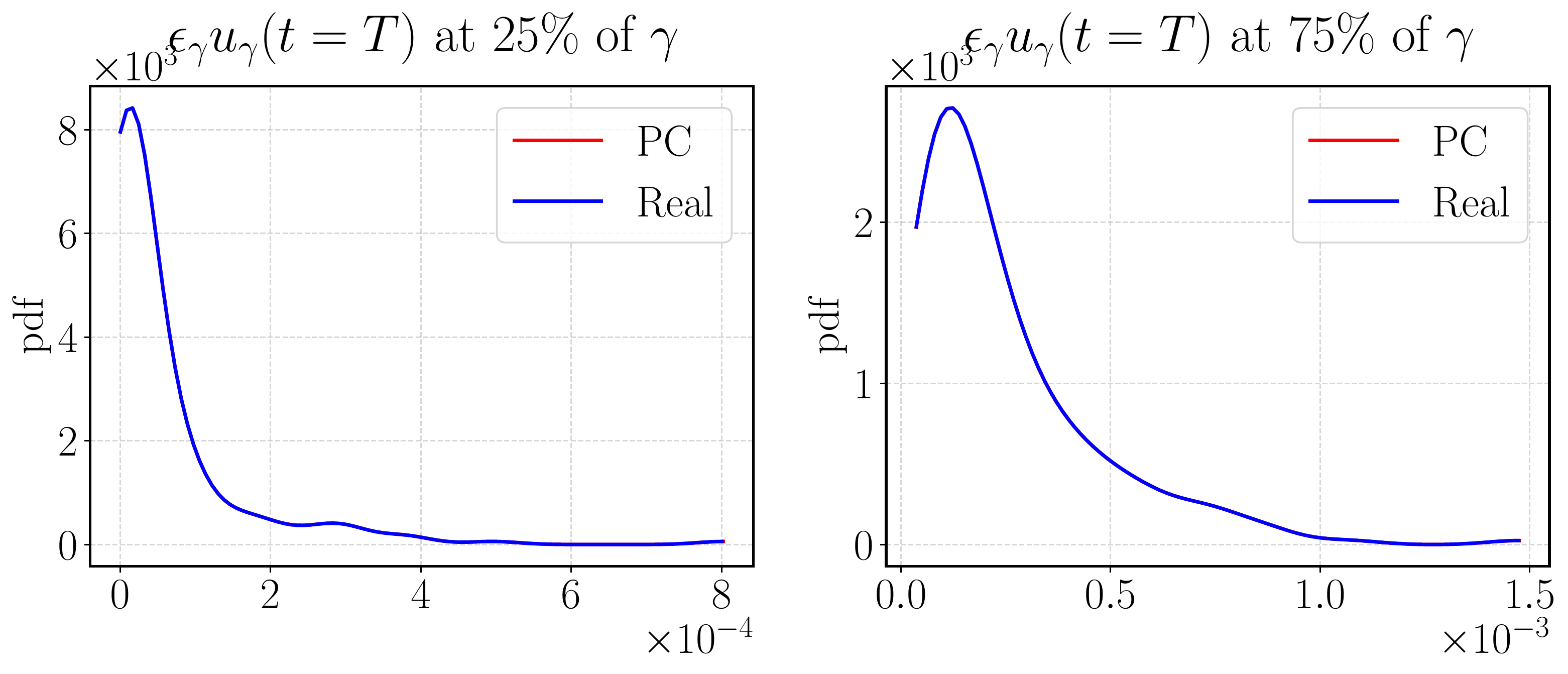}
    \caption{Probability distribution function of $\epsilon_\gamma u_\gamma$ for
    level 5 at two points in $\gamma$, on the
    top at time $t=0.1T$ and on the bottom at final time $t=T$. Test case of
    Subsection \ref{subsec:case1_pdf}.}
    \label{fig:case1_pdf_lev5}
\end{figure}

We consider now the probability density functions (PDFs) of some variables induced
by the uncertain data, which are uniformly distributed. Figure \ref{fig:case1_pdf_lev2} shows, for level 2, the distribution of
$\epsilon_\gamma u_\gamma$ at two points along the fracture, and for both times.
We notice that at the beginning of the simulation the PDFs are more spread
showing a high variability of the considered variable. However, at the end of
the simulation the uncertainty tends to become much smaller and the value is
more concentrated. Another important aspect is that the PC expansion outcomes
might not fulfill physical bounds, in this case we can get negative values of
$\epsilon_\gamma u_\gamma$ which are not correct. The situation improves by
considering a higher level of the sparse grid, indeed as represented in Figure
\ref{fig:case1_pdf_lev5} this phenomena is not present any more and the PDFs
constructed with the PC expansion are in good agreement with the one computed by
the full order model.

\subsection{Multiple fractures network}\label{subsec:example_multiple_fracture}

We consider now a test case with a network composed of multiple-fractures.
The geometry is given by the Benchmark 3 of \cite{Flemisch2016a}, where
fractures at $t=0$ are now considered all highly permeable with material
properties and problem data equal to the previous test case. A graphical representation of the
computational domain is given in Figure \ref{fig:case2_error} on the left.
Fractures $\gamma_3$ and $\gamma_8$ will be considered later for a specific
analysis.

\begin{figure}[tbp]
    \centering
    \includegraphics[width=0.32\textwidth]{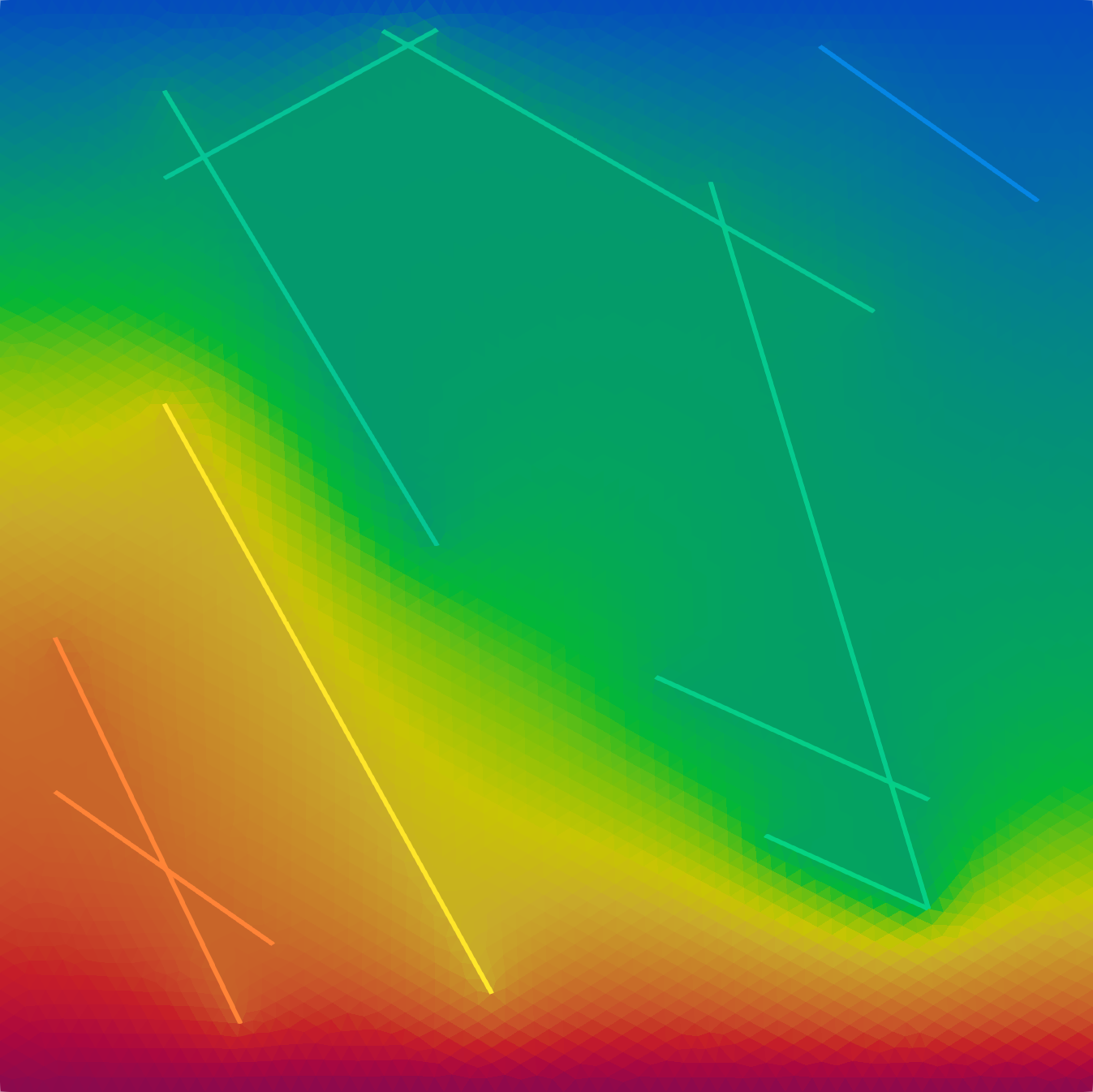}%
    \hspace*{0.02\textwidth}%
    \includegraphics[width=0.32\textwidth]{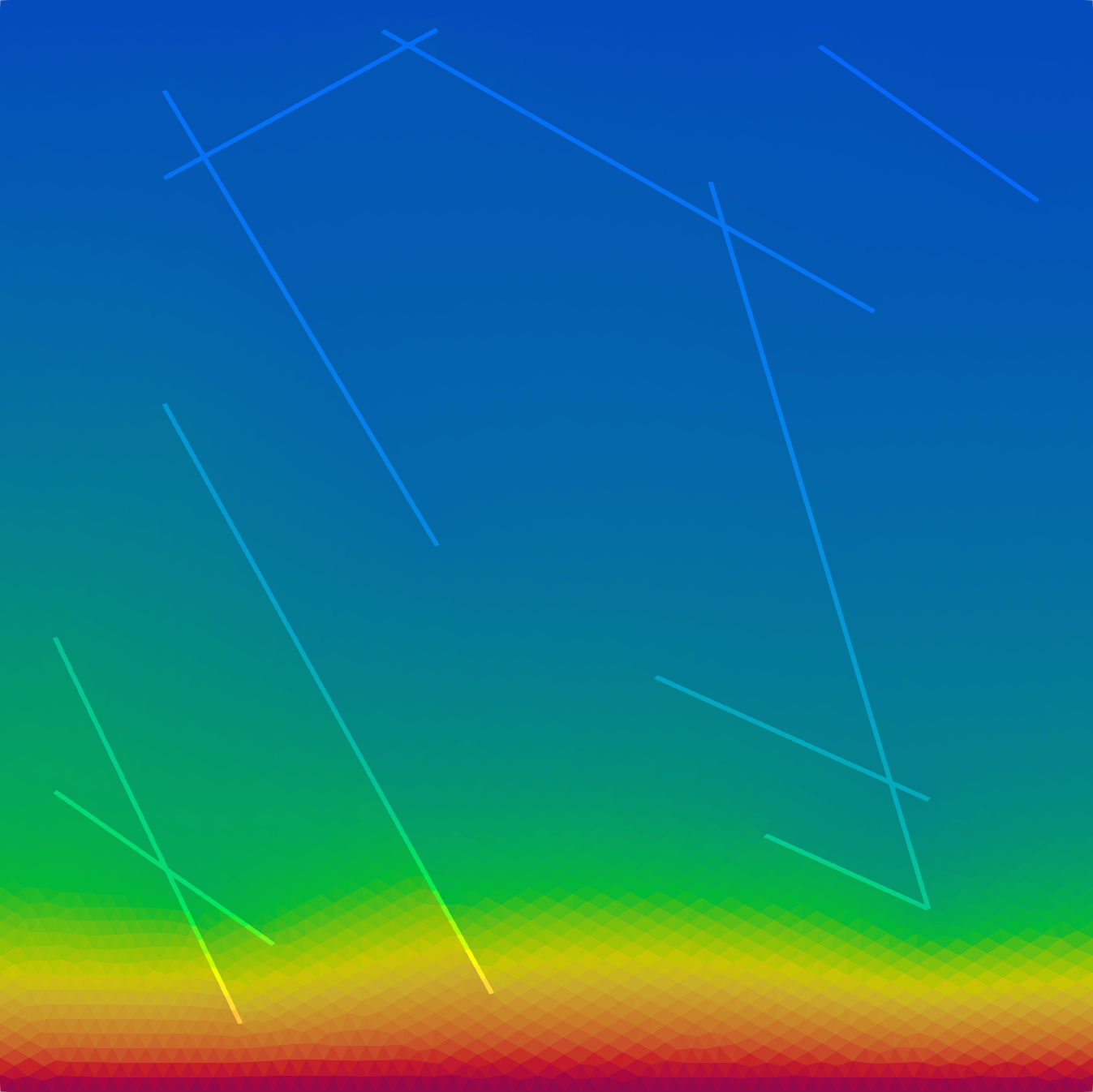}%
    \hspace*{0.02\textwidth}%
    \includegraphics[width=0.32\textwidth]{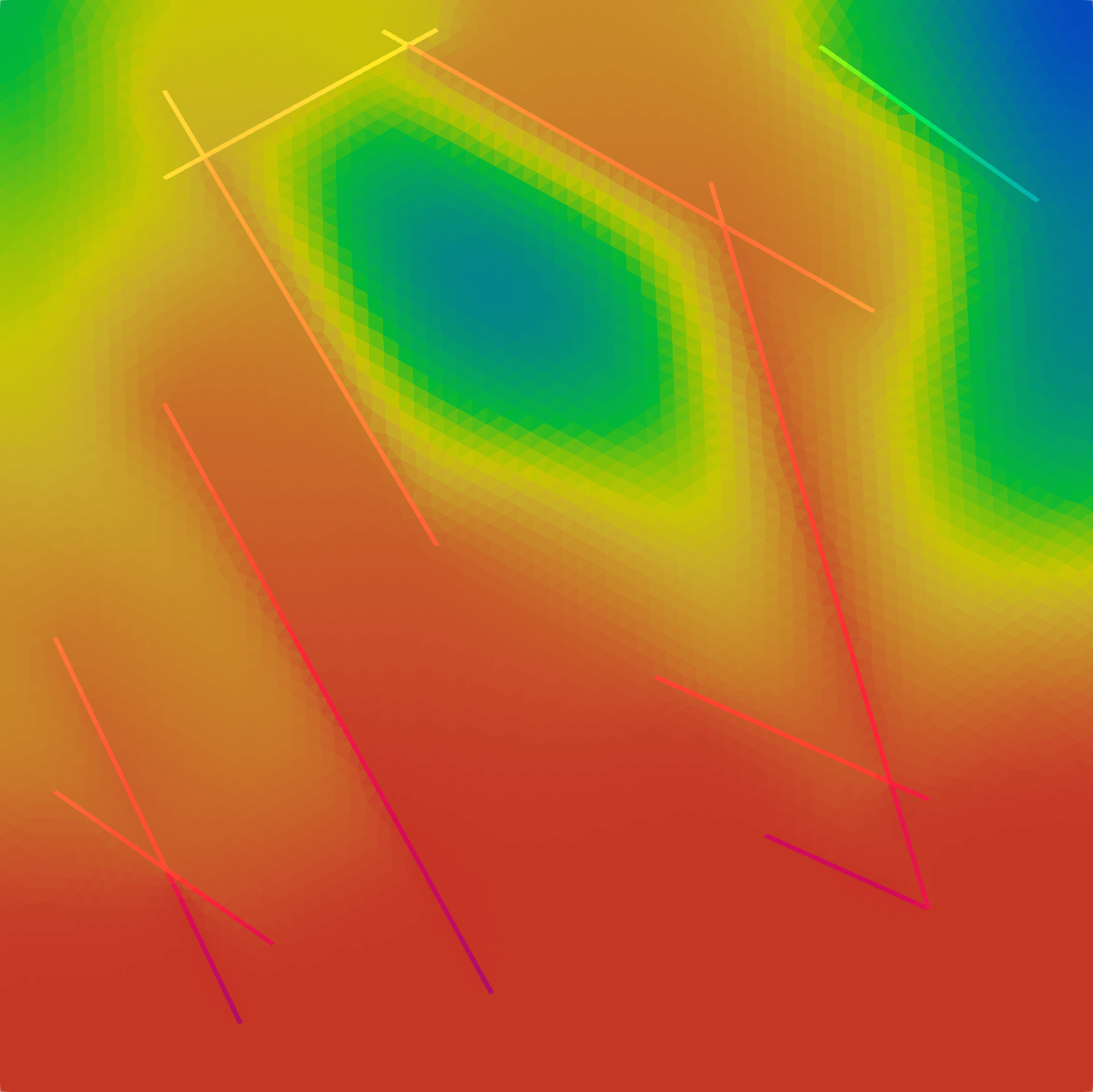}
    \hspace*{0.18\textwidth}%
    \includegraphics[width=0.32\textwidth]{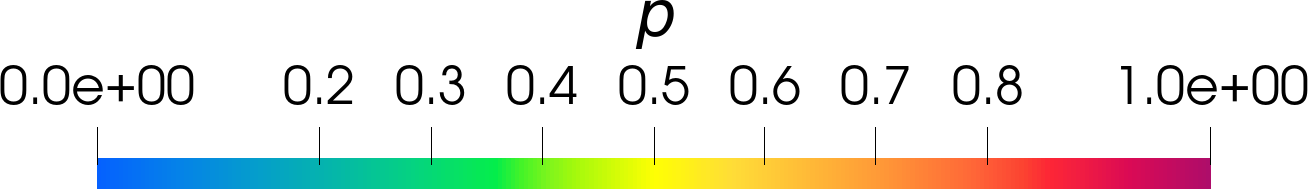}%
    \hspace*{0.18\textwidth}%
    \includegraphics[width=0.32\textwidth]{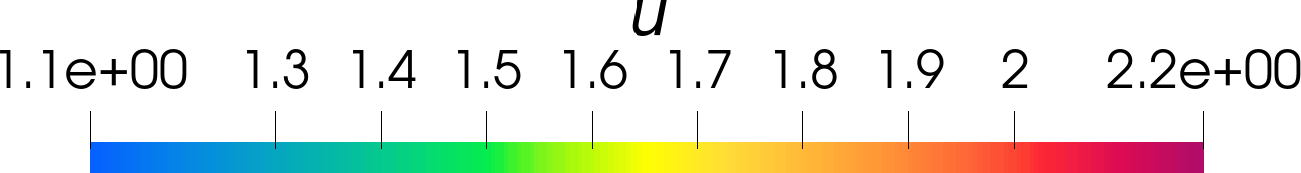}
    \caption{Reference solutions with mean value of the uncertain parameters.
    On the left pressure at time $t=0.1T$ and on the centre for $t=T$, on the
    right the solute for $t=0.1T$.  Test case of
    Subsection \ref{subsec:example_multiple_fracture}.}
    \label{fig:case2_solution}
\end{figure}

Some quantities from the reference numerical solution are reported in Figure
\ref{fig:case2_solution}, where it is possible to notice the variation of the pressure
distribution over time due to the sealing of the fractures and the transport of
the solute when the fractures are still highly permeable.

\subsubsection{Convergence}\label{subsec:case2_conv}

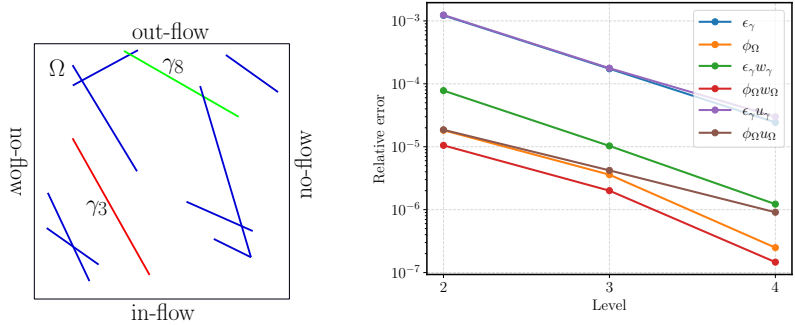
\begin{figure}[tbp]
    \centering
    \resizebox{0.33\textwidth}{!}{\fontsize{1cm}{2cm}\selectfont
\begingroup%
  \makeatletter%
  \providecommand\color[2][]{%
    \errmessage{(Inkscape) Color is used for the text in Inkscape, but the package 'color.sty' is not loaded}%
    \renewcommand\color[2][]{}%
  }%
  \providecommand\transparent[1]{%
    \errmessage{(Inkscape) Transparency is used (non-zero) for the text in Inkscape, but the package 'transparent.sty' is not loaded}%
    \renewcommand\transparent[1]{}%
  }%
  \providecommand\rotatebox[2]{#2}%
  \newcommand*\fsize{\dimexpr\f@size pt\relax}%
  \newcommand*\lineheight[1]{\fontsize{\fsize}{#1\fsize}\selectfont}%
  \ifx\svgwidth\undefined%
    \setlength{\unitlength}{311.96932983bp}%
    \ifx\svgscale\undefined%
      \relax%
    \else%
      \setlength{\unitlength}{\unitlength * \real{\svgscale}}%
    \fi%
  \else%
    \setlength{\unitlength}{\svgwidth}%
  \fi%
  \global\let\svgwidth\undefined%
  \global\let\svgscale\undefined%
  \makeatother%
  \begin{picture}(1,0.98417801)%
    \lineheight{1}%
    \setlength\tabcolsep{0pt}%
    \put(0,0){\includegraphics[width=\unitlength,page=1]{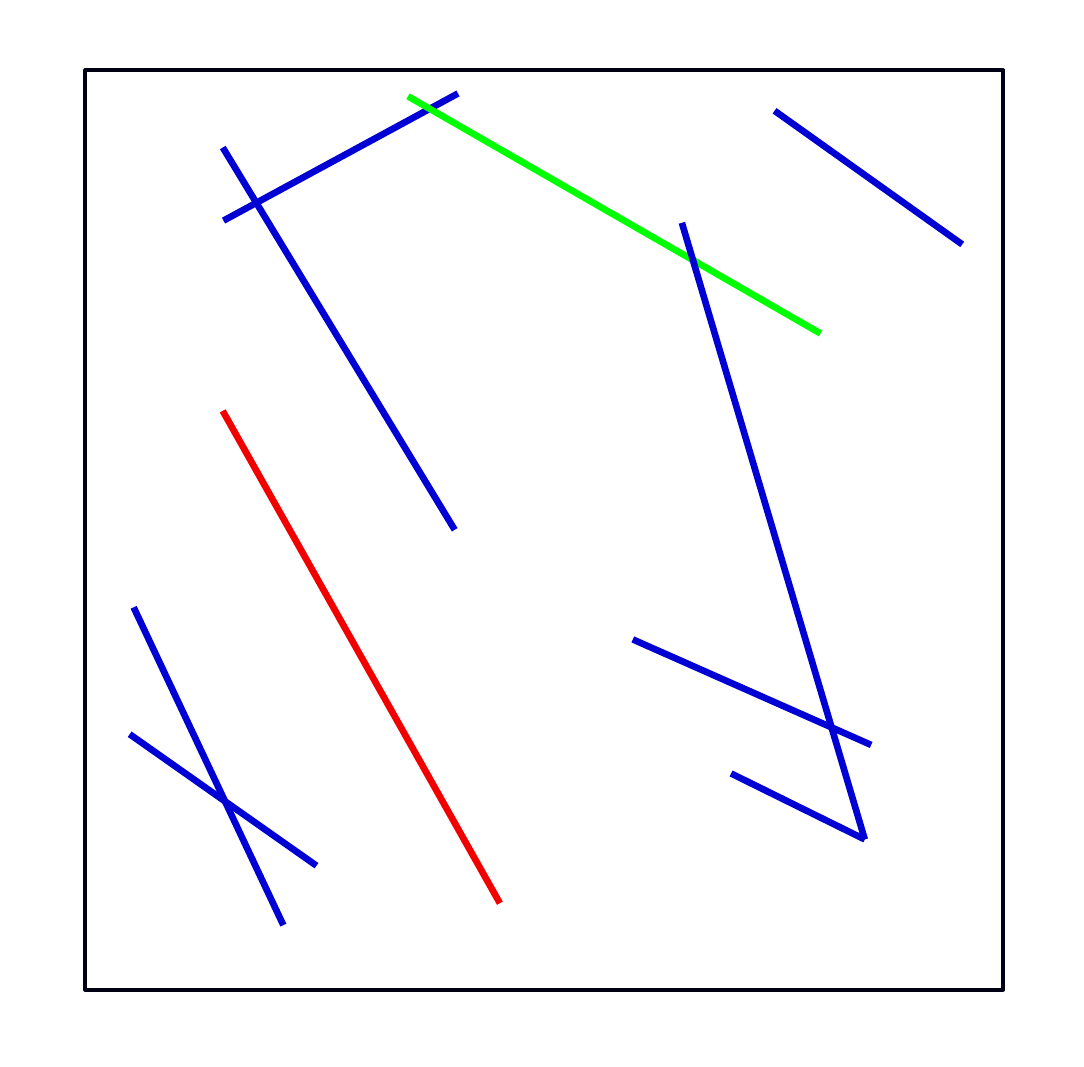}}%
    \put(0.40318924,0.93547031){\color[rgb]{0,0,0}\makebox(0,0)[lt]{\lineheight{1.25}\smash{\begin{tabular}[t]{l}out-flow\end{tabular}}}}%
    \put(0.40068498,0.0009078){\color[rgb]{0,0,0}\makebox(0,0)[lt]{\lineheight{1.25}\smash{\begin{tabular}[t]{l}in-flow\end{tabular}}}}%
    \put(0.9990922,0.40283447){\color[rgb]{0,0,0}\rotatebox{90}{\makebox(0,0)[lt]{\lineheight{1.25}\smash{\begin{tabular}[t]{l}no-flow\end{tabular}}}}}%
    \put(0.00090779,0.64233492){\color[rgb]{0,0,0}\rotatebox{-90}{\makebox(0,0)[lt]{\lineheight{1.25}\smash{\begin{tabular}[t]{l}no-flow\end{tabular}}}}}%
    \put(0.12781315,0.81059003){\color[rgb]{0,0,0}\makebox(0,0)[lt]{\lineheight{1.25}\smash{\begin{tabular}[t]{l}$\Omega$\end{tabular}}}}%
    \put(0.250403155,0.36449245){\color[rgb]{0,0,0}\makebox(0,0)[lt]{\lineheight{1.25}\smash{\begin{tabular}[t]{l}$\gamma_3$\end{tabular}}}}%
    \put(0.50403155,0.836449245){\color[rgb]{0,0,0}\makebox(0,0)[lt]{\lineheight{1.25}\smash{\begin{tabular}[t]{l}$\gamma_8$\end{tabular}}}}%
  \end{picture}%
\endgroup%
}
    \hspace*{0.05\textwidth}%
    \includegraphics[width=0.475\textwidth]{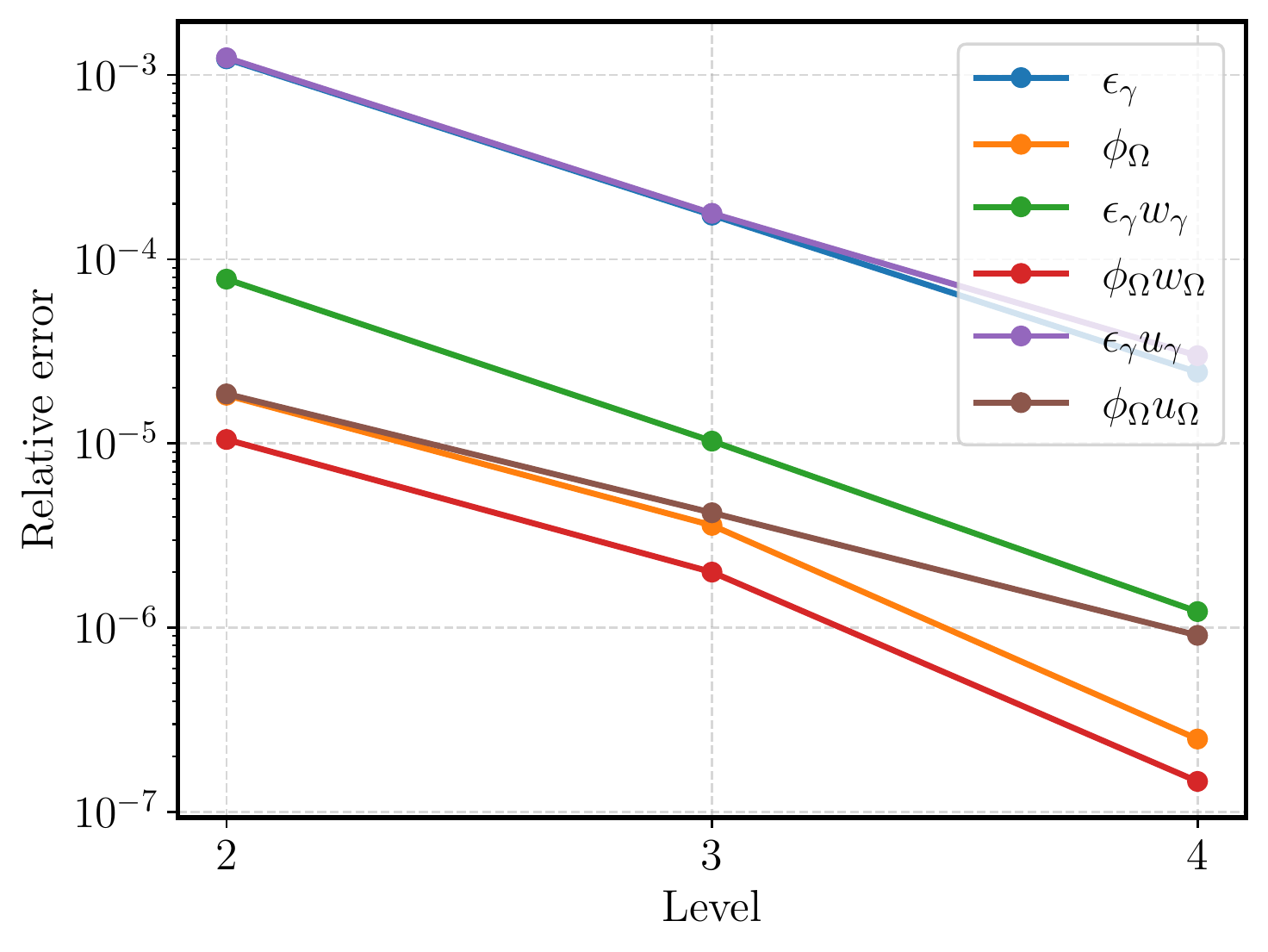}%
    \caption{On the right the computational domain with boundary conditions and
    on the left the error convergence for increasing level of the sparse grids. Test case of
    Subsection \ref{subsec:case2_conv}.}
    \label{fig:case2_error}
\end{figure}

We discuss now the convergence properties of the PC expansion. 
On the right of Figure \ref{fig:case2_error}, we plot the error decay for increasing sparse grid level for multiple variables. Also in this case, the exponential
decay expected is confirmed for all the variables.

\begin{figure}[tbp]
    \centering
    \includegraphics[width=0.32\textwidth]{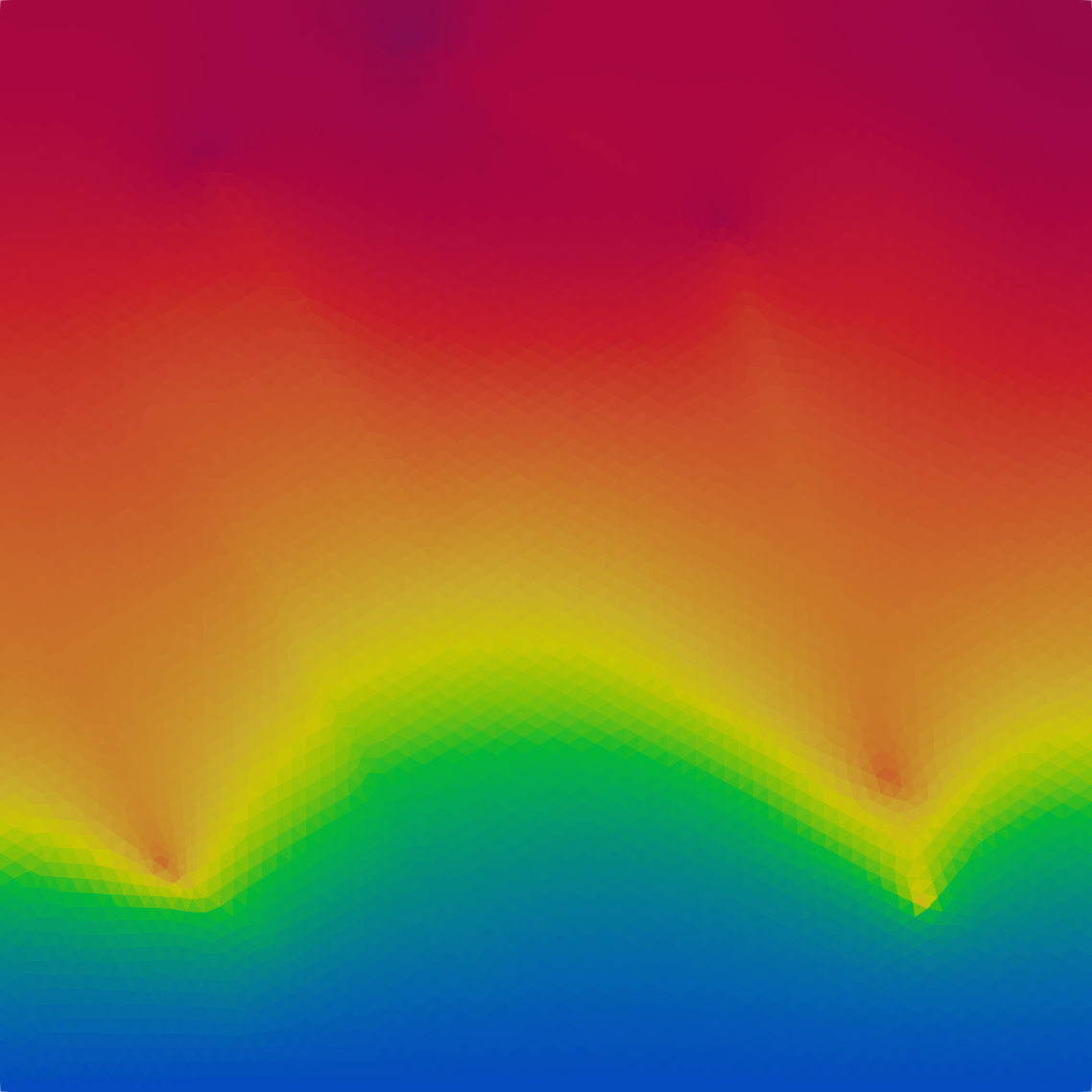}%
    \hspace*{0.02\textwidth}%
    \includegraphics[width=0.32\textwidth]{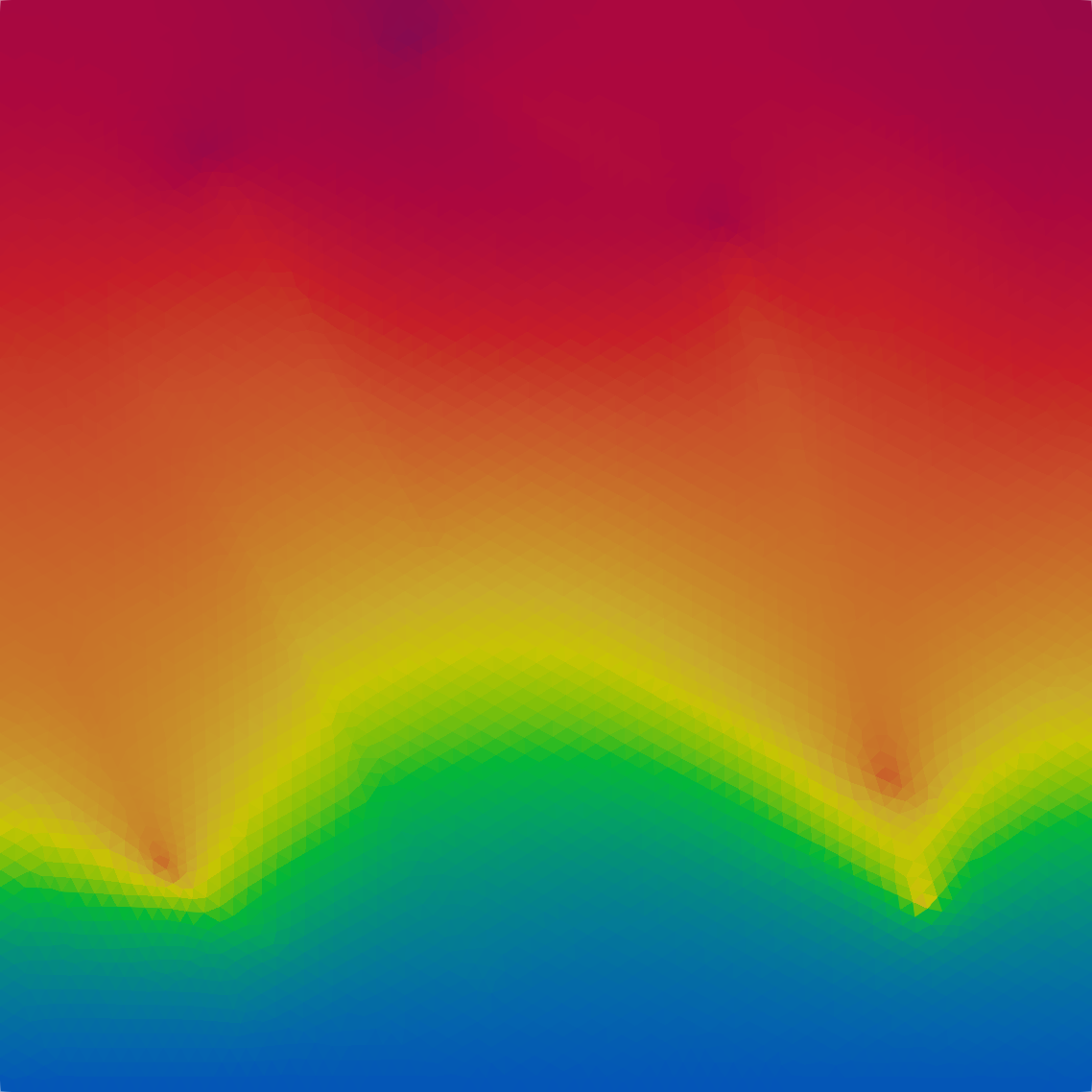}%
    \hspace*{0.02\textwidth}%
    \includegraphics[width=0.32\textwidth]{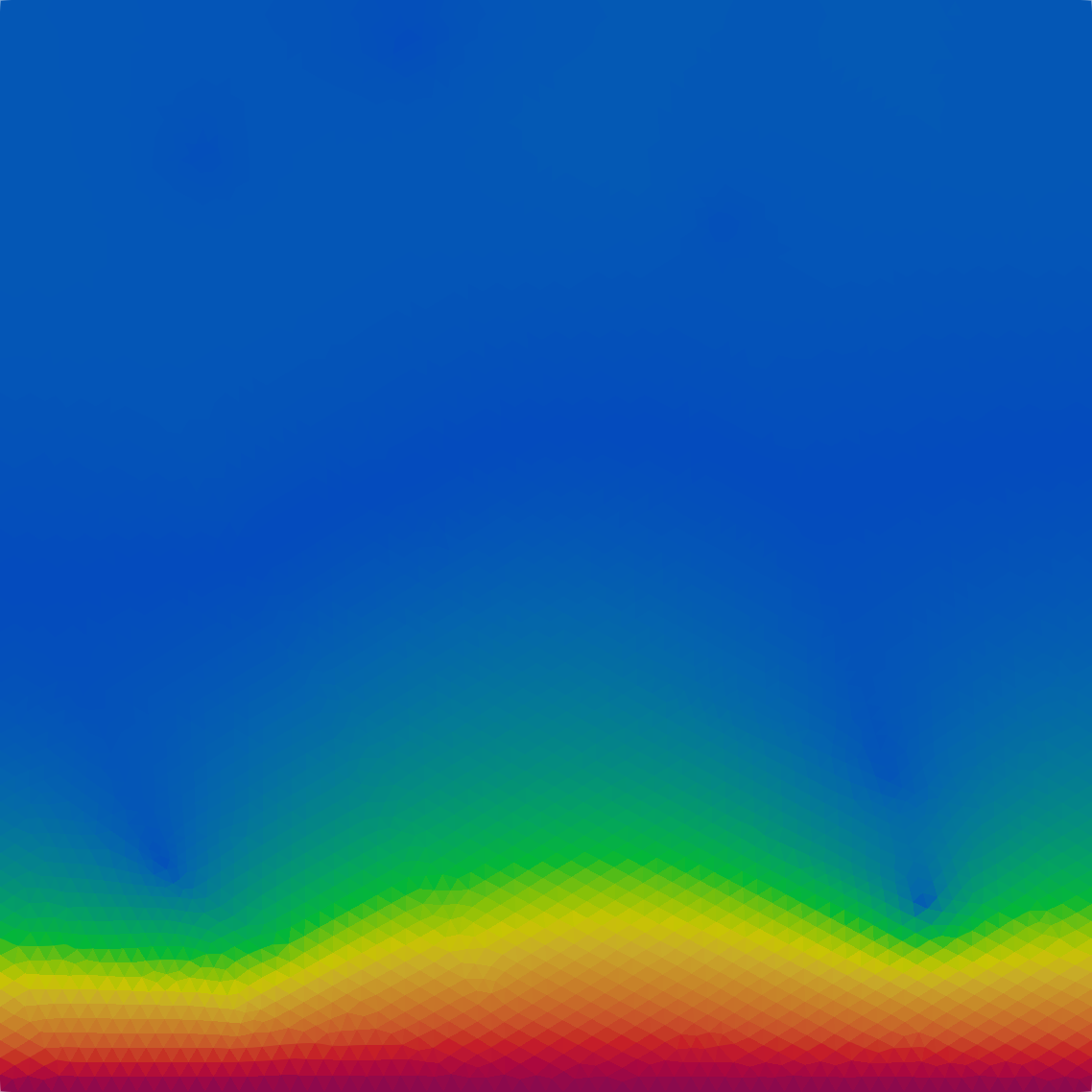}
    \includegraphics[width=0.32\textwidth]{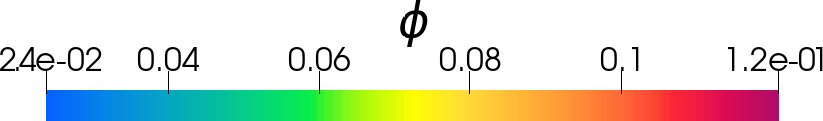}%
    \hspace*{0.02\textwidth}%
    \includegraphics[width=0.32\textwidth]{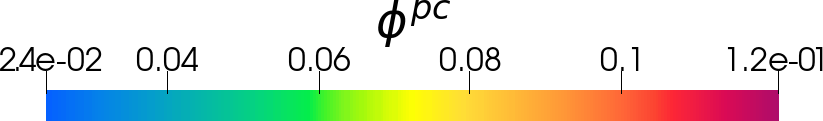}%
    \hspace*{0.02\textwidth}%
    \includegraphics[width=0.32\textwidth]{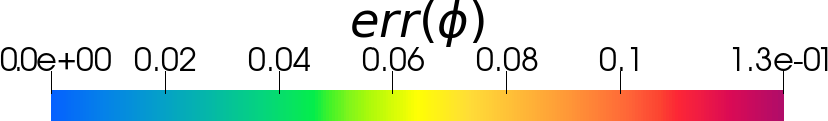}
    \caption{On the left porosity in the media computed with the original model
    and on the centre with the polynomial chaos expansion, on the right the
    error between them.  Test case of
    Subsection \ref{subsec:case2_conv}.}
    \label{fig:case2_porosity}
\end{figure}

In Figure \ref{fig:case2_porosity} we compare the porosity computed by the
differential model with the one computed by the PC expansion. On the right we also
represent the relative error. The two solutions are in good agreement with a maximum error of $13\%$ confined at the
bottom of the domain, the error is much smaller in the other parts being of the
order of $5\%$ or less. The surrogate model provided by the PC expansion gives a
satisfactory result at (almost) no additional computational cost.

\begin{figure}[tbp]
    \centering
    \includegraphics[width=1\textwidth]{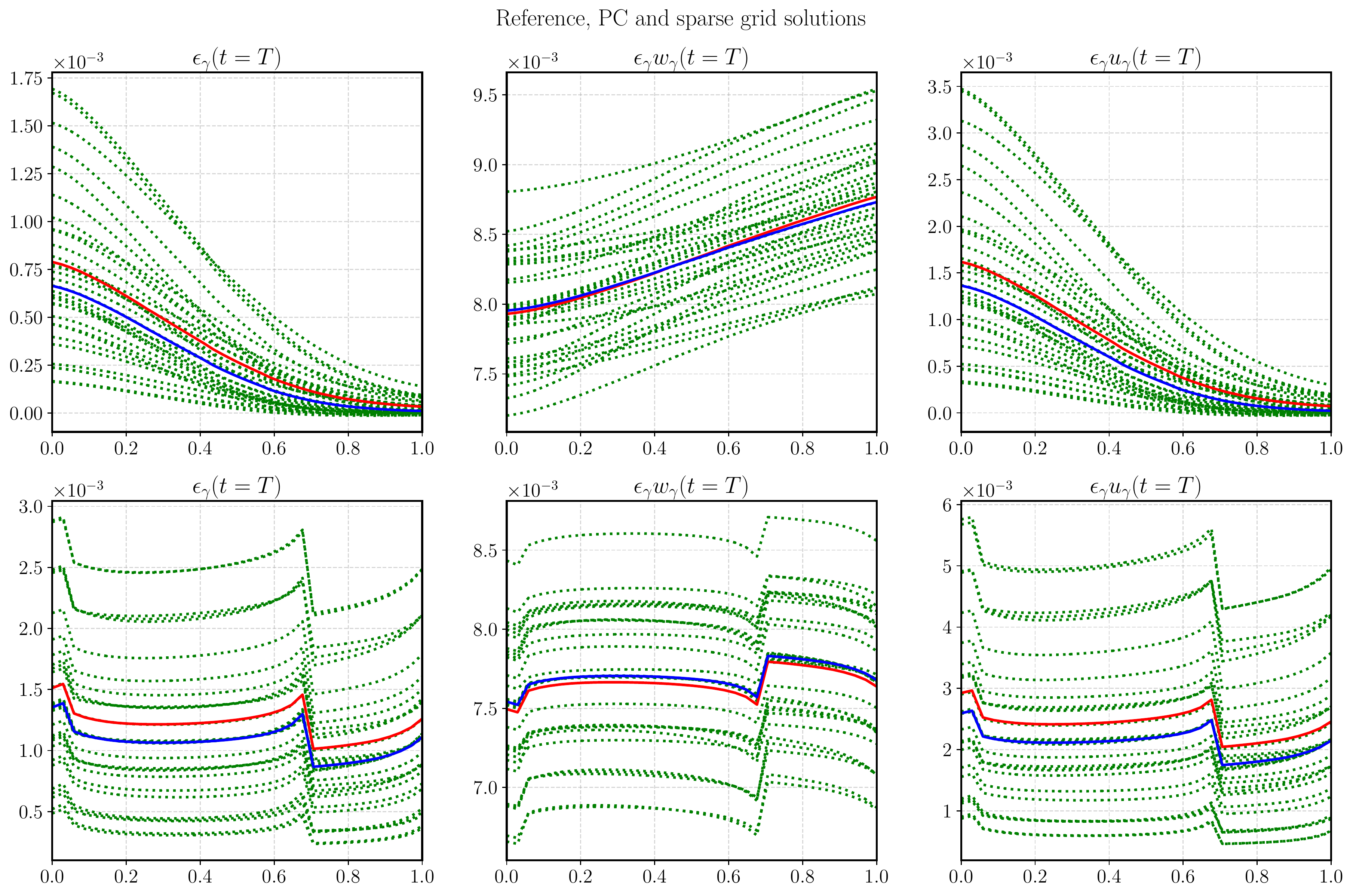}
    \caption{Solutions along two $\gamma_i$ and in blue computed with the original model,
    in red with the polynomial chaos expansion and in green computed by the original model
    on each sparse grid node. On the top $\gamma_3$ and on the bottom for
    $\gamma_8$. Level considered 2. Test case of
    Subsection \ref{subsec:case2_conv}.}
    \label{fig:case2_1d}
\end{figure}

Figure \ref{fig:case2_1d} compares some of the variables along the two fractures
$\gamma_3$ and $\gamma_8$ computed by the full order model and by the PC
expansion. Also in the fractures, we observe a high quality for the solutions
computed with the PC expansion even for $\gamma_8$ that has two intersections with other fractures. The jump across the intersection is properly captured, confirming also in this case that the surrogate model from the PC expansion yields a good approximation of the solutions.

\subsubsection{Analysis of variance and correlations}\label{subsec:case2_var}

\begin{figure}[tbp]
    \centering
    \includegraphics[width=1\textwidth]{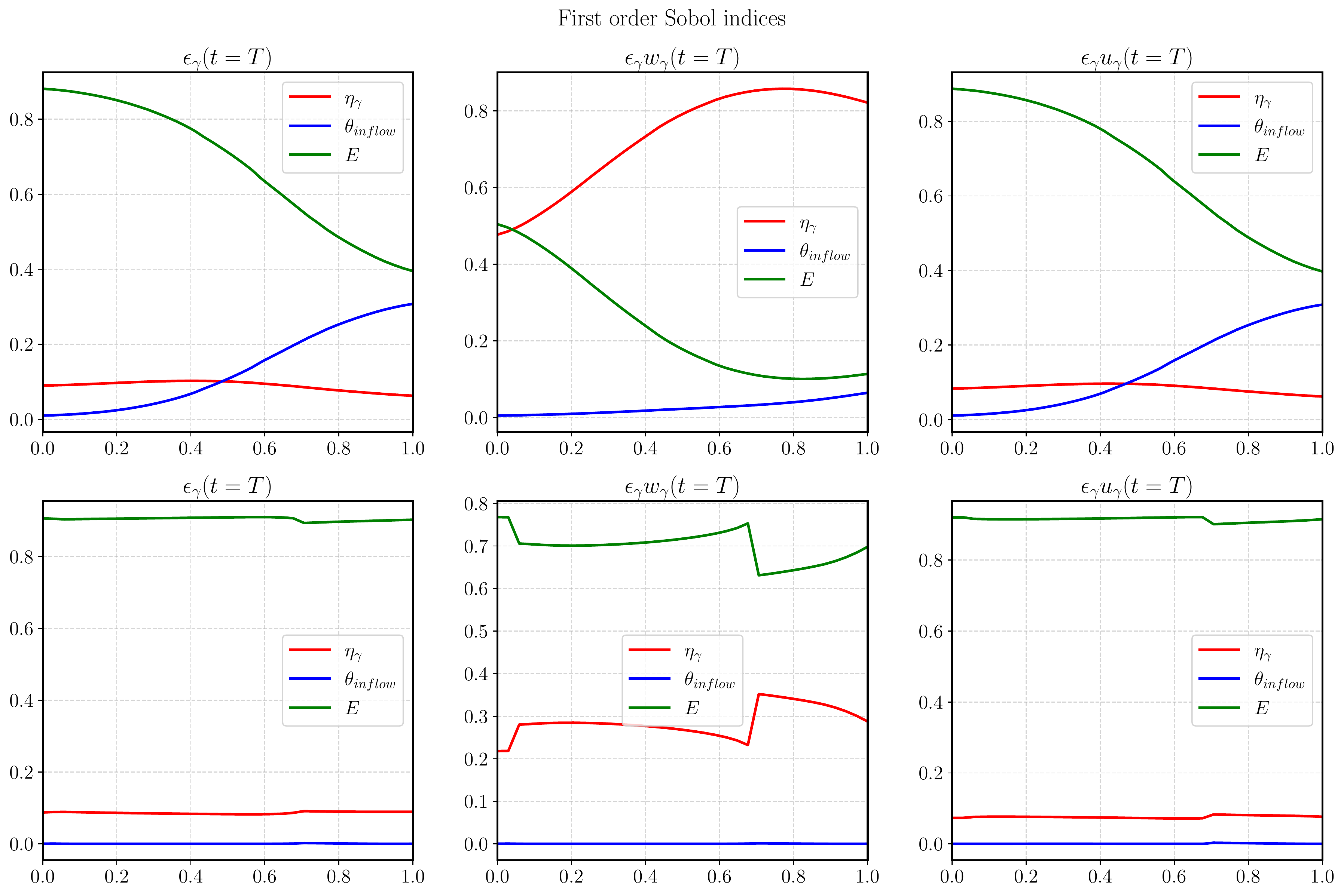}
    \caption{First order Sobol index for different solutions along two
    $\gamma_i$.
    On the top $\gamma_3$ and on the bottom
    for $\gamma_8$. Test case of
    Subsection \ref{subsec:case2_var}.}
    \label{fig:case2_sobol}
\end{figure}

In this section, we present and analyze the impact of the uncertainty on some of the computed
variables. In particular, Figure \ref{fig:case2_sobol} shows the Sobol indices
for some of the variables of interest in the fractures $\gamma_3$ and $\gamma_8$. Since
$\gamma_3$ is closer to the inflow boundary than $\gamma_8$, the associated Sobol indices behave similarly with respect to the ones of the previous test case.
The effect of the inflow is more evident at the lowest tip of the fracture $\gamma_3$ with increased  impact of $\theta^\mathrm{inflow}$ for $\epsilon_\gamma$ and $\epsilon_\gamma u_\gamma$ compared to
the other variables. 
For $\epsilon_\gamma w_\gamma$ the relation expressed by the Sobol index is less clear. 
Since fracture $\gamma_8$ is more distant form the inflow boundary, the high temperature front at the end of the simulation does not fully reach it. Some effect are
still visible since warmer water has been transported by the fractures, especially for $\epsilon_\gamma w_\gamma$ where the activation energy $E$ becomes less important than $\eta_\gamma$ compared with the other two variables under investigation.

\begin{figure}[tbp]
    \centering
    \includegraphics[width=0.32\textwidth]{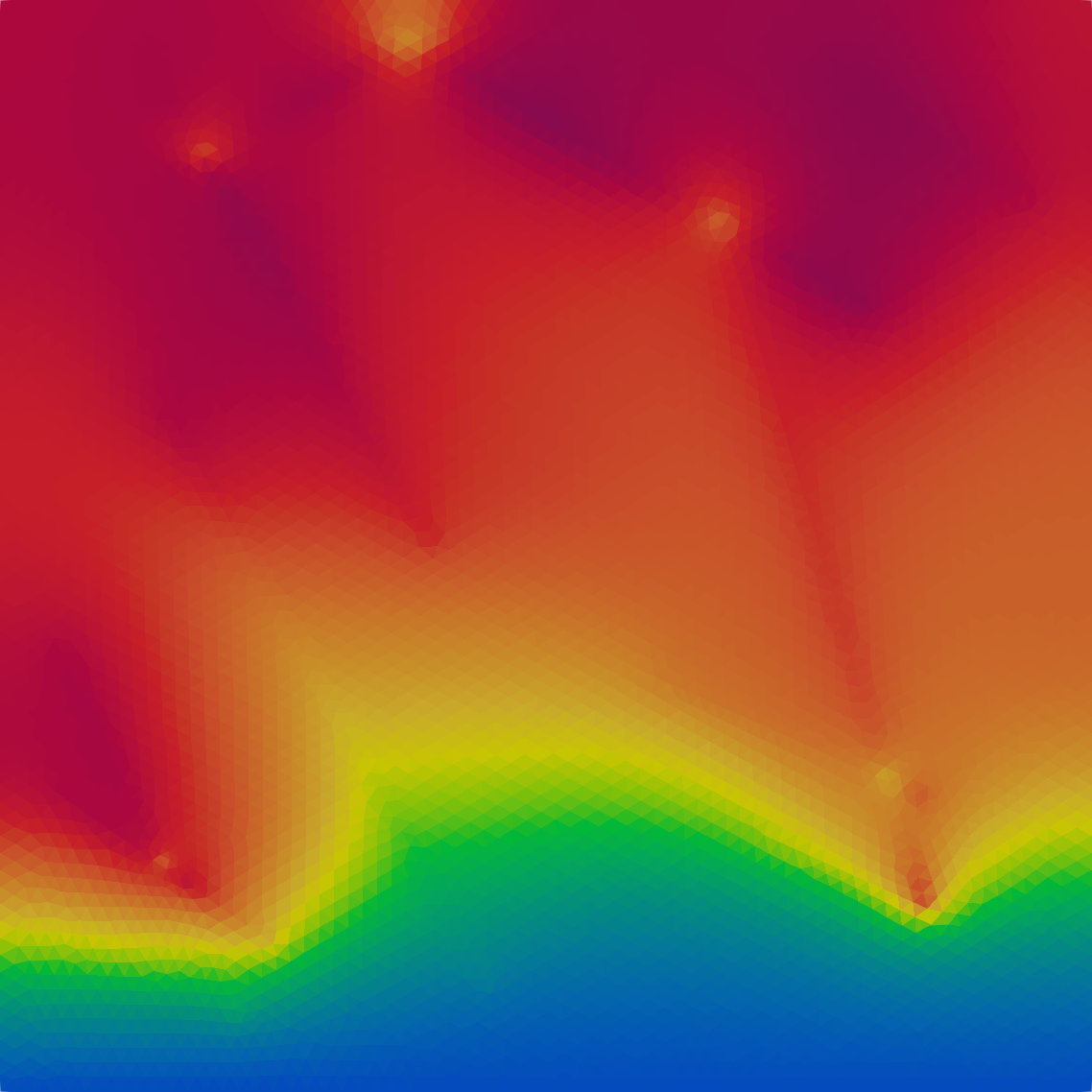}%
    \hspace*{0.02\textwidth}%
    \includegraphics[width=0.32\textwidth]{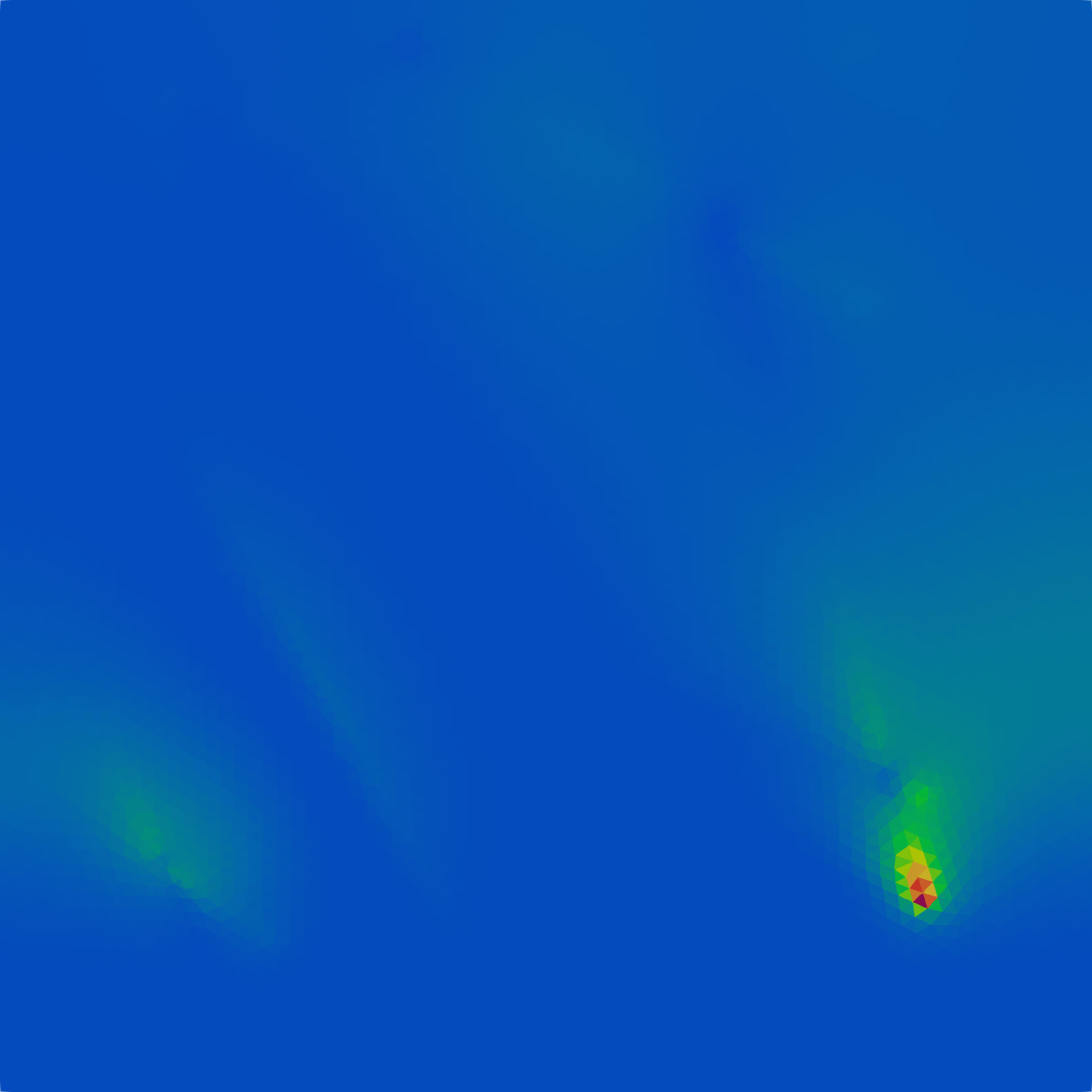}%
    \hspace*{0.02\textwidth}%
    \includegraphics[width=0.32\textwidth]{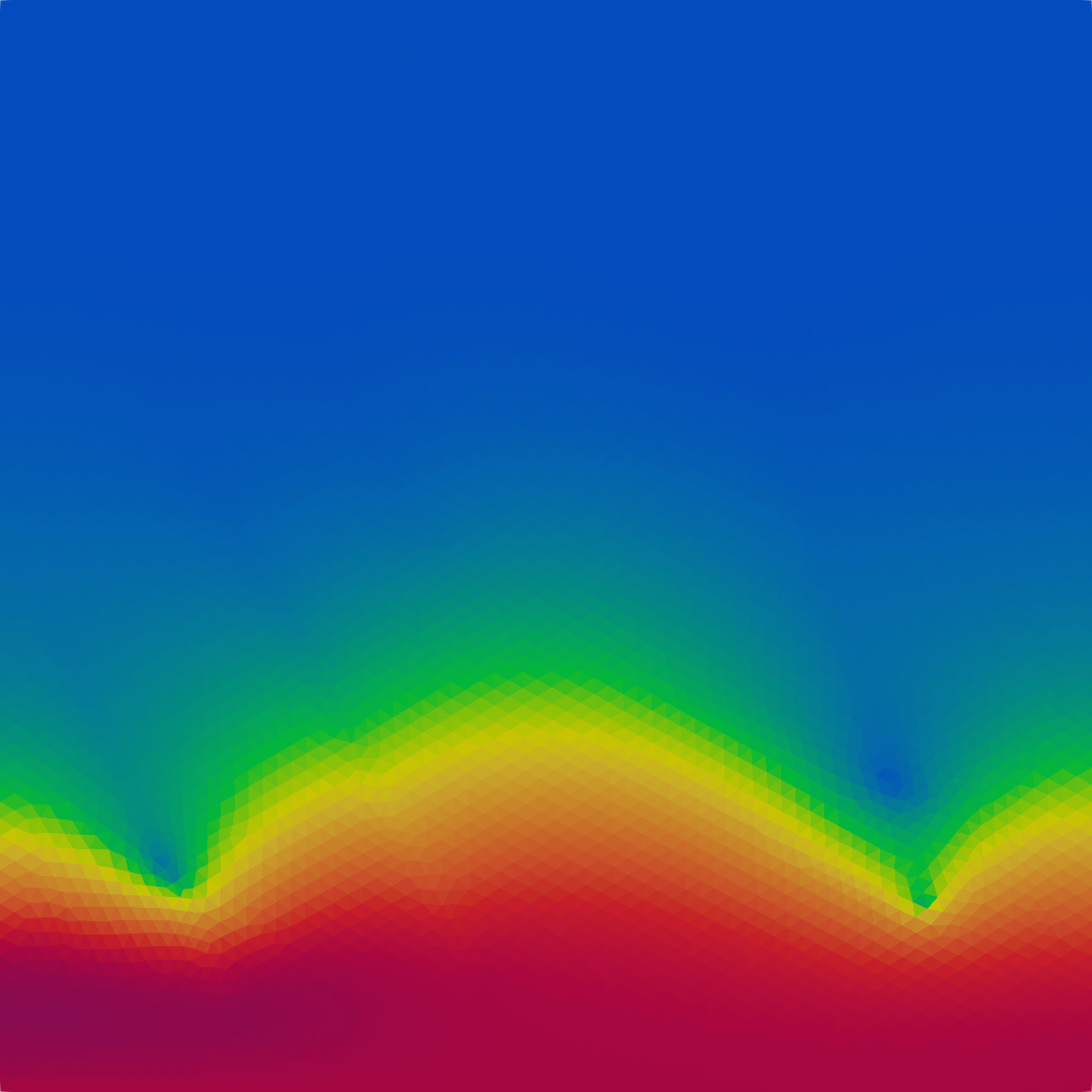}
    \includegraphics[width=0.32\textwidth]{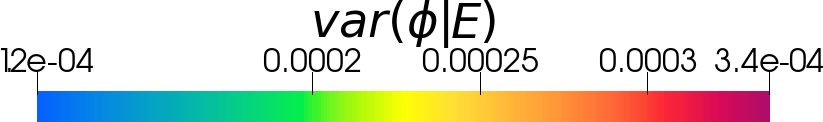}%
    \hspace*{0.02\textwidth}%
    \includegraphics[width=0.32\textwidth]{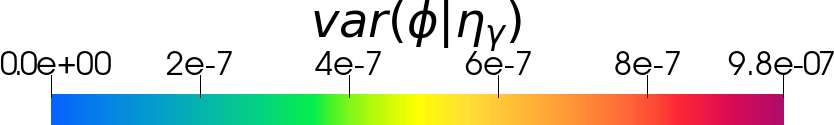}%
    \hspace*{0.02\textwidth}%
    \includegraphics[width=0.32\textwidth]{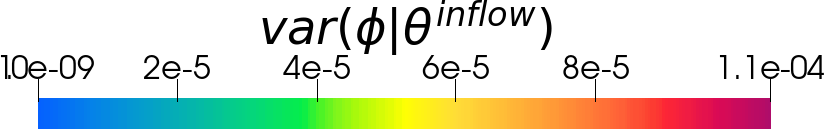}
    \caption{On the left the variance of porosity conditioned to the activation
    energy $E$, on the centre conditioned with $\eta_\gamma$ expansion, and on
    the right conditioned with the temperature inflow $\theta^{inflow}$. Test case of
    Subsection \ref{subsec:case2_var}.}
    \label{fig:case2_var}
\end{figure}

In Figure \ref{fig:case2_var} we present the partial variances of the porosity with respect to the uncertain data at final simulation time. We notice a small impact of $\eta_\gamma$, and a much more significant relevance of $E$ and $\theta^{inflow}$. 
Note that the effect of these two uncertain parameters on reaction speed is opposite. Moreover, we notice that, on the bottom of the domain, the effect of the temperature is more pronounced since it is close to the inflow boundary, 
while for the activation energy $E$ the impact is more predominant away from the inflow and close to the fractures. This can be motivated by the fact that the water and solute get more channelized into the fractures and transported upward. 
Since the fractures do not touch the outflow boundary, the solute flows again into the rock matrix and then alters the value of the porosity by creating more precipitate.

\begin{figure}[tbp]
    \centering
    \includegraphics[width=0.32\textwidth]{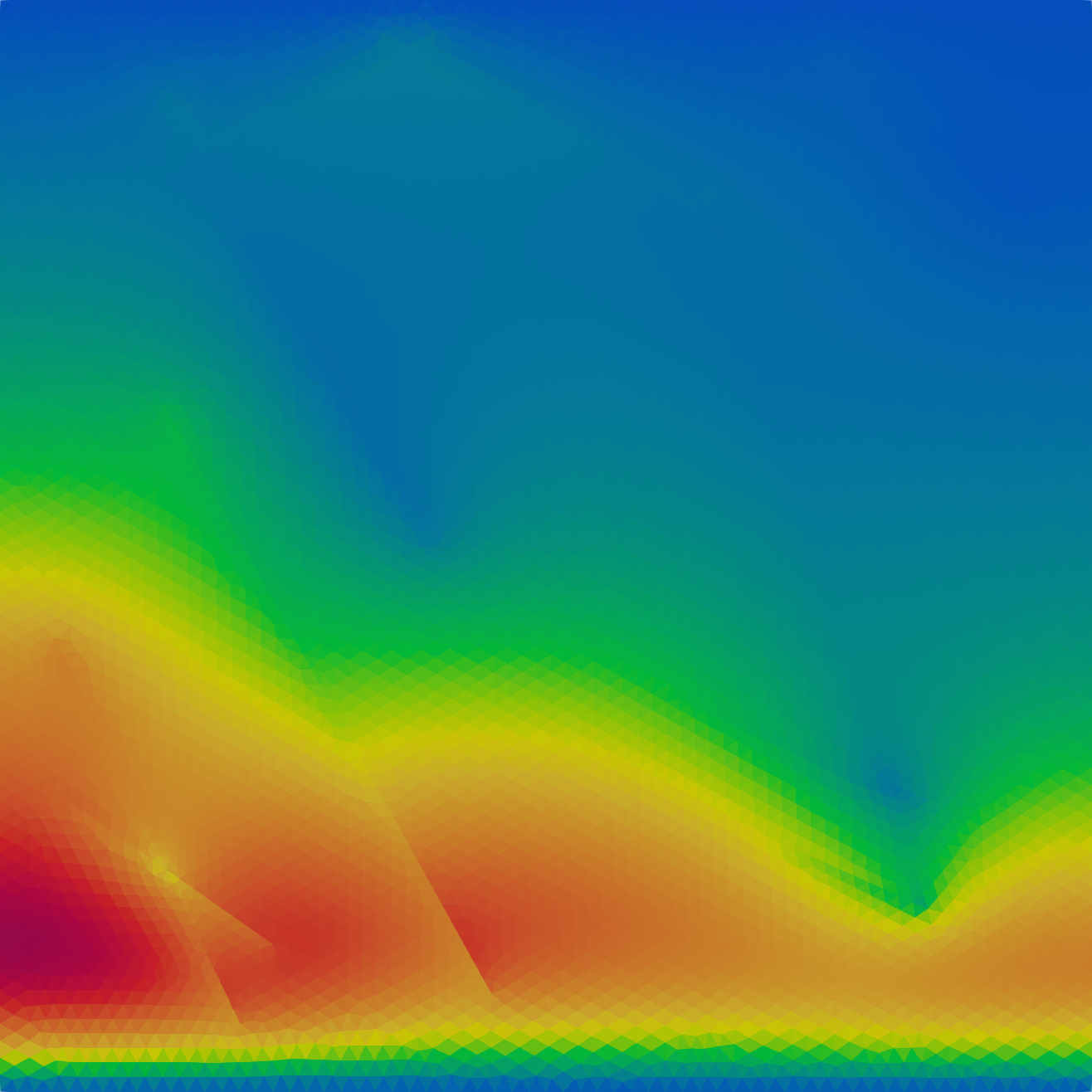}%
    \hspace*{0.02\textwidth}%
    \includegraphics[width=0.32\textwidth]{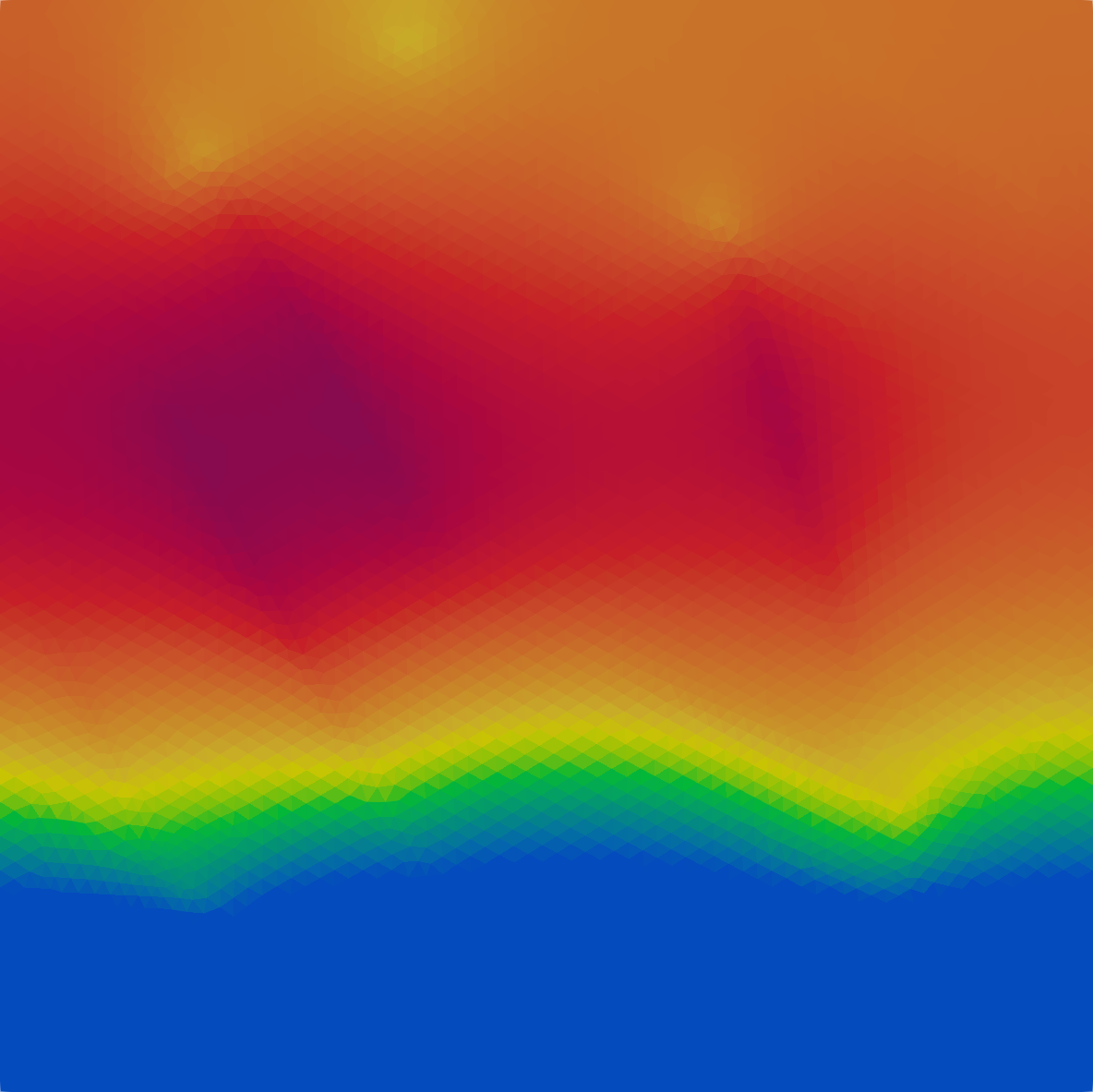}
    \includegraphics[width=0.32\textwidth]{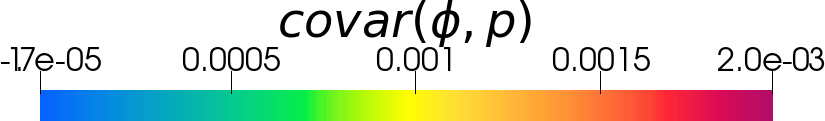}%
    \hspace*{0.02\textwidth}%
    \includegraphics[width=0.32\textwidth]{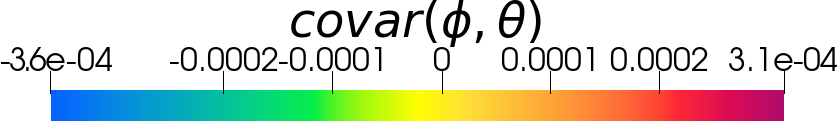}
    \caption{Covariances between different solutions in the porous media. Test case of
    Subsection \ref{subsec:case2_var}.}
    \label{fig:case2_covar}
\end{figure}

The covariances between some of the computed variables are reported in Figure \ref{fig:case2_covar}. 
The correlation between $\phi_\Omega$ and $\theta_\Omega$ is expected: below the warm water front, the increased temperature facilitates the
chemical reaction and thus lowers the porosity; above the front, the solute is lower
than the equilibrium value and the temperature is lower, hence precipitation may not occur once fractures have been sealed and stopped supplying reactant to the upper part of the domain.
The correlation between $\phi_\Omega u_\Omega$ and $\phi_\Omega w_\Omega$ is not included here since they are, as expected, equal to -1 in the whole domain. 
The same applies for $\epsilon_\gamma u_\gamma$ and $\epsilon_\gamma w_\gamma$.

\subsubsection{Probability density functions}\label{subsec:case2_pdf}

Finally, we present the PDFs of some of the
variables at 25\% and 75\% of fractures $\gamma_3$ and $\gamma_8$ for level 2 of the sparse grid. We notice a phenomenon similar to the one observed in the
previous test case. For $\gamma_3$ we might obtain unphysical values of
$\epsilon_\gamma u_\gamma$ when the PDF is
computed by the PC expansion. This issue does not appear for $\gamma_8$, since it is further from the inflow and, as a consequence, is less subject to the uncertainty.
In both cases the probability distribution functions computed with the original model and with the PC expansion are quite similar.
\begin{figure}[tbp]
    \centering
    \includegraphics[width=1\textwidth]{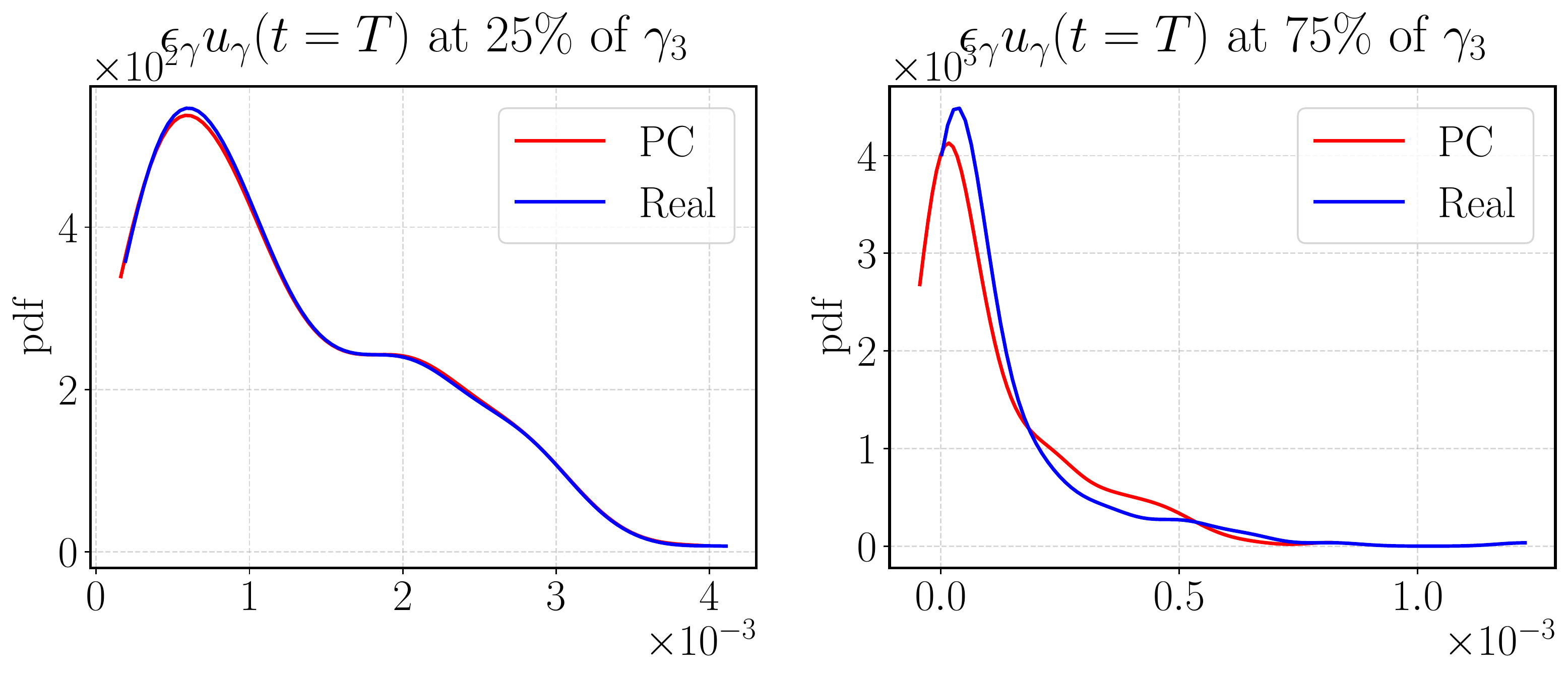}
    \includegraphics[width=1\textwidth]{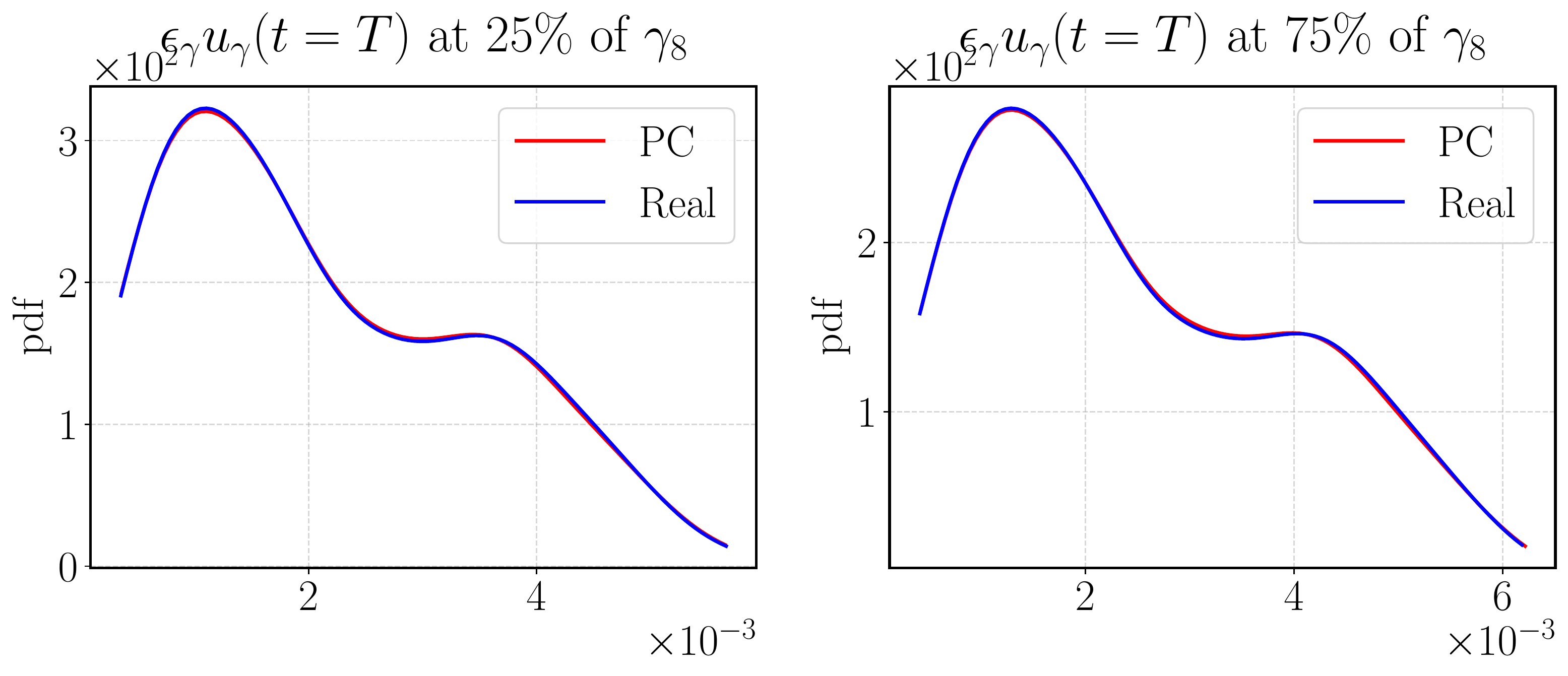}
    \caption{Probability distribution function of $\epsilon_{\gamma_i}
    u_{\gamma_i}$ for
    level 2 at two points in $\gamma_i$. On the top $\gamma_3$ and on the bottom
    for
    $\gamma_8$. Level considered 2. Test case of
    Subsection \ref{subsec:case2_pdf}.}
    \label{fig:case2_pdf_lev2}
\end{figure}

\section{Conclusions}

In this work we have presented a mathematical model to describe the evolution in
time of reactive transport in fractured porous media. We have adopted a simplified model, where only one solute and one precipitate are involved in the
chemical reactions. In several meaningful applications, the chemical processes can be triggered and accelerated by high temperatures, thus we have included in the model also an additional equation to model the thermal effects. 
The solute, transported by a liquid, might form precipitate and alter the porosity in the media, forming a fully coupled and non-linear system. 
In order to obtain a more realistic approximation of the physical process, we have considered a mixed-dimensional model for the fractures where the latter are represented as object of lower dimension. 
We have introduced appropriate equations and derived coupling conditions with the surrounding porous media. Additionally, we have applied a splitting strategy to numerically solve the problem by using standard
discretization schemes.

In a real scenario, several parameters might be affected by uncertainty. 
To quantify their effect on the system, we have considered a polynomial chaos expansion constructed by resorting to spectral projection methods on sparse grids. 
This strategy proved to be very effective, since it provides high quality approximations at low computational cost. Indeed, a limited amount of
simulations are needed to construct a surrogate model which can then be used to perform multiple-simulations. Moreover, from the PC expansion one can easily compute useful statistical quantities such as partial variances and Sobol indexes to investigate the impact of input parameters on the model unknowns and gain insight into the complex model couplings.
This technique has been applied to two numerical examples by increasing the geometrical complexity of the fracture network. 
The results obtained showed the validity of the proposed approach.

\section*{Fundings}

Michele Botti acknowledges funding from the European Union’s Horizon 2020 research and innovation programme under the Marie Skłodowska-Curie grant agreement No. 896616 ``PDGeoFF''. All authors are members of the INdAM Research group GNCS. 

\section*{Conflict of interest}

The authors have no conflicts of interest to declare that are relevant to the content of this article.

\section*{Authors' contribution}

Michele Botti has provided the uncertainty quantification tools used in the work; Alessio Fumagalli and Anna Scotti have devised the mathematical models and all authors contributed to the critical discussion of the test cases, as well as to the writing of the manuscript. Moreover, they all have read and approved the final version. 

\section*{Availability of data and material}

The source code at the basis of the numerical tests is available at \url{https://github.com/pmgbergen/porepy}.

\bibliographystyle{plain}
\bibliography{biblio}

\end{document}